\documentclass[twoside,11pt]{article}

\usepackage{jmlr2e}

\usepackage{appendix}
\usepackage{amsmath}
\usepackage{enumitem}
\usepackage{amsfonts}
\usepackage{multirow}
\usepackage{multicol}

\usepackage{color}
\usepackage{algorithm}
\usepackage{algorithmic}

\usepackage{xcolor,colortbl}
\definecolor{LightCyan}{rgb}{0.88,1,1}

\usepackage{graphicx}
\usepackage{subfig}
\usepackage{hyperref}

\newtheorem{assumption}{Assumption}

\usepackage{lastpage}
\jmlrheading{22}{2021}{1-\pageref{LastPage}}{8/20; Revised
11/21}{12/21}{20-924}{Feihu Huang, Shangqian Gao, Jian Pei, Heng Huang}
\ShortHeadings{Accelerated Zeroth-Order and First-Order Momentum Methods }{Huang,Gao,Pei,Huang}

\begin{document}

\title{ Accelerated Zeroth-Order and First-Order Momentum Methods from Mini to Minimax Optimization }

\author{\name Feihu Huang \email huangfeihu2018@gmail.com \\
       \addr Department of Electrical and Computer Engineering \\
       University of Pittsburgh,
       Pittsburgh, USA
       \AND
       \name Shangqian Gao \email shg84@pitt.edu \\
       \addr Department of Electrical and Computer Engineering \\
       University of Pittsburgh,
       Pittsburgh, USA
       \AND
       \name Jian Pei \email jpei@cs.sfu.ca \\
       \addr School of Computing Science \\
        Simon Fraser University,
        Vancouver, Canada
       \AND
       \name Heng Huang \email heng.huang@pitt.edu \\
       \addr Department of Electrical and Computer Engineering \\
        University of Pittsburgh,
        Pittsburgh, USA
       }

\editor{Tong Zhang}

\maketitle

\begin{abstract}
In the paper, we propose a class of accelerated zeroth-order and first-order momentum methods for both  nonconvex mini-optimization and minimax-optimization.
Specifically, we propose a new accelerated zeroth-order momentum (Acc-ZOM) method for black-box mini-optimization where only function values can be obtained. 
Moreover, we prove that our Acc-ZOM method achieves a lower query complexity of $\tilde{O}(d^{3/4}\epsilon^{-3})$ for finding an $\epsilon$-stationary point, which improves the best known result by a factor of $O(d^{1/4})$ where $d$ denotes the variable dimension.
In particular, our  Acc-ZOM does not need large batches required in the existing zeroth-order stochastic algorithms.
Meanwhile, we propose an accelerated zeroth-order momentum descent ascent (Acc-ZOMDA) method for black-box minimax  optimization, where only function values can be obtained. Our Acc-ZOMDA obtains a low query complexity of $\tilde{O}((d_1+d_2)^{3/4}\kappa_y^{4.5}\epsilon^{-3})$ without requiring large batches 
for finding an $\epsilon$-stationary point, where $d_1$ and $d_2$ denote variable dimensions 
and $\kappa_y$ is condition number. 
Moreover, we propose an accelerated first-order momentum descent ascent (Acc-MDA) method for minimax optimization,  whose explicit gradients are accessible. Our Acc-MDA achieves a low  gradient complexity of $\tilde{O}(\kappa_y^{4.5}\epsilon^{-3})$ without requiring large batches for finding an $\epsilon$-stationary point. In particular, our Acc-MDA can 
obtain a lower gradient complexity of $\tilde{O}(\kappa_y^{2.5}\epsilon^{-3})$ with a batch size $O(\kappa_y^4)$, which improves the best known result by a factor of $O(\kappa_y^{1/2})$.  
Extensive experimental results on black-box adversarial attack to deep neural networks and
poisoning attack to logistic regression demonstrate efficiency of our algorithms.
\end{abstract}


\section{Introduction}
In the paper, we consider solving the following stochastic mini-optimization problem:
\begin{align} \label{eq:1}
 \min_{x \in \mathcal{X}} f(x)=\mathbb{E}_{\xi\sim \mathcal{D}}[f(x;\xi)],
\end{align}
where $f(x): \mathcal{X} \rightarrow \mathbb{R}$ is a differentiable and possibly nonconvex function, and $\mathcal{X}\subseteq \mathbb{R}^d$ is a convex closed set,
and $\xi$ is a random variable following an unknown distribution $\mathcal{D}$.
In machine learning, the expected loss minimization is generally expressed as the problem \eqref{eq:1}.
Stochastic Gradient Descent (SGD) is a standard algorithm for solving the problem \eqref{eq:1}.
However, it suffers from large variances resulting in a high gradient complexity of $O(\epsilon^{-4})$ \citep{ghadimi2013stochastic}
for finding an $\epsilon$-stationary point, \emph{i.e.}, $\mathbb{E}\|\nabla f(x)\|\leq \epsilon$.
Thus, many variance-reduced algorithms \citep{allen2016variance,reddi2016stochastic,zhou2018stochastic,fang2018spider,wang2019spiderboost}
have been developed to improve the gradient complexity of the SGD.
Specifically, \cite{allen2016variance,reddi2016stochastic} proposed the nonconvex version of SVRG algorithm \citep{johnson2013accelerating},
which reaches an improved gradient complexity of $O(\epsilon^{-10/3})$. Subseqently, the SNVRG/SPIDER methods \citep{zhou2018stochastic,fang2018spider,wang2019spiderboost}
have been proposed to obtain a near-optimal gradient complexity of $O(\epsilon^{-3})$.
More recently, the momentum-based variance reduced methods \citep{cutkosky2019momentum,tran2019hybrid} achieved the best known complexity of $\tilde{O}(\epsilon^{-3})$.
At the same time, \cite{arjevani2019lower} established a lower bound of complexity $O(\epsilon^{-3})$ for variance reduced algorithms.

The above first-order methods need to use gradients of the objective function to update the variables.
In many machine learning problems, however, the explicit gradients of their objective functions are difficult or infeasible to access. For example, in the reinforcement learning \citep{malik2020derivative,kumar2020zeroth,huang2020momentum},
it is difficult to calculate the explicit gradients of their objective functions.
Even worse, in the black-box adversarial attack to deep neural networks (DNNs) \citep{chen2018frank},
only prediction labels can be obtained.
To solve such back-box problem \eqref{eq:1} where only the objective function values can be obtained, 
the zeroth-order methods \citep{ghadimi2013stochastic,duchi2015optimal} have been widely used with only querying values of the function $f(x)$
and not accessing to its explicit formation.
Recently, some zeroth-order stochastic algorithms \citep{ghadimi2013stochastic,duchi2015optimal,Nesterov2017RandomGM,chen2019zo}
have been presented by using the smoothing techniques such as Gaussian-distribution and Uniform-distribution smoothing.
Similarly, these zeroth-order stochastic algorithms also suffer from large variances resulting in a high query
complexity of $O(d\epsilon^{-4})$ \citep{ghadimi2013stochastic} for finding an $\epsilon$-stationary point.
To reduce the  query complexity, \citet{fang2018spider,ji2019improved} recently proposed some 
accelerated zeroth-order stochastic  algorithms (\emph{i.e.}, SPIDER-SZO and ZO-SPIDER-Coord)
based on the variance reduced technique of SPIDER \citep{fang2018spider}.
Although these accelerated zeroth-order  methods obtain a lower query complexity of $O(d\epsilon^{-3})$,
these methods require large batches in both inner and outer loops of algorithms.
At the same time, the practical performances of these methods are
not consistent with this low query complexity,
since they  require large batches and
strict learning rates to achieve it.

In the paper, thus, we propose a new accelerated zeroth-order momentum (Acc-ZOM) method to
solve the black-box problem \eqref{eq:1}, which builds on both generic uniform smoothing gradient estimator and
 momentum-based variance reduction technique of STORM/Hybrid-SGD \citep{cutkosky2019momentum,tran2019hybrid}. 
Moreover, we prove that our  Acc-ZOM method achieves a lower function query complexity of $O(d^{3/4}\epsilon^{-3})$ without large batches for finding an $\epsilon$-stationary point, which improves the best known complexity by a factor of $O(d^{1/4})$
(please see Table \ref{tab:1} for query complexity comparison of different non-convex zeroth-order methods).

\begin{table}
  \centering
  \caption{ \textbf{Query complexity} comparison of the representative non-convex zeroth-order methods for finding
  an $\epsilon$-stationary point of the \textbf{black-box} mini-optimization problem \eqref{eq:1} and minimax-optimization problem \eqref{eq:2}, respectively.
  GauGE, UniGE and CooGE are abbreviations of Gaussian, Uniform and Coordinate-Wise smoothing gradient estimators, respectively.
  Here $\kappa_y$ denotes the condition number for function $f(\cdot,y)$. Note that Appendix \ref{appendix:B} provides a comparison of assumptions used in the zeroth-order methods, and  
  Appendix \ref{appendix:C} provides a detailed proof to obtain a correct query complexity of ZO-Min-Max algorithm \citep{liu2019min}.
   } \label{tab:1}
  \resizebox{\textwidth}{!}{
\begin{tabular}{c|c|c|c|c|c}
  \hline
  \textbf{Problem} & \textbf{Algorithm} & \textbf{Reference}  & \textbf{Estimator} & \textbf{Batch Size} & \textbf{Complexity} \\ \hline
  \multirow{6}*{Mini} & ZO-SGD & \citet{ghadimi2013stochastic}   & GauGE & $O(1)$ & $O(d\epsilon^{-4})$  \\ \cline{2-6}
  &ZO-AdaMM & \citet{chen2019zo}  & UniGE & $O(\epsilon^{-2})$ & $O(d^2\epsilon^{-4})$ \\ \cline{2-6}
  &ZO-SVRG &  \citet{ji2019improved}  & CooGE & $O(\epsilon^{-2})$ &  $O(d\epsilon^{-10/3})$  \\ \cline{2-6}
  &ZO-SPIDER-Coord & \citet{ji2019improved}  & CooGE & $O(\epsilon^{-2})$ &  $O(d\epsilon^{-3})$  \\ \cline{2-6}
  &SPIDER-SZO & \citet{fang2018spider}  & CooGE & $O(\epsilon^{-2})$ &  $O(d\epsilon^{-3})$ \\ \cline{2-6}
  \rowcolor{LightCyan} &Acc-ZOM & Ours  & UniGE & {\color{red}{ $O(1)$}} & {\color{red}{ $O(d^{3/4}\epsilon^{-3})$ }}\\ \hline
  \multirow{7}*{Minimax} & ZO-Min-Max &\citet{liu2019min} & UniGE & $O((d_1\!+\!d_2)\kappa_y^2\epsilon^{-2})$ & $O((d_1\!+\!d_2)\kappa_y^6\epsilon^{-6})$ \\ \cline{2-6}
  &ZO-SGDA  & \citet{wang2020zeroth} & GauGE & $O((d_1\!+\!d_2)\epsilon^{-2})$ & $O((d_1\!+\!d_2)\kappa_y^5\epsilon^{-4})$ \\  \cline{2-6}
  &ZO-SGDMSA  & \citet{wang2020zeroth} & GauGE & $O((d_1\!+\!d_2)\epsilon^{-2})$ & $\tilde{O}((d_1\!+\!d_2)\kappa_y^2\epsilon^{-4})$ \\  \cline{2-6}
  &ZO-SREDA-Boost  & \citet{xu2020enhanced} & CooGE & $O(\max(\kappa_y\epsilon^{-1},d_1+d_2)\kappa_y\epsilon^{-1})$ & $O((d_1\!+\!d_2)\kappa_y^3\epsilon^{-3})$ \\  \cline{2-6}
  \rowcolor{LightCyan} & Acc-ZOMDA & Ours & UniGE & {\color{red}{ $O(1)$}} & {\color{red}{ $\tilde{O}((d_1\!+\!d_2)^{3/4}\kappa_y^{4.5}\epsilon^{-3})$ }} \\
  \hline
\end{tabular}
}
\end{table}

Besides the mini-optimization  problem  \eqref{eq:1} is  widely used in machine learning,
there also exist many machine learning applications \citep{shapiro2002minimax,nouiehed2019solving,zhao2020primal}
such as adversarial training \citep{goodfellow2014generative}, reinforcement learning \citep{wai2019variance,wai2018multi}, distributionally robust optimization \citep{qi2020practical} and AUC maximization \citep{ying2016stochastic},
which can be modeled as a minimax optimization problem. In the paper, we further focus on solving the following stochastic minimax optimization problem:
\begin{align} \label{eq:2}
 \min_{x \in \mathcal{X}} \max_{y\in \mathcal{Y}} f(x,y)=\mathbb{E}_{\xi\sim \mathcal{D}'}[f(x,y;\xi)],
\end{align}
where function $f(x,y): \mathcal{X}\times \mathcal{Y} \rightarrow \mathbb{R}$ is strongly concave in variable $y$ but possibly nonconvex in variable $x$, and $\xi$ is a random variable 
following an unknown distribution $\mathcal{D}'$. Here the constraint sets $\mathcal{X}\subseteq \mathbb{R}^{d_1}$
and $ \mathcal{Y}\subseteq \mathbb{R}^{d_2}$
are compact and convex. In fact, the problem \eqref{eq:2} can be seen as a zero-sum game between two players.
The goal of the first player is to minimize $f(x,y)$ by varying $x$, while the other player's aim is to maximize
$f(x,y)$ by varying $y$.
When the problem \eqref{eq:2} is black-box where only noise stochastic function values can be obtained,
we propose an accelerated zeroth-order momentum descent ascent (Acc-ZOMDA) method based on 
the generic uniform smoothing gradient estimator and the variance reduced technique of STORM.
When the problem \eqref{eq:2} is transparent where noise stochastic gradients can be accessed, 
we present an accelerated first-order momentum descent ascent (Acc-MDA) method based on the variance reduced technique of STORM. 

\noindent\textbf{Contributions:}
Our main contributions are summarized as follows:
\begin{itemize}
\setlength{\itemsep}{0pt}
\item[1)] We propose a new accelerated zeroth-order momentum (Acc-ZOM) method to solve the \textbf{black-box mini-optimization} problem \eqref{eq:1}, where only noise stochastic function values can be obtained.
    Moreover, we prove that our Acc-ZOM method achieves a lower query complexity of $O(d^{3/4}\epsilon^{-3})$ for finding an $\epsilon$-stationary point  without requiring large batches, which improves the best known result by a factor of $O(d^{1/4})$.
\item[2)] We propose an accelerated zeroth-order momentum descent ascent (Acc-ZOMDA) method to solve the \textbf{black-box minimax-optimization} problem \eqref{eq:2}, where only noise stochastic function values can be obtained. Moreover, we prove that our Acc-ZOMDA method obtains a low query complexity of $O\big((d_1+d_2)^{3/4}\kappa_y^{4.5}\epsilon^{-3}\big)$ without requiring large batches for finding an $\epsilon$-stationary point (Please see Table \ref{tab:1}).
\item[3)] We further present propose an accelerated first-order momentum descent ascent (Acc-MDA) method to solve the \textbf{transparent minimax-optimization} problem \eqref{eq:2}, whose explicit gradients are accessible. We prove that our Acc-MDA algorithm has a low gradient complexity of $\tilde{O}\big(\kappa_y^{4.5}\epsilon^{-3}\big)$ without requiring large batches for finding an $\epsilon$-stationary point. Our Acc-MDA algorithm reaches the best known gradient complexity of $\tilde{O}\big(\kappa_y^{3}\epsilon^{-3}\big)$ with batch size $O(\kappa_y^3)$ for finding an $\epsilon$-stationary point. Moreover,
our Acc-MDA algorithm obtains a lower gradient complexity of $\tilde{O}\big(\kappa_y^{2.5}\epsilon^{-3}\big)$ with batch size $O(\kappa_y^4)$ for finding an $\epsilon$-stationary point (Please see Table \ref{tab:2}).
\item[4)] We present a class of accelerated zeroth-order and first-order momentum framework for both mini-optimization and minimax-optimization. Moreover, we study the convergence properties of our methods for both \textbf{constrained} and \textbf{unconstrained} optimization, respectively.
\end{itemize}

\begin{table}
  \centering
  \caption{ \textbf{Gradient complexity} comparison of the representative first-order methods for finding
  an $\epsilon$-stationary point of the minimax problem \eqref{eq:2}.
  Here Y denotes the fact that there exists a convex constraint on variable, otherwise is N. Note that our theoretical results
  do not rely on any assumption on convex constraint sets $\mathcal{X}$ and $\mathcal{Y}$,
  so it can be easily extend to the unconstrained setting.
   }
  \label{tab:2}
  \resizebox{\textwidth}{!}{
\begin{tabular}{c|c|c|c|c|c}
  \hline
   \textbf{Algorithm} & \textbf{Reference}  & \textbf{Constraint} on $x,y$ & \textbf{Loop(s)} & \textbf{Batch Size} & \textbf{Complexity} \\ \hline
  PGSVRG & \citet{rafique2018non} & N,\ N & Double & $O(\epsilon^{-2})$ & $O(\kappa_y^3\epsilon^{-4})$ \\  \hline
  SGDA  & \citet{lin2019gradient} & N,\ Y & Single & $O(\kappa_y\epsilon^{-2})$ & $O(\kappa_y^3\epsilon^{-4})$ \\  \hline
  SREDA  & \citet{luo2020stochastic} & N,\ Y & Double& $O(\kappa_y^2\epsilon^{-2})$ & $O(\kappa_y^3\epsilon^{-3})$ \\  \hline
  SREDA-Boost  & \citet{xu2020enhanced} & N,\ N & Double & $O(\kappa_y^2\epsilon^{-2})$ & $O(\kappa_y^3\epsilon^{-3})$ \\  \hline
  \rowcolor{LightCyan} Acc-MDA & Ours & Y (N),\ Y & Single & {\color{red}{ $O(1)$}} & {\color{red}{ $\tilde{O}(\kappa_y^{4.5}\epsilon^{-3})$ }} \\ \hline
  \rowcolor{LightCyan} Acc-MDA & Ours & Y (N),\ Y & Single & {\color{red}{ $O(\kappa_y^{\nu}), \ \nu>0$}} & {\color{red}{ $\tilde{O}(\kappa_y^{(4.5-\nu/2)}\epsilon^{-3})$ }} \\
  \hline
\end{tabular}
}
\end{table}

The remainder of the paper is structured as follows. In Section \ref{section:2}, we review some related works about zeroth-order and first-order methods for mini and minimax optimization. Section \ref{section:3} introduces some preliminaries
 about zeroth-order and first-order methods for mini and minimax optimization. 
We introduce our Acc-ZOM, Acc-ZOMDA and Acc-MDA methods in Sections \ref{section:4}, \ref{section:5} and \ref{section:6}, respectively. In Section \ref{section:7}, we give the  convergence properties of our methods. In Section \ref{section:8}, we apply black-box adversarial attack to DNNs and poisoning attack to logistic regression to verify efficiency of our methods. Conclusions are provided in Section \ref{section:9}. The proofs
of the main results are given in the appendix. 
\section{Related Works}
\label{section:2}
In this section, we recap some zeroth-order and first-order methods for solving the mini-optimization and
minimax-optimization problems, respectively.

\subsection{ Zeroth-Order Mini-Optimization }
Zeroth-order (i.e., gradient-free) methods are a class of powerful optimization tools to solve many complex  machine learning problems,
whose explicit gradients are difficult or even infeasible to access.
Recently, the zeroth-order methods have been widely proposed. 
For example, \cite{ghadimi2013stochastic,duchi2015optimal,Nesterov2017RandomGM} proposed several zeroth-order algorithms
based on the Gaussian smoothing technique.
Subsequently, some accelerated zeroth-order stochastic methods \citep{liu2018zeroth,ji2019improved} have been proposed
by using the variance reduced techniques.
To solve the constrained optimization, the zeroth-order projected method \citep{liu2018projected} and
the zeroth-order Frank-Wolfe methods \citep{balasubramanian2018zeroth,chen2018frank,sahu2019towards,huang2020accelerated}
have been recently proposed.
More recently, \cite{chen2019zo} have  proposed a zeroth-order adaptive momentum method to solve the constrained optimization problems.
To solve the nonsmooth optimization, several zeroth-order proximal algorithms \citep{ghadimi2016mini,huang2019faster,ji2019improved}
and zeroth-order ADMM-based algorithms \citep{gao2018information,liu2018admm,huang2019zeroth,huang2019nonconvex}
have been proposed.

\subsection{ Zeroth-Order Minimax Optimization }
The above zeroth-order methods only focus on the mini-optimization problems.
In fact, many machine learning problems such as reinforcement learning \citep{wai2019variance,wai2018multi},
black-box adversarial attack \citep{liu2019min}, and adversarial training \citep{goodfellow2014generative,liu2019towards}
can be expressed as the minimax-optimization problems.
For the black-box minimax problems where we can only access function values, more recently, some zeroth-order descent ascent
methods \citep{liu2019min,wang2020zeroth,xu2020enhanced} have been presented
to solve the minimax-optimization problem \eqref{eq:2}. In addition,  online zeroth-order
extra-gradient algorithms \citep{roy2019online}
have been proposed to solve the  (strongly) convex-concave minimax problems.

\subsection{ First-Order Minimax Optimization }
For the transparent minimax problems whose explicit gradients are accessible, more recently, some first-order minimax  methods have been widely studied in  \citep{rafique2018non,jin2019local,nouiehed2019solving,thekumparampil2019efficient,lin2019gradient,yang2020global,ostrovskii2020efficient,
yan2020sharp,lin2020near,xu2020unified,boct2020alternating}. 
For example, \cite{lin2019gradient} proposed a class of gradient descent ascent methods (i.e., GDA and SGDA) for nonconvex-(strongly) concave minimax problems. \cite{rafique2018non} studied a class of weakly-convex concave minimax problems and proposed an efficient stochastic gradient descent ascent method (i.e., PGSVRG) based on the variance reduced technique of SVRG.
\cite{luo2020stochastic,xu2020enhanced} proposed a class of faster SGDA  methods (i.e., SREDA and SREDA-Boost) to solve the nonconvex-strongly-concave minimax problems
based on the variance reduced technique of SARAH/SPIDER.
In addition, \cite{tran2020hybrid} presented a hybrid variance-reduced SGD algorithm
for a special case of nonconvex-concave stochastic minimax problems, which are equivalent to a class of stochastic compositional problems studied in \citep{qi2020practical}.

\section{Preliminaries}
\label{section:3}
In this section, we introduce zeroth-order gradient estimators and some mild assumptions for mini-optimization
problem \eqref{eq:1} and minimax-optimization problem \eqref{eq:2}, respectively. 

\subsection{Notations}
$\langle x,y\rangle$ denotes the inner product of two vectors $x$ and $y$. $\|\cdot\|$ denotes the $\ell_2$ norm
for vectors and spectral norm for matrices. $I_{d}$ denotes a $d$-dimensional identity matrix.
Given function $f(x,y)$, $f(x,\cdot)$ denotes  function \emph{w.r.t.} the second variable with fixing $x$,
and $f(\cdot,y)$ denotes function \emph{w.r.t.} the first variable
with fixing $y$. Let $\nabla f(x,y)=(\nabla_x f(x,y),\nabla_y f(x,y))$, where $\nabla_x f(x,y)$ and $\nabla_y f(x,y)$
denote the partial gradients \emph{w.r.t.} variables $x$ and $y$, respectively. Define two
increasing $\sigma$-algebras $\mathcal{F}^1_t := \{\xi_1,\xi_2,\cdots,\xi_{t-1}\}$
and $\mathcal{F}^2_t := \{u^1,u^2,\cdots,u^{t-1}\}$ for all $t\geq 2$, where $\{u^i\}_{i=1}^{t-1}$ is a vector generated
from the uniform distribution over the unit sphere, then
let $\mathbb{E}[\cdot]=\mathbb{E}[\cdot|\mathcal{F}^1_t,\mathcal{F}^2_t]$.
We denote $a=O(b)$ if $a\leq C b$ for some constant $C>0$. The notation $\tilde{O}(\cdot)$ hides logarithmic terms.
Given a convex closed set $\mathcal{X}$, we define a projection operation to $\mathcal{X}$ as
$\mathcal{P}_{\mathcal{X}}(x_0) = \arg\min_{x\in \mathcal{X}} \|x-x_0\|^2$.

\subsection{ Preliminaries for Mini-Optimization}
For solving the mini-optimization problem \eqref{eq:1},
we apply the \textbf{Uni}form smoothing \textbf{G}radient \textbf{E}stimator (UniGE) \citep{gao2018information,ji2019improved}
to generate stochastic zeroth-order gradients.
Specifically, given the stochastic function $f(x;\xi): \mathbb{R}^d \rightarrow \mathbb{R}$,
the UniGE can generate a stochastic zeroth-order gradient, defined as
\begin{align} \label{eq:3}
 \hat{\nabla} f(x;\xi) = \frac{ f(x + \mu u;\xi)-f(x;\xi)}{\mu/d}u,
\end{align}
where $u \in \mathbb{R}^d$ is a vector generated from the uniform distribution over the unit sphere,
and $\mu$ is a smoothing parameter.
Let $f_{\mu}(x;\xi)=\mathbb{E}_{u\sim U_B}[f(x+\mu u;\xi)]$ be a smooth approximation of $f(x;\xi)$,
where $U_B$ is the uniform distribution over the $d$-dimensional unit Euclidean ball $B$.
Further let $\nabla f_{\mu}(x)=\mathbb{E}_{\xi}[\nabla f_{\mu}(x;\xi)]$.
According to Lemma 5 in \citep{ji2019improved}, 
we have $\mathbb{E}_{(\xi,u)}[\hat{\nabla} f(x;\xi)]=\nabla f_{\mu}(x)$. Next, we give some mild assumptions about the problem \eqref{eq:1}.
\begin{assumption} \label{ass:1}
The variance of stochastic zeroth-order gradient is bounded, \emph{i.e.,} there exists a constant $\sigma >0$ such that for all $x$,
it follows $\mathbb{E}\|\hat{\nabla} f(x;\xi) - \nabla f_{\mu}(x)\|^2 \leq \sigma^2$.
\end{assumption}
Assumption \ref{ass:1} is similar to the upper bound of variance of stochastic gradient in \citep{ghadimi2013stochastic,cutkosky2019momentum}.
In the following, we further give some mild conditions about the problem \eqref{eq:1}.
\begin{assumption} \label{ass:3}
The component function $f(x;\xi)$ is $L$-smooth such that
\begin{align}
 \|\nabla f(x;\xi)-\nabla f(x';\xi)\| \leq L \|x - x'\|, \ \forall x,x' \in \mathcal{X}.  \nonumber
\end{align}
\end{assumption}
\begin{assumption} \label{ass:4}
The function $f(x)$ is bounded from below in $\mathcal{X}$, \emph{i.e.,} $f^* = \inf_{x\in \mathcal{X}} f(x)$.
\end{assumption}
Assumption \ref{ass:3} imposes smoothness on each component loss function, which is widely used
in the nonconvex algorithms \citep{fang2018spider,wang2019spiderboost,cutkosky2019momentum}.
Assumptions \ref{ass:4} guarantees the feasibility of the problem \eqref{eq:1}.

\subsection{ Preliminaries for Minimax-Optimization}
For solving the minimax-optimization problem \eqref{eq:2},
we also apply the UniGE to generate stochastic zeroth-order partial gradients.
Specifically, for the stochastic function $f(x,y;\xi): \mathbb{R}^{d_1}\times \mathbb{R}^{d_2} \rightarrow \mathbb{R}$,
given $\mathcal{B}=\{\xi_1,\cdots,\xi_b\}$ drawn i.i.d. from an unknown distribution, the UniGE can generate stochastic
zeroth-order partial gradients, defined as
\begin{align}
 & \hat{\nabla}_x f(x,y;\mathcal{B}) = \frac{1}{b}\sum_{i=1}^b\hat{\nabla}_x f(x,y;\xi_i) = \frac{1}{b}\sum_{i=1}^b\frac{ f(x + \mu_1 \hat{u}_i, y;\xi_i)-f(x,y;\xi_i)}{\mu_1/d_1}\hat{u}_i, \label{eq:5} \\
 & \hat{\nabla}_y f(x,y;\mathcal{B}) = \frac{1}{b}\sum_{i=1}^b\hat{\nabla}_y f(x,y;\xi_i) = \frac{1}{b}\sum_{i=1}^b\frac{ f(x, y+ \mu_2 \tilde{u}_i;\xi_i)-f(x,y;\xi_i)}{\mu_2/d_2}\tilde{u}_i, \label{eq:6}
\end{align}
where $\mu_1$ and $\mu_2$ are the smoothing parameters, and $\hat{U}=\{\hat{u}_i \in \mathbb{R}^{d_1}\}_{i=1}^b$
and $\tilde{U}=\{\tilde{u}_2 \in \mathbb{R}^{d_2}\}_{i=1}^b$ are
generated from the uniform distribution over the unit sphere $U_{B_1}$ and
$U_{B_2}$, respectively. Here $U_{B_1}$ and
$U_{B_2}$ denote the uniform distributions over the $d_1$-dimensional unit Euclidean ball $B_1$ and
$d_2$-dimensional unit Euclidean ball $B_2$, respectively.
The smoothed functions associated to function $f(x,y;\xi)$ can be defined as:
\begin{align}
 f_{\mu_1}(x,y;\xi) = \mathbb{E}_{\hat{u}}\big[f(x+\mu_1\hat{u},y;\xi)\big], \quad f_{\mu_2}(x,y;\xi) = \mathbb{E}_{\tilde{u}}\big[f(x,y+\mu_2\tilde{u};\xi)\big].
\end{align}
Following Lemma 5 in \citep{ji2019improved},
we have $\mathbb{E}_{(\hat{u},\xi)}[\hat{\nabla}_x f(x,y;\xi)]=\nabla_x f_{\mu_1}(x,y)$ and $\mathbb{E}_{(\tilde{u},\xi)}[\hat{\nabla}_y f(x,y;\xi)]
=\nabla_y f_{\mu_2}(x,y)$. Similarly, we have $\mathbb{E}_{(\hat{U},\mathcal{B})}[\hat{\nabla}_x f(x,y;\mathcal{B})]=\nabla_x f_{\mu_1}(x,y)$ and $\mathbb{E}_{(\tilde{U},\mathcal{B})}[\hat{\nabla}_y f(x,y;\mathcal{B})]
=\nabla_y f_{\mu_2}(x,y)$.
Next, we give some mild assumptions about the problem \eqref{eq:2}.
\begin{assumption} \label{ass:2}
The \textbf{variance of zeroth-order  stochastic  gradient} is bounded, \emph{i.e.,} there exists a constant $\delta_1 >0$ such that for all $x$,
it follows $\mathbb{E}\|\hat{\nabla}_x f(x,y;\xi) - \nabla_x f_{\mu_1}(x,y)\|^2 \leq \delta_1^2$, and for all $y$,
it follows $\mathbb{E}\|\hat{\nabla}_y f(x,y;\xi) - \nabla_y f_{\mu_2}(x,y)\|^2 \leq \delta_1^2$.
The \textbf{variance of stochastic gradient} is bounded, \emph{i.e.,} there exists a constant $\delta_2 >0$ such that for all $x$,
it follows $\mathbb{E}\|\nabla_x f(x,y;\xi) - \nabla_x f(x,y)\|^2 \leq \delta_2^2$; There exists a constant $\delta_2 >0$
such that for all $y$, it follows $\mathbb{E}\|\nabla_y f(x,y;\xi) - \nabla_y f(x,y)\|^2 \leq \delta_2^2$.
\end{assumption}
Assumption \ref{ass:2} is similar to the upper bound of variance of stochastic partial gradients in \citep{luo2020stochastic,wang2020zeroth}.
For notational simplicity, let $\delta=\max(\delta_1,\delta_2)$.
By using Assumption \ref{ass:2}, we have $\mathbb{E}\|\hat{\nabla}_x f(x,y;\mathcal{B}) - \nabla f_{\mu_1}(x,y)\|^2 \leq \delta^2/b$ and
$\mathbb{E}\|\hat{\nabla}_y f(x,y;\mathcal{B}) - \nabla f_{\mu_2}(x,y)\|^2 \leq \delta^2/b$.

\begin{assumption} \label{ass:5}
Each component function $f(x,y;\xi)$ has a $L_f$-Lipschitz gradient, i.e.,
 for all $x,x' \in \mathcal{X}$ and $y,y' \in \mathcal{Y}$
\begin{align}
\|\nabla f(x,y;\xi)-\nabla f(x',y';\xi)\| \leq L_f \|(x,y) - (x',y')\|,
\end{align}
where $\nabla f(x,y;\xi)=\big(\nabla_x f(x,y;\xi), \nabla_y f(x,y;\xi)\big)$.
\end{assumption}
\begin{assumption} \label{ass:6}
The objective function $f(x,y)$ is $\tau$-strongly concave in variable $y$, i.e.,
\begin{align}
 \|\nabla_yf(x,y)-\nabla_yf(x,y')\| \geq \tau\|y-y'\|, \ \forall x\in \mathcal{X}, \ y, y'\in \mathcal{Y}.
\end{align}
Then the following inequality holds
\begin{align}
 f(x,y) \leq f(x,y') + \langle\nabla_y f(x,y'), y-y'\rangle - \frac{\tau}{2}\|y-y'\|^2.
\end{align}
\end{assumption}
Assumption \ref{ass:5} also implies the partial gradients $\nabla_xf(x,y)=\mathbb{E}_{\xi}[\nabla_x f(x,y;\xi)]$ and $\nabla_yf(x,y)=\mathbb{E}_{\xi}[\nabla_y f(x,y;\xi)]$ are
$L_f$-Lipschiz continuous. Since $f(x,y)$ is strongly concave in $y\in \mathcal{Y}$, there exists a unique solution to
the problem $\max_{y\in \mathcal{Y}}f(x,y)$ for any $x$, and we define the solution as
$y^*(x) = \arg\max_{y \in \mathcal{Y}} f(x,y)$, and let $F(x) = \max_{y\in \mathcal{Y}} f(x,y) = f(x,y^*(x))$.
\begin{assumption} \label{ass:7}
The function $F(x)$ is bounded from below in $\mathcal{X}$, \emph{i.e.,} $F^* = \inf_{x\in \mathcal{X}} F(x)$.
\end{assumption}

\section{ Accelerated  Zeroth-Order Momentum Method for Mini-Optimization }
\label{section:4}
In this section, we propose a new accelerated zeroth-order momentum (Acc-ZOM) method to solve
the \textbf{black-box} mini-optimization problem \eqref{eq:1},
where only noise stochastic function values can be obtained.
Although our Acc-ZOM method builds on the momentum-based variance reduction technique of STORM \citep{cutkosky2019momentum},
our Acc-ZOM method is the first to  extend the original STORM method to the constrained optimization.
Algorithm \ref{alg:1} summarizes the algorithmic framework of our Acc-ZOM method.

\setlength{\textfloatsep}{2pt}
\begin{algorithm}[!t]
\caption{ Acc-ZOM Algorithm for  Mini Optimization}
\label{alg:1}
\begin{algorithmic}[1] 
\STATE {\bfseries Input:}  $T$, parameters $\{\gamma, k,m,c\}$ and initial input $x_1\in \mathcal{X}$; \\
\STATE {\bfseries initialize:} Draw a sample $\xi_1$, and sample a vector $u\in \mathbb{R}^d$ from uniform distribution
            over unit sphere, then compute $v_1 = \hat{\nabla} f(x_1;\xi_1)$, where the zeroth-order gradient
is estimated from \eqref{eq:3};\\
\FOR{$t = 1, 2, \ldots, T$}
\STATE Compute $\eta_t = \frac{k}{(m+t)^{1/3}}$;
\IF{$\mathcal{X}=\mathbb{R}^d$}
\STATE Update $x_{t+1} = x_t - \gamma \eta_t v_t$;
\ELSE
\STATE Update $\tilde{x}_{t+1} = \mathcal{P}_{\mathcal{X}}(x_t - \gamma v_t)$, and $x_{t+1} = x_t + \eta_t(\tilde{x}_{t+1}-x_t)$;
\ENDIF
\STATE Compute $\alpha_{t+1} = c\eta_t^2$;
\STATE Draw a sample $\xi_{t+1}$, and sample a vector $u\in \mathbb{R}^d$ from uniform distribution
            over unit sphere, then compute $v_{t+1} = \hat{\nabla} f(x_{t+1};\xi_{t+1}) + (1-\alpha_{t+1})\big[v_t - \hat{\nabla} f(x_t;\xi_{t+1})\big]$,
where the zeroth-order gradients are estimated from \eqref{eq:3}; \\
\ENDFOR
\STATE {\bfseries Output:} (for theoretical) $x_{\zeta}$ chosen uniformly random from $\{x_t\}_{t=1}^{T}$.
\STATE {\bfseries Output:} (for practical) $x_T$.
\end{algorithmic}
\end{algorithm}

In Algorithm \ref{alg:1}, we use the zeroth-order variance-reduced stochastic gradients as follows:
\begin{align}
  v_t &= \alpha_t \hat{\nabla} f(x_t;\xi_t)+ (1 - \alpha_t)\big(\hat{\nabla} f(x_t;\xi_t) - \hat{\nabla} f(x_{t-1};\xi_t) + v_{t-1}\big),
\end{align}
where $\alpha_t\in (0,1]$. When $\alpha_t=1$, $v_t$ will degenerate a vanilla zeroth-order stochastic gradient; When $\alpha_t=0$,
$v_t$ will degenerate a  zeroth-order  stochastic   gradient based on variance-reduced technique of SPIDER \citep{fang2018spider}.
When the constraint set $\mathcal{X}=\mathbb{R}^d$,
\emph{i.e.}, the problem \eqref{eq:1} is an unconstrained problem, we use a common metric $\mathbb{E}\|\nabla f(x_t)\|$ used in the nonconvex optimization \citep{fang2018spider,ji2019improved}
to measure the convergence of Algorithm \ref{alg:1}.

When the constraint set  $\mathcal{X}\subset \mathbb{R}^d$,
at the step 8 of Algorithm \ref{alg:1}, we use $0<\eta_t\leq 1$ to ensure the variable $x_t$ for all $t\geq 1$
in the convex constraint set $\mathcal{X}$.
At the same time, we provide a useful metric $\mathbb{E}[\mathcal{G}_t]$ to measure the convergence properties of our Acc-ZOM for constrained optimization, defined as
\begin{align} \label{eq:12}
 \mathcal{G}_t = \frac{1}{\gamma}\|\tilde{x}_{t+1} - x_t\| + \|\nabla f(x_t) - v_t\|.
\end{align}
In fact, our metric $\mathbb{E}[\mathcal{G}_t]$ is tighter than standard gradient mapping metric $\mathbb{E} \|G_{\mathcal{X}}(x_t,\nabla f(x_t),\gamma)\|$ used in \citep{ghadimi2016mini}, i.e., $\mathcal{G}_t \geq \|G_{\mathcal{X}}(x_t,\nabla f(x_t),\gamma)\|$, where
\begin{align}
 & G_{\mathcal{X}}(x_t,\nabla f(x_t),\gamma) = \frac{1}{\gamma}\big(x_t - \mathcal{P}_{\mathcal{X}}(x_t - \gamma \nabla f(x_t))\big), \nonumber \\
 & \mathcal{P}_{\mathcal{X}}(x_t - \gamma \nabla f(x_t)) = \arg\min_{x\in \mathcal{X}} \bigg\{ \langle\nabla f(x_t), x-x_t\rangle + \frac{1}{2\gamma}\|x-x_t\|^2\bigg\}.
\end{align}
Let $w(x) = \frac{1}{2}\|x\|^2$, as in \citep{ghadimi2016mini}, we give a prox-function associated with $w(x)$, defined as
\begin{align}
 V(x,x_t) = w(x) - \big( w(x_t) + \langle \nabla w(x_t), x-x_t\rangle \big) = \frac{1}{2}\|x-x_t\|^2.
\end{align}
At the same time, the step 8 of Algorithm \ref{alg:1} can be rewritten as
\begin{align}
 \tilde{x}_{t+1}=\mathcal{P}_{\mathcal{X}}(x_t - \gamma v_t) = \arg\min_{x\in \mathcal{X}} \bigg\{ \langle v_t, x-x_t\rangle + \frac{1}{2\gamma}\|x-x_t\|^2\bigg\}.
\end{align}
Then we also can obtain a gradient mapping $ G_{\mathcal{X}}(x_t,v_t,\gamma) = \frac{1}{\gamma}\big(x_t - \mathcal{P}_{\mathcal{X}}(x_t - \gamma v_t)\big)=\frac{1}{\gamma}\big(x_t - \tilde{x}_{t+1}\big)$. Since the function $w(x)= \frac{1}{2}\|x\|^2$ is 1-strongly convex,
we have
\begin{align} \label{eq:MG}
\|G_{\mathcal{X}}(x_t,\nabla f(x_t),\gamma)\| & = \|G_{\mathcal{X}}(x_t,\nabla f(x_t),\gamma) - G_{\mathcal{X}}(x_t,v_t,\gamma) + G_{\mathcal{X}}(x_t,v_t,\gamma) \| \nonumber \\
& \leq \|G_{\mathcal{X}}(x_t,\nabla f(x_t),\gamma) - G_{\mathcal{X}}(x_t,v_t,\gamma)\| + \|G_{\mathcal{X}}(x_t,v_t,\gamma)\| \nonumber \\
& \mathop{\leq}^{(i)} \|\nabla f(x_t)-v_t\| + \|G_{\mathcal{X}}(x_t,v_t,\gamma)\| \nonumber \\
& = \|\nabla f(x_t)-v_t\| + \frac{1}{\gamma}\|x_t - \tilde{x}_{t+1}\|,
\end{align}
where the above inequality $(i)$ holds by Proposition 1 of \citep{ghadimi2016mini}.

In fact, the original STORM method \citep{cutkosky2019momentum} is only competent to
  \textbf{unconstrained} optimization.
In Algorithm \ref{alg:1}, when using  stochastic gradient instead of stochastic zeroth-order gradient for solving
the problem \eqref{eq:1}, our Acc-ZOM algorithm will reduce to a new version of STORM method
for \textbf{constrained} optimization.

\section{Accelerated Zeroth-Order Momentum Descent Ascent Method for Minimax Optimization}
\label{section:5}
In the section, we propose an accelerated zeroth-order momentum descent ascent (Acc-ZOMDA) method
to solve the \textbf{black-box} minimax problem \eqref{eq:2}, where only stochastic function values can be obtained.
In fact, we extend the above Acc-ZOM method to solve the minimax problem and then  obtain the Acc-ZOMDA method.
Algorithm \ref{alg:2} describes the algorithmic framework of our  Acc-ZOMDA method.

\begin{algorithm}[tb]
\caption{ Acc-ZOMDA Algorithm for  Minimax Optimization}
\label{alg:2}
\begin{algorithmic}[1] 
\STATE {\bfseries Input:}  $T$, parameters $\{\gamma, \lambda, k,m,c_1,c_2\}$ and initial input $x_1\in \mathcal{X}$ and $y_1\in\mathcal{Y}$; \\
\STATE {\bfseries initialize:}  Draw a mini-batch samples $\mathcal{B}_1=\{\xi_i^1\}_{i=1}^b$,
and draw vectors $\{\hat{u}_i\in \mathbb{R}^{d_1}\}_{i=1}^b$ and
$\{\tilde{u}_i\in \mathbb{R}^{d_2}\}_{i=1}^b$ from uniform distribution over unit sphere,
then compute $v_1 = \hat{\nabla}_x f(x_1,y_1;\mathcal{B}_1)$ and $w_1 = \hat{\nabla}_y f(x_1,y_1;\mathcal{B}_1)$,
where the zeroth-order gradients are estimated from \eqref{eq:5} and \eqref{eq:6}; \\
\FOR{$t = 1, 2, \ldots, T$}
\STATE Compute $\eta_t = \frac{k}{(m+t)^{1/3}}$;
\IF{$\mathcal{X}=\mathbb{R}^{d_1}$}
\STATE Update $x_{t+1} = x_t - \gamma \eta_t v_t$;
\ELSE
\STATE Update $\tilde{x}_{t+1} = \mathcal{P}_{\mathcal{X}}(x_t - \gamma v_t)$ and $x_{t+1} = x_t + \eta_t(\tilde{x}_{t+1}-x_t)$;
\ENDIF
\STATE Update $\tilde{y}_{t+1} = \mathcal{P}_{\mathcal{Y}}(y_t + \lambda w_t)$ and $y_{t+1} = y_t + \eta_t(\tilde{y}_{t+1}-y_t)$;
\STATE Compute $\alpha_{t+1} = c_1\eta_t^2$ and $\beta_{t+1} = c_2\eta_t^2$;
\STATE Draw a mini-batch samples $\mathcal{B}_{t+1}=\{\xi_i^{t+1}\}_{i=1}^b$, and draw vectors $\{\hat{u}_i\in \mathbb{R}^{d_1}\}_{i=1}^b$ and
$\{\tilde{u}_i\in \mathbb{R}^{d_2}\}_{i=1}^b$
from uniform distribution over unit sphere; \\
\STATE Compute
$v_{t+1} = \hat{\nabla}_x f(x_{t+1},y_{t+1};\mathcal{B}_{t+1}) + (1-\alpha_{t+1})\big[v_t - \hat{\nabla}_x f(x_t,y_t;\mathcal{B}_{t+1})\big]$
and $w_{t+1} = \hat{\nabla}_y f(x_{t+1},y_{t+1};\mathcal{B}_{t+1}) + (1-\beta_{t+1})\big[w_t - \hat{\nabla}_y f(x_t,y_t;\mathcal{B}_{t+1})\big]$,
where the zeroth-order gradients are estimated from \eqref{eq:5} and \eqref{eq:6}. \\
\ENDFOR
\STATE {\bfseries Output:} (for theoretical) $x_{\zeta}$ and $y_{\zeta}$ chosen uniformly random from $\{x_t, y_t\}_{t=1}^{T}$.
\STATE {\bfseries Output:} (for practical) $x_T$ and $y_T$.
\end{algorithmic}
\end{algorithm}

In Algorithm \ref{alg:2}, we use the momentum-based variance reduced technique of STORM
to estimate the stochastic zeroth-order partial gradients $v_t$ and $w_t$.
When the constraint set $\mathcal{X}=\mathbb{R}^{d_1}$,
\emph{i.e.}, the problem \eqref{eq:2} is an unconstrained problem w.r.t. variable $x$,
we use a common metric $\mathbb{E}\|\nabla F(x_t)\|$ used in \citep{lin2019gradient,wang2020zeroth}
to measure the convergence of Algorithm \ref{alg:2}, where the function $F(x)=\max_{y\in \mathcal{Y}}f(x,y)$.

When the constraint set $\mathcal{X}\subset \mathbb{R}^{d_1}$,
we define a useful metric $\mathbb{E}[\mathcal{H}_t]$ to measure the convergence properties of our Acc-ZOMDA Algorithm,
\begin{align} \label{eq:14}
 \mathcal{H}_t =  \frac{1}{\gamma}\|\tilde{x}_{t+1} - x_t\| + \|\nabla_x f(x_t,y_t) - v_t\| + L_f\|y_t-y^*(x_t)\|,
\end{align}
where the first two terms of $\mathcal{H}_t$ measure  convergence of the iteration solutions $\{x_t\}_{t=1}^T$,
and the last term measures  convergence of the iteration solutions $\{y_t\}_{t=1}^T$. In fact, our new metric  $\mathbb{E}[\mathcal{H}_t]$ is tighter than the generic gradient mapping metric $\mathbb{E}\|G_{\mathcal{X}}(x_t,\nabla F(x_t),\gamma)\|$, i.e.,
 $\mathcal{H}_t \geq \|G_{\mathcal{X}}(x_t,\nabla F(x_t),\gamma)\|$, where $G_{\mathcal{X}}(x_t,\nabla F(x_t),\gamma)$ is a gradient mapping, defined as
\begin{align}
 & G_{\mathcal{X}}(x_t,\nabla F(x_t),\gamma) = \frac{1}{\gamma}\big(x_t - \mathcal{P}_{\mathcal{X}}(x_t - \gamma \nabla F(x_t))\big), \nonumber \\
 & \mathcal{P}_{\mathcal{X}}(x_t - \gamma \nabla F(x_t)) = \arg\min_{x\in \mathcal{X}} \bigg\{ \langle\nabla F(x_t), x-x_t\rangle + \frac{1}{2\gamma}\|x-x_t\|^2\bigg\}, \label{eq:7}
\end{align}
where $F(x_t)=f(x_t,y^*(x_t))=\min_{y\in \mathcal{Y}}f(x_t,y)$. At the same time, the step 8 of Algorithm \ref{alg:2} can be rewritten as
\begin{align}
 \tilde{x}_{t+1}=\mathcal{P}_{\mathcal{X}}(x_t - \gamma v_t) = \arg\min_{x\in \mathcal{X}} \bigg\{ \langle v_t, x-x_t\rangle + \frac{1}{2\gamma} \|x-x_t\|^2\bigg\}.
\end{align}
Then we also can obtain a gradient mapping  $G_{\mathcal{X}}(x_t,v_t,\gamma) = \frac{1}{\gamma}\big(x_t - \mathcal{P}_{\mathcal{X}}(x_t - \gamma v_t)\big)=\frac{1}{\gamma}\big(x_t - \tilde{x}_{t+1}\big)$.
Since the function $w(x)=\frac{1}{2}\|x\|^2$ is 1-strongly convex,
we have
\begin{align} \label{eq:MH}
\|G_{\mathcal{X}}(x_t,\nabla F(x_t),\gamma)\| & = \|G_{\mathcal{X}}(x_t,\nabla F(x_t),\gamma) - G_{\mathcal{X}}(x_t,v_t,\gamma) + G_{\mathcal{X}}(x_t,v_t,\gamma)\| \nonumber \\
& \leq \|G_{\mathcal{X}}(x_t,\nabla F(x_t),\gamma) - G_{\mathcal{X}}(x_t,v_t,\gamma)\| + \|G_{\mathcal{X}}(x_t,v_t,\gamma)\|  \nonumber \\
& \mathop{\leq}^{(i)} \|\nabla F(x_t)-v_t\| + \|G_{\mathcal{X}}(x_t,v_t,\gamma)\| \nonumber \\
& = \|\nabla F(x_t)-\nabla_xf(x_t,y_t) + \nabla_xf(x_t,y_t) - v_t\| + \frac{1}{\gamma}\|x_t - \tilde{x}_{t+1}\| \nonumber \\
& \leq \|\nabla_x f(x_t,y^*(x_t))\!-\!\nabla_xf(x_t,y_t)\| + \|\nabla_xf(x_t,y_t) - v_t\| \!+\! \frac{1}{\gamma}\|x_t - \tilde{x}_{t+1}\| \nonumber \\
& \mathop{\leq}^{(ii)} L_f\|y^*(x_t)-y_t\| + \|\nabla_xf(x_t,y_t)-v_t\| + \frac{1}{\gamma}\|x_t - \tilde{x}_{t+1}\|,
\end{align}
where the above inequality $(i)$ holds by Proposition 1 of \citep{ghadimi2016mini}, and the above inequality $(ii)$ is due to Assumption 5.

\begin{algorithm}[tb]
\caption{ Acc-MDA Algorithm for  Minimax Optimization}
\label{alg:4}
\begin{algorithmic}[1] 
\STATE {\bfseries Input:} $T$, parameters $\{\gamma, \lambda, k,m,c_1,c_2\}$
and initial input $x_1 \in \mathcal{X}$ and $y_1 \in \mathcal{Y}$; \\
\STATE {\bfseries initialize:} Draw a mini-batch samples $\mathcal{B}_1=\{\xi_i^1\}_{i=1}^b$,
and then compute stochastic gradients $v_1 = \nabla_x f(x_1,y_1;\mathcal{B}_1)$ and $w_1 = \nabla_y f(x_1,y_1;\mathcal{B}_1)$;  \\
\FOR{$t = 1, 2, \ldots, T$}
\STATE Compute $\eta_t = \frac{k}{(m+t)^{1/3}}$;
\IF {$\mathcal{X}=\mathbb{R}^{d_1}$}
\STATE Update $x_{t+1} = x_t - \gamma \eta_t v_t$;
\ELSE
\STATE Update $\tilde{x}_{t+1} = \mathcal{P}_{\mathcal{X}}(x_t - \gamma v_t)$ and $x_{t+1} = x_t + \eta_t(\tilde{x}_{t+1}-x_t)$;
\ENDIF
\STATE Update $\tilde{y}_{t+1} = \mathcal{P}_{\mathcal{Y}}(y_t + \lambda w_t)$ and $y_{t+1} = y_t + \eta_t(\tilde{y}_{t+1}-y_t)$;
\STATE Compute $\alpha_{t+1} = c_1\eta_t^2$ and $\beta_{t+1} = c_2\eta_t^2$;
\STATE Draw a mini-batch samples $\mathcal{B}_{t+1}=\{\xi^{t+1}_i\}_{i=1}^b$, and then compute stochastic gradients
  $v_{t+1} = \nabla_x f(x_{t+1},y_{t+1};\mathcal{B}_{t+1}) + (1-\alpha_{t+1})\big[v_t - \nabla_x f(x_t,y_t;\mathcal{B}_{t+1})\big]$ and
          $w_{t+1} = \nabla_y f(x_{t+1},y_{t+1};\mathcal{B}_{t+1}) + (1-\beta_{t+1})\big[w_t - \nabla_y f(x_t,y_t;\mathcal{B}_{t+1})\big]$; \\
\ENDFOR
\STATE {\bfseries Output:} (for theoretical) $x_{\zeta}$ and $y_{\zeta}$ chosen uniformly random from $\{x_t, y_t\}_{t=1}^{T}$.
\STATE {\bfseries Output:} (for practical) $x_T$ and $y_T$.
\end{algorithmic}
\end{algorithm}

\section{ Accelerated First-Order Momentum Descent Ascent Method for Minimax Optimization }
\label{section:6}
In this section, we propose an  accelerated  first-order momentum descent ascent (Acc-MDA) method
to solve the \textbf{transparent} minimax problem \eqref{eq:2},
whose explicit stochastic gradients are accessible. 
Algorithm \ref{alg:4} gives the algorithmic framework of our Acc-MDA method.
In Algorithm \ref{alg:4}, we use the stochastic gradients instead of the stochastic zeroth-order gradients
used in Algorithm \ref{alg:2}.
In our Acc-MDA algorithm, we use the momentum-based variance-reduced technique of STORM to estimate the partial derivatives $v_t$ and $w_t$ on variables $x$ and $y$, respectively. Moreover, our Acc-MDA algorithm also uses the momentum iteration to update variables $x$ and $y$ as follows:
\begin{align}
&\tilde{x}_{t+1} = \mathcal{P}_{\mathcal{X}}(x_t - \gamma v_t), \quad  x_{t+1} = x_t + \eta_t(\tilde{x}_{t+1}-x_t), \\
& \tilde{y}_{t+1} = \mathcal{P}_{\mathcal{Y}}(y_t + \lambda w_t), \quad y_{t+1} = y_t + \eta_t(\tilde{y}_{t+1}-y_t). 
\end{align}
At the same time, at step 6 of Algorithm \ref{alg:4}, i.e., $x_{t+1} = x_t - \gamma \eta_t v_t$ also can be rewritten as $\tilde{x}_{t+1}=x_t - \gamma v_t$ and $x_{t+1} = x_t + \eta_t(\tilde{x}_{t+1}-x_t)$. 

By combining Algorithms \ref{alg:2} and \ref{alg:4},
we can propose an accelerated semi-zeroth-order momentum descent ascent (Acc-Semi-ZOMDA) method
to solve one-sided black-box problem \eqref{eq:2} studied in \citep{liu2019min},
where the explicit stochastic partial gradients in variable $x$ can not be accessible.
Specifically, in the Acc-Semi-ZOMDA algorithm, we only use the stochastic partial gradients $w_t$ instead of
the stochastic zeroth-order partial gradients $w_t$ in Algorithm \ref{alg:2}.

\section{ Convergence Analysis}
\label{section:7}
In this section, we study the convergence properties of our algorithms (Acc-ZOM, Acc-ZOMDA and Acc-MDA) under some mild conditions.
\subsection{Convergence Analysis of the Acc-ZOM Algorithm}
In this subsection, we analyze convergence of our \textbf{Acc-ZOM} algorithm for solving the \emph{constrained} and
\emph{unconstrained} \textbf{mini-optimization} problem \eqref{eq:1}, respectively.
\subsubsection{Convergence Analysis of the Acc-ZOM Algorithm for Constrained Mini-Optimization }
In the subsection, we analyze convergence properties of the Acc-ZOM algorithm for solving the  \textbf{constrained}
 problem \eqref{eq:1}, i.e., $\mathcal{X}\subset \mathbb{R}^d$.
The following convergence results build on a new metric $\mathbb{E}[\mathcal{G}_t]$,
where $\mathcal{G}_t$ is defined in \eqref{eq:12}.
The related proofs of these convergence analysis are provided in Appendix \ref{Appendix:A1}.

We begin with defining a
function $f_{\mu}(x)=\mathbb{E}_{u\sim U_B}[f(x+\mu u)]$,
which is a smooth approximation of function $f(x)$,
where $U_B$ is the uniform distribution over the $d$-dimensional unit Euclidean ball $B$.

\begin{theorem} \label{th:1}
Suppose the sequence $\{x_t\}_{t=1}^T$ be generated from Algorithm \ref{alg:1}. When $\mathcal{X}\subset \mathbb{R}^d$, and let $\eta_t = \frac{k}{(m+t)^{1/3}}$
for all $t\geq 0$, $0< \gamma \leq \min\big(\frac{m^{1/3}}{2Lk},\frac{1}{2\sqrt{6d}L}\big)$,
$c\geq \frac{2}{3k^3} + \frac{5}{4}$, $k>0$, $m\geq \max\big(2, (ck)^3, k^3\big)$ and $0<\mu \leq \frac{1}{d(m+T)^{2/3}}$,
we have
 \begin{align}
  \frac{1}{T}\sum_{t=1}^T\mathbb{E} \|G_{\mathcal{X}}(x_t,\nabla f(x_t),\gamma)\| \leq \frac{1}{T}\sum_{t=1}^T\mathbb{E}[\mathcal{G}_t] \leq  \frac{\sqrt{2M}m^{1/6}}{T^{1/2}} + \frac{\sqrt{2M}}{T^{1/3}} + \frac{L}{2(m+T)^{2/3}},
 \end{align}
where $M=\frac{f_{\mu}(x_1) - f^*}{k\gamma} + \frac{m^{1/3}\sigma^2}{k^2} + \frac{9L^2}{4k^2} + 2k^2c^2\sigma^2\ln(m+T)$.
\end{theorem}

\begin{remark}
Without loss of generality, let $m\geq \max \big( 2, (ck)^3, k^3, (\frac{k}{\sqrt{6d}})^3\big)$,
we have $\frac{m^{1/3}}{2Lk} \geq \frac{1}{2\sqrt{6d}L}$.
It is easy verified that $\gamma=O(\frac{1}{\sqrt{d}})$, $c=O(1)$ and $m=O(1)$.
Then we have $M=O\big(\sqrt{d}+\ln(m+T)\big)=\tilde{O}\big(\sqrt{d}\big)$.
Thus, the Acc-ZOM algorithm has $\tilde{O}\big(\frac{d^{1/4}}{T^{1/3}}\big)$ convergence rate.
By $\frac{d^{1/4}}{T^{1/3}} \leq \epsilon$, i.e., $\mathbb{E}[\mathcal{G}_\zeta] \leq \epsilon$, we
choose $T \geq d^{3/4}\epsilon^{-3}$. In Algorithm \ref{alg:1}, we require to query four function values for estimating the zeroth-order gradients
$v_t$ at each iteration, and need $T$ iterations.
Thus, the Acc-ZOM algorithm has a query complexity of $4T=\tilde{O}(d^{3/4}\epsilon^{-3})$ for finding an $\epsilon$-stationary point.
\end{remark}

\subsubsection{ Convergence Analysis of Acc-ZOM Algorithm for Unconstrained Mini-Optimization }
In this subsection, we study the convergence properties of our Acc-ZOM algorithm for solving the 
\textbf{unconstrained} problem \eqref{eq:1}, i.e., $\mathcal{X}= \mathbb{R}^d$.
The following convergence analysis builds on the common metric $\mathbb{E}\|\nabla f(x)\|$ used in nonconvex optimization \citep{ji2019improved}.
The related proofs of these convergence analysis are provided in Appendix \ref{Appendix:A2}.

 \begin{theorem} \label{th:01}
 Suppose the sequence $\{x_t\}_{t=1}^T$ be generated from Algorithm \ref{alg:1}.
 When $\mathcal{X}=\mathbb{R}^d$, and let $\eta_t = \frac{k}{(m+t)^{1/3}}$
 for all $t\geq 0$, $0< \gamma \leq \min\big(\frac{m^{1/3}}{2Lk},\frac{1}{2\sqrt{6d}L}\big)$,
 $c\geq \frac{2}{3k^3} + \frac{5}{4}$, $k>0$, $m\geq \max\big(2, k^3, (ck)^3\big)$ and $0<\mu \leq \frac{1}{d(m+T)^{2/3}}$,
 we have
 \begin{align}
  \frac{1}{T} \sum_{t=1}^T\mathbb{E}\|\nabla f(x_t)\|
  \leq  \frac{\sqrt{2M}m^{1/6}}{T^{1/2}} + \frac{\sqrt{2M}}{T^{1/3}} + \frac{L}{2(m+T)^{2/3}},
\end{align}
 where $M=\frac{f_{\mu}(x_1) - f^*}{k\gamma} + \frac{m^{1/3}\sigma^2}{k^2} + \frac{9L^2}{4k^2} + 2k^2c^2\sigma^2\ln(m+T)$.
 \end{theorem}

\begin{remark}
Since the conditions of Theorem \ref{th:01} are the same conditions of Theorem \ref{th:1},
Theorem \ref{th:01} also show that
our Acc-ZOM algorithm has a lower query complexity of $\tilde{O}(d^{3/4}\epsilon^{-3})$ for finding an $\epsilon$-stationary point.
\end{remark}

\subsection{Convergence Analysis of the Acc-ZOMDA Algorithm}
In this subsection, we analyze convergence of our \textbf{Acc-ZOMDA} algorithm for solving the
 \emph{constrained} and
\emph{unconstrained} \textbf{minimax-optimization} problem \eqref{eq:2}, respectively.
\subsubsection{Convergence Analysis of the Acc-ZOMDA Algorithm for Constrained Minimax Optimization}
In the subsection, we provide the convergence properties of our  Acc-ZOMDA algorithm
for solving the  \textbf{constrained} minimax  problem \eqref{eq:2}, i.e.,
$\mathcal{X}\subset \mathbb{R}^{d_1}$ and $\mathcal{Y}\subset \mathbb{R}^{d_2}$ (or $\mathcal{Y}= \mathbb{R}^{d_2}$). The following results build on
new convergence metric $\mathbb{E}[\mathcal{H}_t]$, where $\mathcal{H}_t$ is defined as in
\eqref{eq:14}. The related proofs of these convergence analysis are provided in Appendix \ref{Appendix:A3}.

We first define a function $F_{\mu_1}(x) = \mathbb{E}_{u_1\sim U_{B_1}}[F(x+\mu_1u_1)]$, which is a smoothing approximation of
the function $F(x)= f(x,y^*(x))=\max_{y \in \mathcal{Y}} f(x,y) $. For notational simplicity, let $\tilde{d}=d_1+d_2$,
$L_g=L_f+\frac{L^2_f}{\tau}$ and $\kappa_y = L_f/\tau$ denote the condition number for function $f(\cdot,y)$.

\begin{theorem} \label{th:2}
Suppose the sequence $\{x_t,y_t\}_{t=1}^T$ be generated from Algorithm \ref{alg:2}. When $\mathcal{X}\subset \mathbb{R}^{d_1}$, and let  $\eta_t = \frac{k}{(m+t)^{1/3}}$
for all $t\geq 0$, $c_1 \geq \frac{2}{3k^3} + \frac{9\tau^2}{4}$ and $c_2 \geq \frac{2}{3k^3} + \frac{625\tilde{d}L^2_f}{3b}$,
$k>0$, $1\leq b\leq \tilde{d}$, $m\geq \max\big( 2, k^3, (c_1k)^3, (c_2k)^3\big)$,  $0<\lambda\leq \min\big(\frac{1}{6L_f},\frac{75\tau}{24}\big)$,
$0< \gamma \leq \min\big( \frac{\lambda\tau}{2L_f}\sqrt{\frac{6b/\tilde{d}}{36 \lambda^2 + 625\kappa_y^2}}, \frac{m^{1/3}}{2L_gk}\big)$, $0<\mu_1\leq \frac{1}{d_1(m+T)^{2/3}}$
and $0<\mu_2\leq \frac{1}{\tilde{d}^{1/2}d_2(m+T)^{2/3}}$,
we have
\begin{align}
 \frac{1}{T} \sum_{t=1}^T \mathbb{E}  \|G_{\mathcal{X}}(x_t,\nabla F(x_t),\gamma)\| \leq \frac{1}{T} \sum_{t=1}^T \mathbb{E}[\mathcal{H}_t] \leq  \frac{2\sqrt{3M'}m^{1/6}}{T^{1/2}} + \frac{2\sqrt{3M'}}{T^{1/3}} + \frac{L_f}{2(m+T)^{2/3}}.
\end{align}
where $\Delta_1=\|y_1-y^*(x_1)\|^2$ and  $M' =  \frac{F_{\mu_1}(x_1) - F^*}{\gamma k} + \frac{ 25\tilde{d}L^2_f\Delta_1}{ k\lambda\tau b}+  \frac{2m^{1/3}\delta^2}{b\tau^2 k^2} + \frac{36\tau^2L_f^2 + 625L^4_f}{8b\tau^2}(m+T)^{-2/3} + \frac{9L^2_f}{4b\tau^2 k^2} + \frac{2(c_1^2+c_2^2)\delta^2 k^2}{b\tau^2}\ln(m+T)$.
\end{theorem}

\begin{remark}
Without loss of generality, let $m\geq \max \big( \big(L_g\lambda\tau k\sqrt{\frac{6b/\tilde{d}}{36\lambda^2 + 625\kappa_y^2}}\big)^3,2, (c_1k)^3, (c_2k)^3, \\ k^3\big)$ and $\tau \leq \frac{1}{L_f}$.  
It is easy verified that $k=O(1)$, $\lambda=O(\tau)$, $\gamma^{-1}=O(\sqrt{\frac{\tilde{d}}{b}}\kappa_y^{3})$, $c_1=O(1)$,
$c_2=O(\frac{\tilde{d}}{b}L_f^2)$ and $m=O(\frac{\tilde{d}^3}{b^3}L_f^6)$. Then we have
$M'=O(\sqrt{\frac{\tilde{d}}{b}}\kappa_y^{3}+ \frac{\tilde{d}}{b}\kappa^2_y+ \frac{\tilde{d}}{b^2}\kappa^2_y +\frac{\kappa^2_y}{b}(m+T)^{-2/3} + \frac{\kappa^2_y}{b} +\frac{\tilde{d}^2}{b^3} \kappa^2_y\ln(m+T))$. Note that in $M'$, we only keep $b$, $\tilde{d}$, $T$ and $\kappa_y$ terms.
\textbf{When $b=1$, we have $M'=\tilde{O}\big(\sqrt{\tilde{d}}\kappa_y^3+\tilde{d}^2\kappa_y^2\big)$.} When $\kappa_y \geq \tilde{d}^{3/2}$,  the Acc-ZOMDA algorithm has a convergence rate of $\tilde{O}\big(\frac{\kappa_y^{3/2}\tilde{d}^{1/4}}{T^{1/3}}\big)$.
By $\frac{\kappa_y^{3/2}\tilde{d}^{1/4}}{T^{1/3}} \leq \epsilon$, i.e., $\mathbb{E}[\mathcal{H}_\zeta]\leq \epsilon$, we
choose $T \geq  \kappa_y^{4.5}\tilde{d}^{3/4}\epsilon^{-3}$. In Algorithm \ref{alg:2}, we need to query eight function values for estimating the zeroth-order gradients $v_t$ and $w_t$ at each iteration, and need $T$ iterations.
Thus, the Acc-ZOMDA algorithm has a query complexity of $8T=\tilde{O}\big(\kappa_y^{4.5}\tilde{d}^{3/4}\epsilon^{-3}\big)$ for finding an $\epsilon$-stationary point.
When $1\leq \kappa_y \leq \tilde{d}^{3/2}$, the Acc-ZOMDA algorithm has a convergence rate of $\tilde{O}\big(\frac{\kappa_y\tilde{d}}{T^{1/3}}\big)$. Similarly, the Acc-ZOMDA algorithm has a query complexity of $8T=\tilde{O}\big(\kappa_y^3\tilde{d}^3\epsilon^{-3}\big)$ for finding an $\epsilon$-stationary point.
\end{remark}

\subsubsection{Convergence Analysis of the Acc-ZOMDA Algorithm for Unconstrained Minimax Optimization}
In the subsection, we further provide the convergence properties of our  Acc-ZOMDA algorithm for solving
the \textbf{unconstrained} minimax  problem \eqref{eq:2},
i.e., $\mathcal{X}= \mathbb{R}^{d_1}$ and $\mathcal{Y}= \mathbb{R}^{d_2}$ (or $\mathcal{Y}\subset \mathbb{R}^{d_2}$). 
The following convergence results build on the common metric $\mathbb{E}\|\nabla F(x)\|$
used in \citep{lin2019gradient,wang2020zeroth}, where $F(x)=\max_{y\in \mathcal{Y}}f(x,y)$.
The related proofs of these convergence analysis are provided in Appendix \ref{Appendix:A4}.

\begin{theorem} \label{th:02}
Suppose the sequence $\{x_t,y_t\}_{t=1}^T$ be generated from Algorithm \ref{alg:2}.
When $\mathcal{X}=\mathbb{R}^{d_1}$, and let $\eta_t = \frac{k}{(m+t)^{1/3}}$
for all $t\geq 0$, $c_1 \geq \frac{2}{3k^3} + \frac{9\tau^2}{4}$ and $c_2 \geq \frac{2}{3k^3} + \frac{625\tilde{d}L^2_f}{3b}$,
$k>0$, $1\leq b\leq \tilde{d}$, $m\geq \max\big( 2, k^3, (c_1k)^3, (c_2k)^3\big)$,  $0<\lambda\leq \min\big(\frac{1}{6L_f},\frac{75\tau}{24}\big)$,
$0< \gamma \leq \min\big( \frac{\lambda\tau}{2L_f}\sqrt{\frac{6b/\tilde{d}}{36 \lambda^2 + 625\kappa_y^2}}, \frac{m^{1/3}}{2L_gk}\big)$, $0<\mu_1\leq \frac{1}{d_1(m+T)^{2/3}}$
and $0<\mu_2\leq \frac{1}{\tilde{d}^{1/2}d_2(m+T)^{2/3}}$,
we have
\begin{align}
 \frac{1}{T} \sum_{t=1}^T \mathbb{E}\|\nabla F(x_t)\| \leq  \frac{\sqrt{2M'}m^{1/6}}{T^{1/2}} + \frac{\sqrt{2M'}}{T^{1/3}} + \frac{L_f}{2(m+T)^{2/3}},
\end{align}
where $\Delta_1=\|y_1-y^*(x_1)\|^2$ and $M' = \frac{F_{\mu_1}(x_1) - F^*}{\gamma k}  + \frac{ 25\tilde{d}L^2_f\Delta_1}{k\lambda\tau b} +\frac{2m^{1/3}\delta^2}{b\tau^2 k^2} + \frac{36\tau^2L_f^2 + 625L^4_f}{8b\tau^2}(m+T)^{-2/3} + \frac{9L^2_f}{4b\tau^2 k^2} + \frac{2(c_1^2+c_2^2)\delta^2 k^2}{b\tau^2}\ln(m+T)$.
\end{theorem}

\begin{remark}
Since the conditions of Theorem \ref{th:02} are the same conditions of Theorem \ref{th:2},
Theorem \ref{th:02} has the same results of Theorem \ref{th:2}.
\textbf{When $b=1$, we have $M'=\tilde{O}\big(\sqrt{\tilde{d}}\kappa_y^3+\tilde{d}^2\kappa_y^2\big)$.} When $\kappa_y \geq \tilde{d}^{3/2}$,  the Acc-ZOMDA algorithm has a convergence rate of $\tilde{O}\big(\frac{\kappa_y^{3/2}\tilde{d}^{1/4}}{T^{1/3}}\big)$.
By $\frac{\kappa_y^{3/2}\tilde{d}^{1/4}}{T^{1/3}} \leq \epsilon$, i.e., $\mathbb{E}\|\nabla F(x_\zeta)\| \leq \epsilon$, we
choose $T \geq  \kappa_y^{4.5}\tilde{d}^{3/4}\epsilon^{-3}$. In Algorithm \ref{alg:2}, we need to query eight function values for estimating the zeroth-order gradients $v_t$ and $w_t$ at each iteration, and need $T$ iterations.
Thus, the Acc-ZOMDA algorithm has a query complexity of $8T=\tilde{O}\big(\kappa_y^{4.5}\tilde{d}^{3/4}\epsilon^{-3}\big)$ for finding an $\epsilon$-stationary point.
When $1\leq \kappa_y \leq \tilde{d}^{3/2}$, the Acc-ZOMDA algorithm has a convergence rate of $\tilde{O}\big(\frac{\kappa_y\tilde{d}}{T^{1/3}}\big)$. Similarly, the Acc-ZOMDA algorithm has a query complexity of $8T=\tilde{O}\big(\kappa_y^3\tilde{d}^3\epsilon^{-3}\big)$ for finding an $\epsilon$-stationary point.
\end{remark}

\subsection{Convergence Analysis of the Acc-MDA Algorithm}
In the subsection, we analyze convergence of our \textbf{Acc-MDA} algorithm for solving the
 \emph{constrained} and
\emph{unconstrained} \textbf{ minimax-optimization} problem \eqref{eq:2}, respectively. 
\subsubsection{Convergence Analysis of the Acc-MDA Algorithm for Constrained Minimax Optimization}
In the subsection, we give the convergence properties of our  Acc-MDA algorithm for solving the  \textbf{constrained} minimax 
problem \eqref{eq:2}, i.e., $\mathcal{X}\subset \mathbb{R}^{d_1}$ and $\mathcal{Y}\subset \mathbb{R}^{d_2}$ (or $\mathcal{Y}= \mathbb{R}^{d_2}$).
The following convergence results build on a new metric $\mathbb{E}[\mathcal{H}_t]$,
where $\mathcal{H}_t$ is defined in \eqref{eq:14}.
The related proofs of these convergence analysis are provided in Appendix \ref{Appendix:A5}.

\begin{theorem} \label{th:4}
Suppose the sequence $\{x_t,y_t\}_{t=1}^T$ be generated from Algorithm \ref{alg:4}. When $\mathcal{X}\subset \mathbb{R}^{d_1}$, and 
$\eta_t = \frac{k}{(m+t)^{1/3}}$
for all $t\geq 0$, $c_1 \geq \frac{2}{3k^3} + \frac{9\tau^2}{4}$
and $c_2 \geq \frac{2}{3k^3} + \frac{75L^2_f}{2}$, $k>0$, $m\geq \max\big(2, k^3, (c_1k)^3, (c_2k)^3\big)$,
$0<\lambda\leq \min\big(\frac{1}{6L_f}, \frac{27b\tau}{16}\big)$ and
$0< \gamma \leq \min\big( \frac{\lambda\tau}{2L_f}\sqrt{\frac{2b}{8\lambda^2 + 75\kappa_y^2b}}, \frac{m^{1/3}}{2L_gk}\big)$,
we have
\begin{align}
  \frac{1}{T} \sum_{t=1}^T \mathbb{E}  \|G_{\mathcal{X}}(x_t,\nabla F(x_t),\gamma)\| \leq \frac{1}{T} \sum_{t=1}^T \mathbb{E}[\mathcal{H}_t] \leq  \frac{2\sqrt{3M''}m^{1/6}}{T^{1/2}} + \frac{2\sqrt{3M''}}{T^{1/3}},
\end{align}
where $\Delta_1=\|y_1-y^*(x_1)\|^2$ and  $M'' =  \frac{F(x_1) - F^*}{\gamma k} +  \frac{9L^2_f\Delta_1}{ k\lambda\tau}+\frac{2m^{1/3}\delta^2}{b \tau^2k^2} + \frac{2(c_1^2+c_2^2)\delta^2 k^2}{b\tau^2}\ln(m+T)$.
\end{theorem}

\begin{remark}
Without loss of generality, let $\frac{\lambda\tau}{2L_f}\sqrt{\frac{2b}{8\lambda^2 + 75\kappa_y^2b}}\leq \frac{m^{1/3}}{2L_gk}$, we have
$m\geq \max \big( 2, k^3, \\(c_1k)^3, (c_2k)^3, \big(\frac{L_g\lambda\tau k}{L_f}\sqrt{\frac{2b}{8\lambda^2 + 75\kappa_y^2b}}\big)^3\big)$. Let
$\gamma = \frac{\lambda\tau}{2L_f}\sqrt{\frac{2b}{8\lambda^2 + 75\kappa_y^2b}} = \frac{\lambda}{2\kappa_y}\sqrt{\frac{2b}{8\lambda^2 + 75\kappa_y^2b}}$ and $\lambda=\min\big(\frac{1}{6L_f}, \frac{27b\tau}{16}\big)$. Without loss of generality, let $\tau \leq
\frac{1}{L_f}$.
\textbf{When $b=1$}, it is easy verified that $k=O(1)$, $\lambda=O(\tau)$, $\gamma^{-1}=O(\kappa_y^3)$,  $c_1=O(1)$,
$c_2=O(L_f^2)$ and $m=O(L_f^6)$. Then we have
$M''=O(\kappa_y^3+ \kappa^2_y + \kappa^2_y + \kappa^2_y\ln(m+T)) = O(\kappa_y^3)$.
Thus, the Acc-MDA algorithm has a  convergence rate of $O\big(\frac{\kappa_y^{3/2}}{T^{1/3}}\big)$.
By $\frac{\kappa_y^{3/2}}{T^{1/3}} \leq \epsilon$, i.e., $\mathbb{E}[\mathcal{H}_\zeta] \leq \epsilon$, we
choose $T \geq  \kappa_y^{4.5}\epsilon^{-3}$. In Algorithm \ref{alg:4}, we need to compute four stochastic partial gradients to obtain gradient estimators $v_t$ and $w_t$ at each iteration, and need $T$ iterations.
Thus, the Acc-MDA algorithm has a gradient complexity of $4\cdot T=\tilde{O}\big(\kappa_y^{4.5}\epsilon^{-3}\big)$
for finding an $\epsilon$-stationary point.
\end{remark}

\begin{corollary}
Under the same conditions of Theorem \ref{th:4}, when $b=O(\kappa_y^{\nu})$ for $\nu>0$ and $\frac{27b\tau}{16} \leq \frac{1}{6L_f}$, i.e., $\kappa_y^{\nu}\leq \frac{8}{81L_f\tau}$, our Acc-MDA algorithm has a lower gradient complexity of
$\tilde{O}\big(\kappa_y^{(3-\nu/2)}\epsilon^{-3}\big)$ for finding an $\epsilon$-stationary point.
\end{corollary}

\begin{proof}
Under the above conditions of Theorem \ref{th:4}, without loss of generality, let $\frac{\lambda\tau}{2L_f}\sqrt{\frac{2b}{8\lambda^2 + 75\kappa_y^2b}} \\ \leq \frac{m^{1/3}}{2L_gk}$, we have
$m\geq \max \big( 2, k^3, (c_1k)^3, (c_2k)^3, \big(\frac{L_g\lambda\tau k}{L_f}\sqrt{\frac{2b}{8\lambda^2 + 75\kappa_y^2b}}\big)^3\big)$. Let
$\gamma = \frac{\lambda\tau}{2L_f}\sqrt{\frac{2b}{8\lambda^2 + 75\kappa_y^2b}} \\ = \frac{\lambda}{2\kappa_y}\sqrt{\frac{2b}{8\lambda^2 + 75\kappa_y^2b}}$ and $\lambda=\min\big(\frac{1}{6L_f}, \frac{27b\tau}{16}\big)$.

Given $b=O(\kappa_y^{\nu})$ for $\nu>0$ and $\frac{27b\tau}{16} \leq \frac{1}{6L_f}$, i.e., $\kappa_y^{\nu}\leq \frac{8}{81L_f\tau}$, it is easy verified that $k=O(1)$, $\lambda=O(b\tau)$, $\gamma^{-1}=O(\frac{\kappa_y^3}{b})$, $c_1=O(1)$ and
$c_2=O(L_f^2)$. Since $L_g=L_f+\frac{L^2_f}{\tau}$, we have $\frac{L_g\lambda\tau k}{L_f}\sqrt{\frac{2b}{8\lambda^2 + 75\kappa_y^2b}} =
(1+\kappa_y)\lambda\tau k\sqrt{\frac{2b}{8\lambda^2 + 75\kappa_y^2b}} = O(\frac{b}{\kappa_y}) $, we have $m=\max(L_f^6,\frac{b^3}{\kappa_y^3})$. Then we have
$M''=O(\frac{\kappa_y^3}{b}+ \frac{\kappa^2_y}{b} + \frac{\kappa^2_y}{b} + \frac{\kappa^2_y}{b}\ln(m+T)) = O(\frac{\kappa_y^3}{b})=O(\kappa_y^{(3-\nu)})$.  Thus, our Acc-MDA algorithm has a convergence rate of $\tilde{O}\big(\frac{\kappa_y^{(3/2-\nu/2)}}{T^{1/3}}\big)$.
By $\frac{\kappa_y^{(3/2-\nu/2)}}{T^{1/3}} \leq \epsilon$, i.e., $\mathbb{E}[\mathcal{H}_\zeta] \leq \epsilon$, we
choose $T \geq  \kappa_y^{(4.5-3\nu/2)}\epsilon^{-3}$.
Thus, our Acc-MDA algorithm reaches a lower gradient complexity of $4b\cdot T=\tilde{O}\big(\kappa_y^{(4.5-\nu/2)}\epsilon^{-3}\big)$ for finding an $\epsilon$-stationary point.
\end{proof}

\subsubsection{ Convergence Analysis of Acc-MDA Algorithm for  Unconstrained Minimax Optimization }
In the subsection, we further give the convergence properties of our  Acc-MDA algorithm for solving
the \textbf{unconstrained} minimax  problem \eqref{eq:2},
i.e., $\mathcal{X}= \mathbb{R}^{d_1}$ and $\mathcal{Y}= \mathbb{R}^{d_2}$ (or $\mathcal{Y}\subset \mathbb{R}^{d_2}$). 
The following convergence results build on the common metric $\mathbb{E}\|\nabla F(x)\|$
used in \citep{lin2019gradient,luo2020stochastic}, where $F(x)=\max_{y\in \mathcal{Y}}f(x,y)$.
The related proofs of these convergence analysis are provided in Appendix \ref{Appendix:A6}.

\begin{theorem} \label{th:04}
Suppose the sequence $\{x_t,y_t\}_{t=1}^T$ be generated from Algorithm \ref{alg:4}.
When $\mathcal{X}=\mathbb{R}^{d_1}$, and let $\eta_t = \frac{k}{(m+t)^{1/3}}$
for all $t\geq 0$, $c_1 \geq \frac{2}{3k^3} + \frac{9\tau^2}{4}$
and $c_2 \geq \frac{2}{3k^3} + \frac{75L^2_f}{2}$, $k>0$, $m\geq \max\big(2, k^3, (c_1k)^3, (c_2k)^3\big)$,  $0<\lambda\leq \min\big(\frac{1}{6L_f},
\frac{27b\tau}{16}\big)$ and
$0< \gamma \leq \min\big( \frac{\lambda\tau}{2L_f}\sqrt{\frac{2b}{8\lambda^2 + 75\kappa_y^2b}}, \frac{m^{1/3}}{2L_gk}\big)$,
we have
\begin{align}
 \frac{1}{T} \sum_{t=1}^T \mathbb{E}\|\nabla F(x_t)\|
  \leq \frac{\sqrt{2M''}m^{1/6}}{T^{1/2}} + \frac{\sqrt{2M''}}{T^{1/3}},
\end{align}
where $\Delta_1=\|y_1-y^*(x_1)\|^2$ and $M'' =  \frac{F(x_1) - F^*}{\gamma k} + \frac{9L^2_f\Delta_1}{k \lambda\tau} +\frac{2m^{1/3}\delta^2}{b\tau^2k^2} + \frac{2(c_1^2+c_2^2)\delta^2 k^2}{b\tau^2}\ln(m+T)$.
\end{theorem}
\begin{remark}
Since the conditions of Theorem \ref{th:04} are the same conditions of Theorem \ref{th:4},
Theorem \ref{th:04} has the same results of Theorem \ref{th:4}.
\textbf{When $b=1$}, our Acc-MDA algorithm has a gradient complexity of $4\cdot T=\tilde{O}\big(\kappa_y^{4.5}\epsilon^{-3}\big)$ for finding an $\epsilon$-stationary point;
\textbf{when $b=O(\kappa_y^{\nu})$ for $\nu>0$ and $\frac{27b\tau}{16} \leq \frac{1}{6L_f}$}, i.e., $\kappa_y^{\nu}\leq \frac{8}{81L_f\tau}$,
our Acc-MDA algorithm also reaches a lower gradient complexity of $4b\cdot T=\tilde{O}\big(\kappa_y^{(4.5-\nu/2)}\epsilon^{-3}\big)$ for finding an $\epsilon$-stationary point. When giving $b=O(\kappa_y^3)$, our Acc-MDA reaches the best known gradient complexity of $\tilde{O}\big(\kappa_y^3\epsilon^{-3}\big)$.
When giving $b=O(\kappa_y^4)$, our Acc-MDA reaches a lower gradient complexity of $\tilde{O}\big(\kappa_y^{2.5}\epsilon^{-3}\big)$.
\end{remark}

\begin{remark}
The above low gradient complexities are obtained when $b=O(\kappa_y^{\nu})$ and $\kappa_y^{\nu}\leq \frac{8}{81L_f\tau}$, where $L_f$ denotes the smooth parameter of objective function $f(x,y)$. Without loss of generality, let $\nu=1$, we have $L_f \leq \frac{2\sqrt{2}}{9}$. 
 Although $L_f$ may be large, we can easily change the original objective function $f(x,y)$ into a new function $\hat{f}(x,y)=r f(x,y), \ 0<r< 1$. Since $\nabla \hat{f}(x,y) = r \nabla f(x,y)$,
the gradient of function $\hat{f}(x,y)$ is $\hat{L}$-Lipschitz continuous ($\hat{L}=r L_f$). Thus, we can choose a suitable hyper-parameter $r$ to let this new objective function $\hat{f}(x,y)$ satisfy the condition $\hat{L} \leq \frac{2\sqrt{2}}{9}$.
\end{remark}

\section{ Numerical Experiments}
\label{section:8}
In this section, we evaluate the performance of our algorithms on two applications:
1) black-box adversarial attack to deep neural networks (DNNs) and 2) poisoning attack to logistic regression.
In the first application, we compare our Acc-ZOM algorithm with the ZO-AdaMM \citep{chen2019zo}, ZO-SPIDER-Coord \citep{ji2019improved},
SPIDER-SZO \citep{fang2018spider} and ZO-SFW \citep{sahu2019towards}.
In the second application, for two-side black-box attack, we compare our Acc-ZOMDA algorithm with ZO-Min-Max \citep{liu2019min}
and ZO-SGDMSA \citep{wang2020zeroth} and ZO-SREDA-Boost \citep{xu2020enhanced}.
For one-side black-box attack, we choose ZO-Min-Max~\citep{liu2019min} as a baseline.
For transparent attack, we compare our Acc-MDA algorithm with SGDA~\citep{lin2019gradient} and SREDA-Boost~\citep{xu2020enhanced}.
Note that the SREDA-Boost~\citep{xu2020enhanced} is an improved version of the SREDA algorithm \citep{luo2020stochastic} and
the difference between SREDA-Boost and SREDA is using different learning rate.
In the transparent attack, thus, we only choose the SREDA-Boost as a comparison method.

\subsection{Black-Box Adversarial Attack to DNNs}
In this subsection, we use our Acc-ZOM algorithm to generate adversarial perturbations to attack the pre-trained black-box DNNs,
whose parameters are hidden and only its outputs are accessible. Let $(a, b)$ denote an image $a$ with its true label $b \in\{1,2, \cdots, K\}$,
where $K$ is the total number of image classes. Given multiple images $\{a_{i}, b_{i}\}_{i=1}^{n}$, we design a universal perturbation $x$ to
a pre-trained black-box DNN.
Following \citep{guo2019simple}, we consider the following untargeted attack problem:
\begin{equation}   \label{eq:ap}
    \min_{x \in \mathcal{X}} \frac{1}{n} \sum_{i=1}^{n} \ \max\big( f_{b_i}(x+a_i)-\max_{j\neq b_i} f_j(x+a_i),0 \big),
    \quad \textrm{s.t.} \ \mathcal{X}=\{ \|x\|_\infty \leq \varepsilon\}
\end{equation}
where $f_j(x+a_i)$ represents the output with $j$-th class, that is, the final output before softmax of DNN.
In the experiment, we normalize the pixel values to $[0, 1]^d$, and  use the following smooth form as in \citep{lee2001ssvm} to approximate the above untargeted attack problem:
 \begin{align}
    \min_{x \in \mathcal{X}} & \ \frac{1}{n} \sum_{i=1}^{n} \bigg\{  f_{b_i}(x+a_i)-\max_{j\neq b_i} f_j(x+a_i) + \ln\big(1+\exp\big(\max_{j\neq b_i} f_j(x+a_i)-f_{b_i}(x+a_i)\big)\big) \bigg\}, \nonumber \\
     \textrm{s.t.} & \ \mathcal{X}=\{ \|x\|_\infty \leq \varepsilon\}. \nonumber
\end{align} 

\begin{figure*}[t]%
    \centering
    \subfloat[MNIST]{{\includegraphics[width=3.7cm]{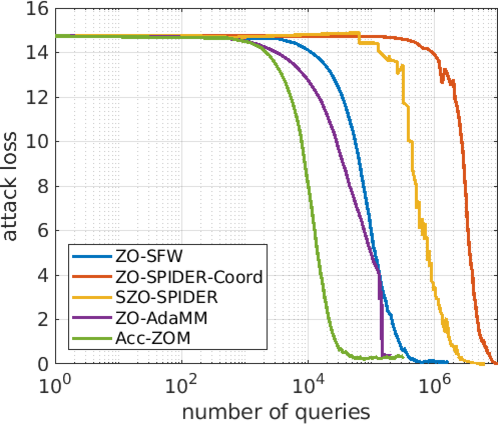} }}%
    \subfloat[FashionMNIST]{{\includegraphics[width=3.7cm]{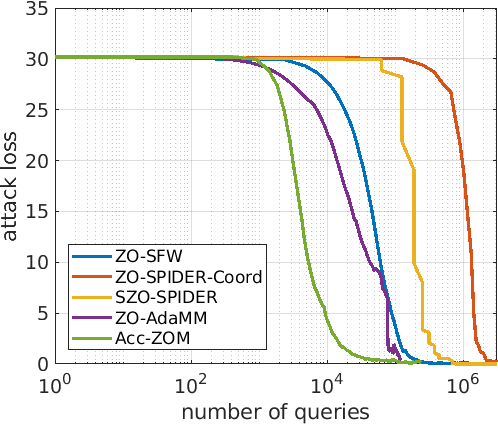} }}
    \subfloat[CIFAR-10]{{\includegraphics[width=3.7cm]{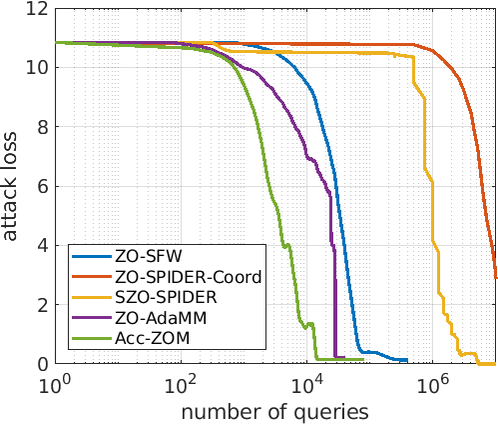} }}%
    \subfloat[SVHN]{{\includegraphics[width=3.7cm]{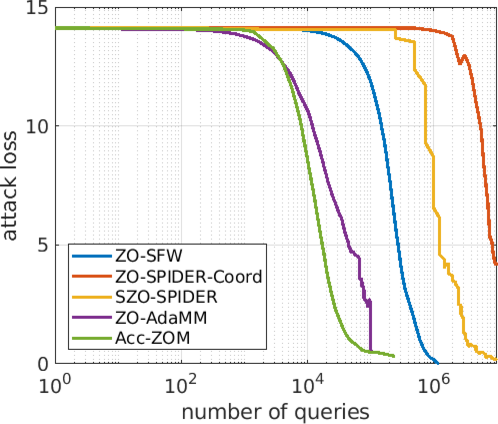} }}%
    \caption{Experimental results of black-box adversarial attack on four datasets: MNIST, FashionMNIST, CIFAR-10 and SVHN.}
    \label{fig:bb-attack}
\end{figure*}

\begin{figure*}[t]%
    \centering
    \subfloat[MNIST]{{\includegraphics[width=3.7cm]{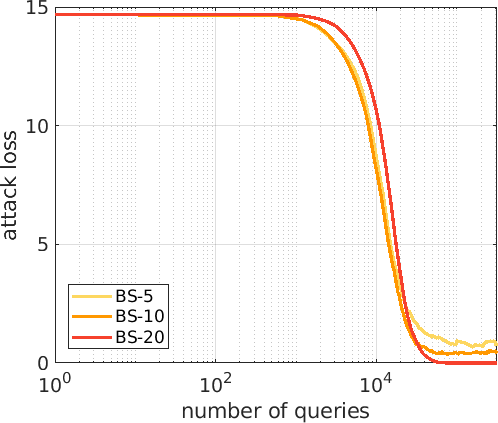} }}%
    \subfloat[FashionMNIST]{{\includegraphics[width=3.7cm]{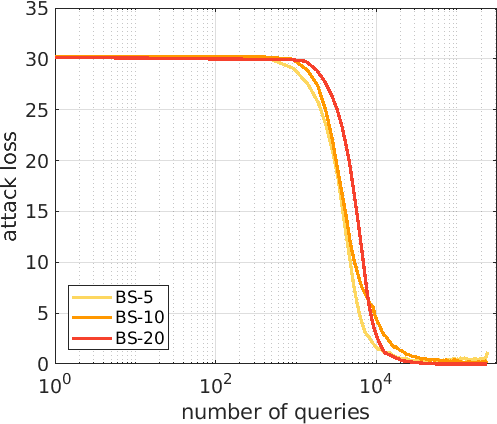} }}%
    \subfloat[CIFAR-10]{{\includegraphics[width=3.7cm]{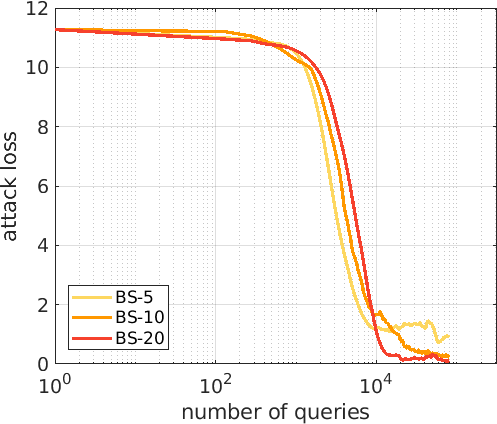} }}%
    \subfloat[SVHN]{{\includegraphics[width=3.7cm]{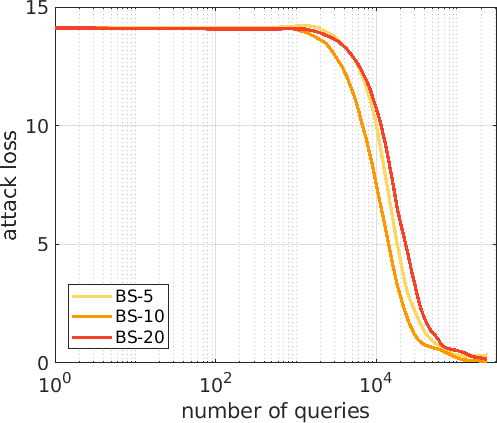} }}%
    \caption{Impact of batch-size on our Acc-ZOM algorithm.}
    \label{fig:bb-attack-bn}
\end{figure*}

In the experiment, we use the pre-trained DNNs on four benchmark datasets: MNIST, FashionMNIST, CIFAR-10, and SVHN,
which attain $99.4\%$, $91.8\%$, $93.2\%$, and $80.8\%$ test accuracy, respectively.
Here, $n$ in problem \eqref{eq:ap} is set to 40 for all datasets.
The batch size of all algorithms is 10. Different datasets require different $\varepsilon$.
Specifically, $\varepsilon$ is set to $0.4$, $0.3$, $0.1$, $0.2$ for MNIST, FashionMNIST, CIFAR-10, and SVHN, respectively.
The hyper-parameters $\gamma, k,m,c$ of the Acc-ZOM are 0.1, 1, 3, 3.
For the other algorithms, we follow the hyper-parameters  in their original paper for a fair comparison.
In Fig.~\ref{fig:bb-attack}, we plot attack loss vs. the number of function queries for each algorithm.
Fig.~\ref{fig:bb-attack} shows that our Acc-ZOM algorithm can largely outperform other algorithms in terms of function queries.
We select hyper-parameters following the theoretic analysis. $k$ is first chosen as 1. Given $k$, $c$ have to be larger than $\frac{2}{3k^3}+\frac{5}{4}$, we then choose $c$ as 3, which is the smallest integer larger than the threshold.
Similarly, $m$ is chosen as 3 to satisfy the condition $m\leq \max((ck)^3,k^3)$.
To study the impact of batch-size, we use three different batch-size settings: 5, 10, 20.
From Fig.~\ref{fig:bb-attack-bn}, we can see that our Acc-ZOM algorithm can work well on a range of batch-size selections.

\subsection{Poisoning Attack to Logistic Regression}
In this subsection, we apply the task of poisoning attack to logistic regression to demonstrate the efficiency of our
Acc-ZOMDA, Acc-Semi-ZOMDA and Acc-MDA.
Let $\{a_i,b_i\}_{i=1}^n$ denote the training dataset, in which $n_0\ll n$ samples are corrupted by a perturbation vector $x$. Following \cite{liu2019min},
this poisoning attack problem is formulated as
\begin{align} \label{eq:23}
 \max_{x\in \mathcal{X}}\min_{y\in \mathcal{Y}} & \ f(x,y)= h(x,y;\mathcal{D}_{p}) + h(0,y;\mathcal{D}_{t}), \\
 \mbox{s.t.}  & \ \mathcal{X}=\{\|x\|_\infty \leq \varepsilon\}, \
   \mathcal{Y}=\{\|y\|^2_2 \leq \lambda_{\textrm{reg}}\} \nonumber
\end{align}
where $\mathcal{D}_{p}$ and $\mathcal{D}_{t}$ are corrupted set and clean set respectively, $y$ is the model parameter,
the corrupted rate $\frac{|\mathcal{D}_{p}|}{|\mathcal{D}_{t}|+|\mathcal{D}_{p}|}$ is set to 0.15.
Here $h(x,y;\mathcal{D})=-\frac{1}{|\mathcal{D}|}\sum_{(a_i,b_i)\in \mathcal{D}} \big[ b_i\log(g(x,y;a_i)) + (1-b_i)\log(1-g(x,y;a_i)) \big]$
with $g(x,y;a_i) = \frac{1}{1+e^{-(x+a_i)^Ty}}$. Note that the above problem \eqref{eq:23} can be written in the form
of \eqref{eq:2}, i.e., $\min_{x\in \mathcal{X}}\max_{y\in\mathcal{Y}} \big\{-f(x,y)\big\}$. 
In the experiment, we generate $n=1000$ samples.
Specifically, we randomly draw the feature vector $a_i\in \mathbb{R}^{100}$ from normal distribution $\mathcal{N}(0,1)$,
and label $b_i=1$ if $\frac{1}{1+e^{-(a_i^T\theta+\nu_i)}}>\frac{1}{2}$, otherwise $b_i=0$.
Here we choose $\theta=(1,1,\cdots,1)$ as the ground-truth model parameters, and $\nu_i\in \mathcal{N}(0,10^{-3})$.
For this experiment, we set $\varepsilon$ and $\lambda_{\textrm{reg}}$ to $2$ and $0.001$.
We also chose the hyper-parameters $\gamma,\lambda,k,m,c_1,c_2$ of our Acc-ZOMDA as $0.2, 0.08, 1, 3, 3, 3$.

\begin{figure*}[t]%
    \centering
    \subfloat[Two-Side Black-Box Attack]{{\includegraphics[width=0.325\textwidth]{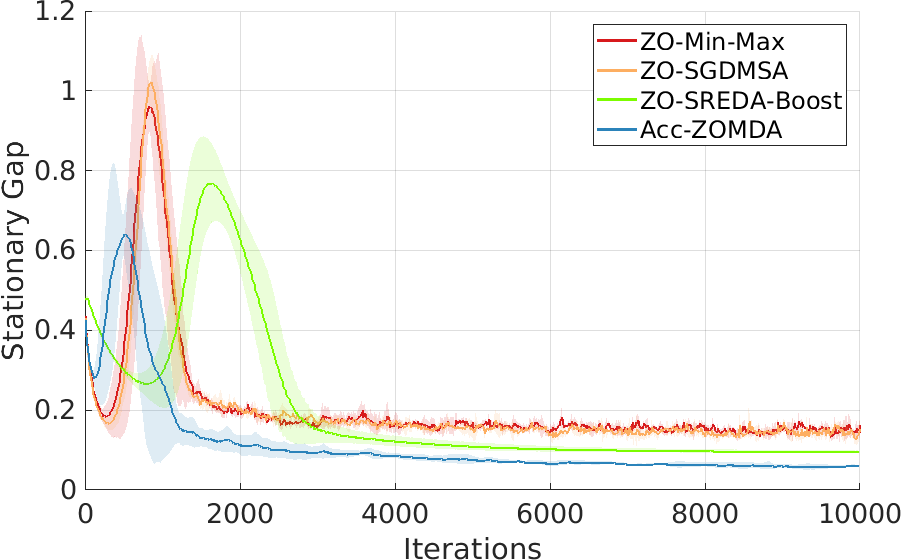} }}%
    \subfloat[One-Side Black-Box Attack]{{\includegraphics[width=0.325\textwidth]{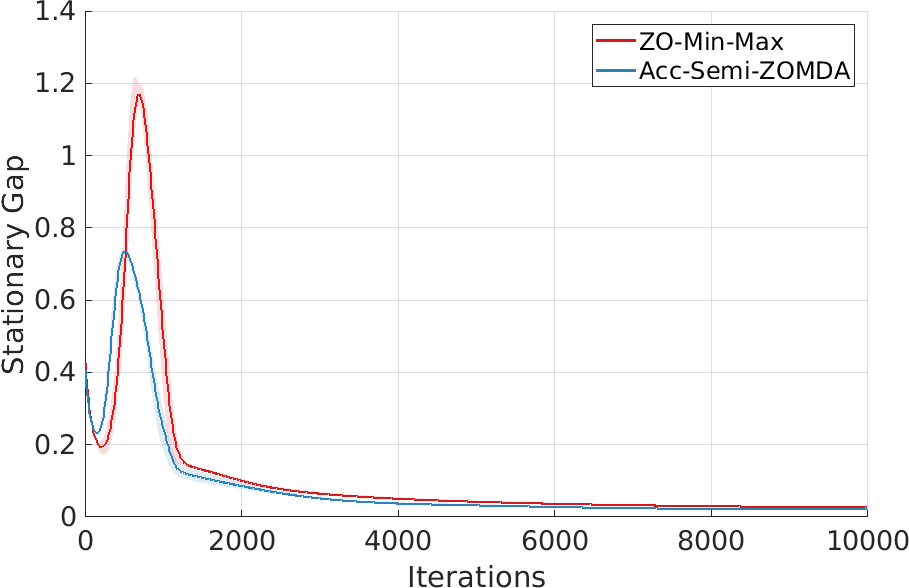} }}%
    \subfloat[Transparent Attack]{{\includegraphics[width=0.325\textwidth]{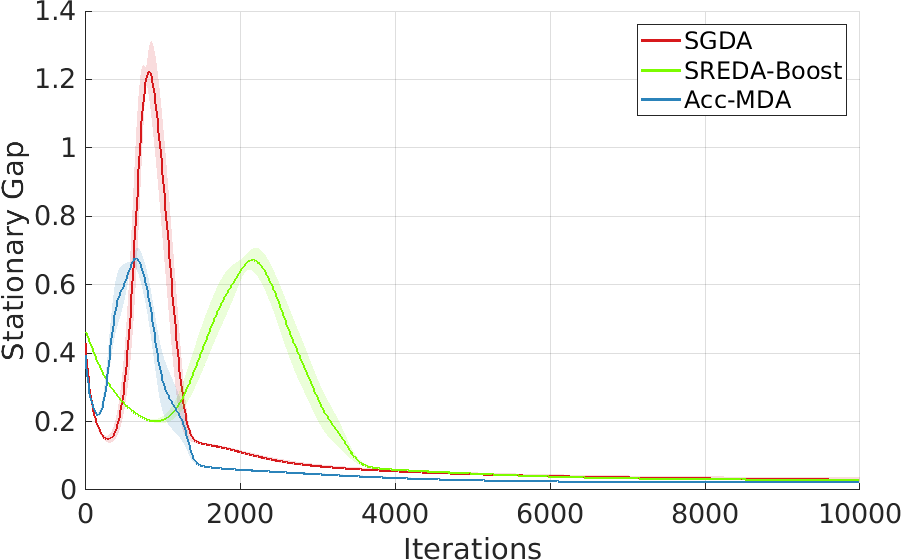} }}%
    \caption{Stationary gap of different methods in two-side black-box scenario, one-side black-box scenario and transparent scenario.}
    \label{fig:pa}
\end{figure*}

\begin{figure*}[!t]%
    \centering
    \subfloat[Two-Side Black-Box Attack]{{\includegraphics[width=0.325\textwidth]{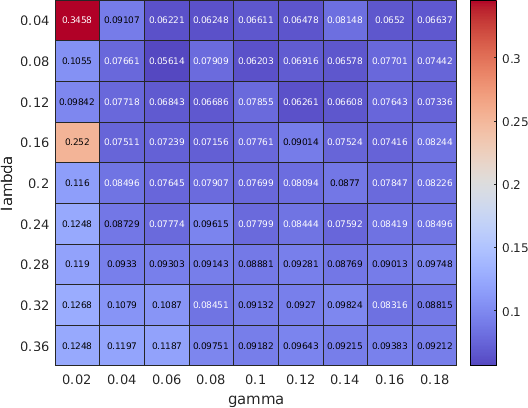} }}%
     \subfloat[One-Side Black-Box Attack]{{\includegraphics[width=0.325\textwidth]{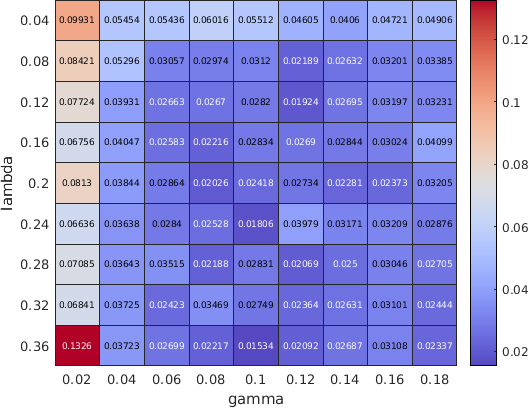} }}%
    \subfloat[Transparent Attack]{{\includegraphics[width=0.325\textwidth]{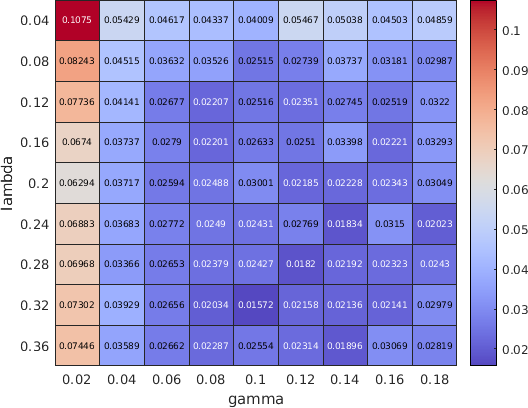} }}%
    \caption{Stationary gap given different combinations of tuning parameters $(\gamma,\lambda)$.}
    \label{fig:hm}
\end{figure*}

From Fig.~\ref{fig:pa}(a), we can find that our Acc-ZOMDA algorithm converges fastest and achieves lowest stationary gap.
The Acc-ZOMDA is also robust to different learning rate pairs of $(\gamma, \lambda)$.
In Fig.~\ref{fig:pa}(b,c), we show the comparison results for one-side black-box (black-box w.r.t attacker)
poison attack and transparent poison attack.
All hyper-parameter settings are the same as two-side black-box attack.
These results demonstrate that our Acc-Semi-ZOMDA and Acc-MDA algorithms compare favorably with other algorithms.

To better understanding the settings of hyper-parameters, we visualize the stationary gap given different combinations of $(\gamma, \lambda)$.
We set $\gamma$ from 0.04 to 0.036 and $\lambda$ from 0.02 to 0.18. 
From Fig.~\ref{fig:hm}, we can see that our method can achieve ideal stationary gap with most combinations of $(\gamma, \lambda)$
across three different scenarios.

\section{Conclusions}
\label{section:9}
In the paper, we proposed a class of accelerated zeroth-order and first-order momentum methods for both nonconvex mini-optimization and minimax-optimization, which build on the momentum-based variance reduced technique of STORM and momentum update. Moreover, we gave an effective convergence analysis framework for our methods. Specifically, we proved that our zeroth-order methods can obtain a low query complexity without requiring any large bathes. Meanwhile, our first-order method also can obtain a low gradient complexity without requiring any large bathes. In particular, 
our methods are the first to extend the STORM algorithm to constrained optimization and minimax optimization. 

\acks{ We thank editor and three anonymous reviewers for their valuable comments. This work was partially supported by NSF IIS 1845666, 1852606, 1838627, 1837956, 1956002, OIA 2040588. }


%

%
%
%
%
%
%
%



\appendix
\section{ Detailed Convergence Analysis }
In this section, we provide the detailed convergence analysis of our algorithms.
We first review some useful lemmas.
\begin{lemma} \label{lem:1}
 \citep{lin2019gradient} Under the above Assumptions \ref{ass:5} and \ref{ass:6}, 
 the function $F(x)=\max_{y\in \mathcal{Y}} f(x,y)$ has $L_g$-Lipschitz continuous gradient, such as
 \begin{align}
 \|\nabla F(x) - \nabla F(x')\| \leq L_g\|x-x'\|, \ \forall x, x' \in \mathcal{X}
 \end{align}
 where $L_g=L_f+\frac{L^2_f}{\tau}$.
\end{lemma}
\begin{lemma} \label{lem:2}
 \citep{lin2019gradient} Under the above Assumptions \ref{ass:5} and \ref{ass:6},
 the mapping $y^*(x)=\arg\max_{y\in \mathcal{Y}}f(x,y)$ is $\kappa_y$-Lipschitz continuous,
 such as
 \begin{align}
  \|y^*(x) - y^*(x')\| \leq \kappa_y \|x-x'\|, \ \forall x,x' \in \mathcal{X}
 \end{align}
 where $\kappa_y = L_f/\tau$ denotes the condition number for function $f(\cdot,y)$.
\end{lemma}
\begin{lemma} \label{lem:A1}
\citep{nesterov2018lectures} Assume that $f(x)$ is a differentiable convex function and $\mathcal{X}$ is a convex set.
 $x^* \in \mathcal{X}$ is the solution of the
constrained problem $\min_{x\in \mathcal{X}}f(x)$, if
\begin{align}
 \langle \nabla f(x^*), x-x^*\rangle \geq 0, \ \forall  x\in \mathcal{X}.
\end{align}
\end{lemma}
\begin{lemma} \label{lem:A2}
\citep{nesterov2018lectures} Assume that the function $f(x)$ is $L$-smooth, i.e., $\|\nabla f(x)-\nabla f(y)\|\leq L\|x-y\|$,
 the following inequality satisfies
\begin{align}
 |f(y)-f(x)-\nabla f(x)^T(y-x)|\leq \frac{L}{2}\|x-y\|^2.
\end{align}
\end{lemma}
\begin{lemma} \label{lem:A3}
 \citep{gao2018information,ji2019improved} Let $f_{\mu}(x)=\mathbb{E}_{u\sim U_B}[f(x+\mu u)]$ be a smooth approximation of function $f(x)$,
where $U_B$ is the uniform distribution over the $d$-dimensional unit Euclidean ball $B$. Given zeroth-order gradient
$\hat{\nabla} f(x) = \frac{ f(x + \mu u)-f(x)}{\mu/d}u$, we have
\begin{itemize}
\item[(1)] If $f(x)$ has $L$-Lipschitz continuous gradient (i.e., $L$-smooth), then $f_{\mu}(x)$ has $L$-Lipschitz continuous gradient; 
\item[(2)] $|f_{\mu}(x) - f(x)| \leq \frac{\mu^2 L}{2}$ and $\|\nabla f_{\mu}(x) - \nabla f(x)\| \leq \frac{\mu Ld}{2}$ for any $x\in \mathbb{R}^d$;
\item[(3)] $\mathbb{E}[\frac{1}{|\mathcal{S}|}\sum_{i\in \mathcal{S}}\hat{\nabla}f(x;\xi_i)] = \nabla f_{\mu}(x)$ for any $x\in \mathbb{R}^d$;
\item[(4)] $\mathbb{E}\|\hat{\nabla}f(x;\xi)-\hat{\nabla}f(x';\xi)\|^2 \leq 3dL^2\|x-x'\|^2 + \frac{3L^2d^2\mu}{2}$ for any $x, x' \in \mathbb{R}^d$.
\end{itemize}
\end{lemma}
\begin{lemma} \label{lem:A4}
For i.i.d. random variables $\{\xi_i\}_{i=1}^n$ with zero mean, we have
$\mathbb{E}\|\frac{1}{n}\sum_{i=1}^n\xi_i\|^2 = \frac{1}{n}\mathbb{E}\|\xi_i\|^2$ for any $i\in [n]$.
\end{lemma}
Note that the above results (1)-(2) of Lemma \ref{lem:A3} come from Lemma 4.1 in \citep{gao2018information}, and
the above results (3)-(4) come from Lemma 5 in \citep{ji2019improved}. In addition, the result (4) of Lemma \ref{lem:A3}
is an extended result from Lemma 5 in \citep{ji2019improved}.

\subsection{ Convergence Analysis of Acc-ZOM Algorithm for Constrained Mini-Optimization }
\label{Appendix:A1}
In this subsection, we study the convergence properties of our Acc-ZOM algorithm for solving the black-box \textbf{constrained}
 problem \eqref{eq:1},i.e., $\mathcal{X}\subset \mathbb{R}^d$. 
We first let $f_{\mu}(x)=\mathbb{E}_{u\sim U_B}[f(x+\mu u)]$
be a smooth approximation of function $f(x)$, where $U_B$ is the uniform distribution over the $d$-dimensional unit Euclidean ball $B$.

\begin{lemma} \label{lem:B1}
 Suppose that the sequence $\{x_t\}_{t=1}^T$ be generated from Algorithm \ref{alg:1}.
 Let $0<\eta_t \leq 1$ and $0< \gamma \leq \frac{1}{2L\eta_t}$,
 then we have
 \begin{align}
  f_{\mu}(x_{t+1}) - f_{\mu}(x_t) \leq \eta_t\gamma\|\nabla f_{\mu}(x_t)-v_t\|^2-\frac{\eta_t}{2\gamma}\|\tilde{x}_{t+1}-x_t\|^2.
 \end{align}
\end{lemma}
\begin{proof}
 According to Assumption \ref{ass:3} and Lemma \ref{lem:A3}, the function $f_{\mu}(x)$ is $L$-smooth.
 Then we have
 \begin{align} \label{eq:B1}
  f_{\mu}(x_{t+1}) & \leq f_{\mu}(x_t) + \langle\nabla f_{\mu}(x_t), x_{t+1}-x_t\rangle + \frac{L}{2}\|x_{t+1}-x_t\|^2  \\
  & = f_{\mu}(x_t) + \eta_t\langle \nabla f_{\mu}(x_t),\tilde{x}_{t+1}-x_t\rangle + \frac{L\eta_t^2}{2}\|\tilde{x}_{t+1}-x_t\|^2 \nonumber \\
  & = f_{\mu}(x_t) + \eta_t\langle \nabla f_{\mu}(x_t)-v_t,\tilde{x}_{t+1}-x_t\rangle + \eta_t\langle v_t,\tilde{x}_{t+1}-x_t\rangle + \frac{L\eta_t^2}{2}\|\tilde{x}_{t+1}-x_t\|^2,\nonumber
 \end{align}
 where the second equality is due to $x_{t+1}=x_t + \eta_t(\tilde{x}_{t+1}-x_t)$.
 By the step 8 of Algorithm \ref{alg:1}, we have $\tilde{x}_{t+1}=\mathcal{P}_{\mathcal{X}}(x_t - \gamma v_t)
 =\arg\min_{x\in\mathcal{X}}\frac{1}{2}\|x- x_t + \gamma v_t\|^2$. Since $\mathcal{X}$ is a convex set
 and the function $\frac{1}{2}\|x- x_t + \gamma v_t\|^2$ is convex, by using Lemma \ref{lem:A1},
 we have
 \begin{align} \label{eq:B2}
  \langle \tilde{x}_{t+1}- x_t + \gamma v_t, x-\tilde{x}_{t+1}\rangle \geq 0, \ \forall x\in \mathcal{X}.
 \end{align}
 In Algorithm \ref{alg:1}, let the initialize solution $x_1 \in \mathcal{X}$, and the sequence $\{x_t\}_{t\geq 1}$ generates as follows:
 \begin{align}
  x_{t+1} = x_{t} + \eta_t(\tilde{x}_{t+1} - x_t) = \eta_t\tilde{x}_{t+1} + (1-\eta_t)x_t,
 \end{align}
 where $0<\eta_t \leq 1$.
 Since $\mathcal{X}$ is convex set and $x_t, \ \tilde{x}_{t+1} \in \mathcal{X}$, we have $x_{t+1} \in \mathcal{X}$
 for any $t \geq 1$.
 Set $x=x_t$ in the inequality \eqref{eq:B2}, we have
 \begin{align} \label{eq:B3}
 \langle v_t, \tilde{x}_{t+1}-x_t\rangle \leq -\frac{1}{\gamma}\|\tilde{x}_{t+1}-x_t\|^2.
 \end{align}
 By using Cauchy-Schwarz inequality and Young's inequality, we have
 \begin{align} \label{eq:B4}
  \langle \nabla f_{\mu}(x_t)-v_t,\tilde{x}_{t+1}-x_t\rangle & \leq \|\nabla f_{\mu}(x_t)-v_t\|\cdot\|\tilde{x}_{t+1}-x_t\| \nonumber \\
  & \leq \gamma\|\nabla f_{\mu}(x_t)-v_t\|^2+\frac{1}{4\gamma}\|\tilde{x}_{t+1}-x_t\|^2.
 \end{align}
 Combining the inequalities \eqref{eq:B1}, \eqref{eq:B3} with \eqref{eq:B4},
 we obtain
 \begin{align}
 & f_{\mu}(x_{t+1}) \nonumber \\
 & \leq f_{\mu}(x_t) + \eta_t\langle \nabla f_{\mu}(x_t)-v_t,\tilde{x}_{t+1}-x_t\rangle + \eta_t\langle v_t,\tilde{x}_{t+1}-x_t\rangle + \frac{L\eta_t^2}{2}\|\tilde{x}_{t+1}-x_t\|^2 \nonumber \\
  & \leq f_{\mu}(x_t) + \eta_t\gamma\|\nabla f_{\mu}(x_t)-v_t\|^2+\frac{\eta_t}{4\gamma}\|\tilde{x}_{t+1}-x_t\|^2 -\frac{\eta_t}{\gamma}\|\tilde{x}_{t+1}-x_t\|^2  + \frac{L\eta_t^2}{2}\|\tilde{x}_{t+1}-x_t\|^2 \nonumber \\
  & = f_{\mu}(x_t) + \eta_t\gamma\|\nabla f_{\mu}(x_t)-v_t\|^2 -\frac{\eta_t}{2\gamma}\|\tilde{x}_{t+1}-x_t\|^2 -\big(\frac{\eta_t}{4\gamma}-\frac{L\eta_t^2}{2}\big)\|\tilde{x}_{t+1}-x_t\|^2 \nonumber \\
  & \leq f_{\mu}(x_t) + \eta_t\gamma\|\nabla f_{\mu}(x_t)-v_t\|^2 -\frac{\eta_t}{2\gamma}\|\tilde{x}_{t+1}-x_t\|^2,
 \end{align}
 where the last inequality is due to $0< \gamma \leq \frac{1}{2L\eta_t}$.
\end{proof}

\begin{lemma} \label{lem:C1}
 Suppose the zeroth-order stochastic gradient $\{v_t\}$ be generated from Algorithm \ref{alg:1}, we have
 \begin{align}
 \mathbb{E} \|\nabla f_{\mu}(x_{t+1}) - v_{t+1}\|^2 & \leq (1-\alpha_{t+1})^2 \mathbb{E} \|\nabla f_{\mu}(x_t) -v_t\|^2 + 6(1-\alpha_{t+1})^2 dL^2 \eta_t^2\mathbb{E} \| \tilde{x}_{t+1} - x_t\|^2 \nonumber \\
 & \quad + 3(1-\alpha_{t+1})^2L^2d^2\mu^2 + 2\alpha_{t+1}^2\sigma^2.
 \end{align}
\end{lemma}
\begin{proof}
 According to the definition of $v_{t+1}$ in Algorithm \ref{alg:1}, we have
 \begin{align}
  v_{t+1}-v_t = -\alpha_{t+1}v_t + (1-\alpha_{t+1})\big(\hat{\nabla} f(x_{t+1};\xi_{t+1}) - \hat{\nabla} f(x_t;\xi_{t+1})\big)
  + \alpha_{t+1}\hat{\nabla} f(x_{t+1};\xi_{t+1}). \nonumber
 \end{align}
 Then we have
 \begin{align}
  & \mathbb{E} \|\nabla f_{\mu}(x_{t+1}) - v_{t+1}\|^2 \nonumber \\
  & = \mathbb{E} \|\nabla f_{\mu}(x_{t+1}) - v_t - (v_{t+1}-v_t)\|^2 \nonumber \\
  & = \mathbb{E} \|\nabla f_{\mu}(x_{t+1}) - v_t + \alpha_{t+1}v_t- \alpha_{t+1}\hat{\nabla} f(x_{t+1};\xi_{t+1}) \nonumber \\
  &\quad- (1-\alpha_{t+1})(\hat{\nabla} f(x_{t+1};\xi_{t+1})
  - \hat{\nabla} f(x_t;\xi_{t+1})) \|^2 \nonumber \\
  & = \mathbb{E} \|(1-\alpha_{t+1})(\nabla f_{\mu}(x_t) -v_t) + \alpha_{t+1}\big(\nabla f_{\mu}(x_{t+1})- \hat{\nabla} f(x_{t+1};\xi_{t+1})\big)\nonumber \\
  & \quad + (1-\alpha_{t+1})\big(\nabla f_{\mu}(x_{t+1})-\nabla f_{\mu}(x_t)-\hat{\nabla} f(x_{t+1};\xi_{t+1}) + \hat{\nabla} f(x_t;\xi_{t+1})\big) \|^2 \nonumber \\
  & = (1-\alpha_{t+1})^2\mathbb{E}\|\nabla f_{\mu}(x_t)-v_t\|^2+ \|\alpha_{t+1}\big(\nabla f_{\mu}(x_{t+1})- \hat{\nabla} f(x_{t+1};\xi_{t+1})\big)\nonumber \\
  & \quad + (1-\alpha_{t+1})\big(\nabla f_{\mu}(x_{t+1})-\nabla f_{\mu}(x_t)-\hat{\nabla} f(x_{t+1};\xi_{t+1}) + \hat{\nabla} f(x_t;\xi_{t+1})\big) \|^2 \nonumber \nonumber \\
  & \leq (1-\alpha_{t+1})^2\mathbb{E} \|\nabla f_{\mu}(x_t) -v_t\|^2 + 2(1-\alpha_{t+1})^2\mathbb{E} \|\nabla f_{\mu}(x_{t+1}) - \nabla f_{\mu}(x_t) - \hat{\nabla} f(x_{t+1};\xi_{t+1})  \nonumber \\
  & \quad + \hat{\nabla} f(x_t;\xi_{t+1})\|^2+ 2\alpha_{t+1}^2 \mathbb{E} \|\nabla f_{\mu}(x_{t+1})- \hat{\nabla} f(x_{t+1};\xi_{t+1})\|^2\nonumber \\
  & \leq (1-\alpha_{t+1})^2 \mathbb{E} \|\nabla f_{\mu}(x_t) \!-\!v_t\|^2 \!+\! 2(1-\alpha_{t+1})^2 \mathbb{E} \| \hat{\nabla} f(x_{t+1};\xi_{t+1}) \!-\! \hat{\nabla} f(x_t;\xi_{t+1})\|^2 \!+\! 2\alpha_{t+1}^2\sigma^2 \nonumber \\
  & \leq (1-\alpha_{t+1})^2 \mathbb{E} \|\nabla f_{\mu}(x_t) -v_t\|^2 + 6(1-\alpha_{t+1})^2 dL^2 \mathbb{E} \| x_{t+1} - x_t\|^2 \nonumber \\
  &\quad + 3(1-\alpha_{t+1})^2L^2d^2\mu^2 + 2\alpha_{t+1}^2\sigma^2 \nonumber \\
  & = (1-\alpha_{t+1})^2 \mathbb{E} \|\nabla f_{\mu}(x_t) -v_t\|^2 + 6(1-\alpha_{t+1})^2 dL^2 \eta_t^2\mathbb{E} \| \tilde{x}_{t+1} - x_t\|^2 \nonumber \\
  &\quad + 3(1-\alpha_{t+1})^2L^2d^2\mu^2 + 2\alpha_{t+1}^2\sigma^2,
 \end{align}
 where the fourth equality follows by $\mathbb{E}_{(u,\xi)}[\hat{\nabla} f(x_{t+1};\xi_{t+1})]=\nabla f_{\mu}(x_{t+1})$ and \\$ \mathbb{E}_{(u,\xi)}[\hat{\nabla} f(x_{t+1};\xi_{t+1}) 
 - \hat{\nabla} f(x_t;\xi_{t+1})]=\nabla f_{\mu}(x_{t+1}) - \nabla f_{\mu}(x_t)$; the first inequality holds by Cauchy-Schwarz inequality; the second inequality holds by the equality $\mathbb{E}\|\zeta-\mathbb{E}[\zeta]\|^2 =\mathbb{E}\|\zeta\|^2 - \|\mathbb{E}[\zeta]\|^2$ and Assumption \ref{ass:1}, and the last inequality holds by Young's inequality and Lemma \ref{lem:A3}.
 \end{proof}

 \begin{theorem} \label{th:A1}
 (Restatement of Theorem 1)
 Suppose the sequence $\{x_t\}_{t=1}^T$ be generated from Algorithm \ref{alg:1}. When $\mathcal{X}\subset \mathbb{R}^d$, and let $\eta_t = \frac{k}{(m+t)^{1/3}}$
 for all $t\geq 0$, $0< \gamma \leq \min\big(\frac{m^{1/3}}{2Lk},\frac{1}{2\sqrt{6d}L}\big)$,
 $c\geq \frac{2}{3k^3} + \frac{5}{4}$, $k>0$, $m\geq \max\big(2, k^3 ,(ck)^3\big)$ and $0<\mu \leq \frac{1}{d(m+T)^{2/3}}$,
 we have
 \begin{align}
  \frac{1}{T}\sum_{t=1}^T\mathbb{E} \|G_{\mathcal{X}}(x_t,\nabla f(x_t),\gamma)\| & \leq \frac{1}{T} \sum_{t=1}^T\mathbb{E} \big[ \|\nabla f(x_{t}) - v_t\| + \frac{1}{\gamma} \|\tilde{x}_{t+1}-x_t\| \big] \nonumber \\
  & \leq  \frac{2\sqrt{2M}m^{1/6}}{T^{1/2}} + \frac{2\sqrt{2M}}{T^{1/3}} + \frac{L}{2(m+T)^{2/3}},
\end{align}
where $M=\frac{f_{\mu}(x_1) - f^*}{k\gamma} + \frac{m^{1/3}\sigma^2}{k^2} + \frac{9L^2}{4k^2} + 2k^2c^2\sigma^2\ln(m+T)$.
 \end{theorem}
 \begin{proof}
 Since $\eta_t=\frac{k}{(m+t)^{1/3}}$ on $t$ is decreasing and $m\geq k^3$, we have $\eta_t \leq \eta_0 = \frac{k}{m^{1/3}} \leq 1$ and $\gamma \leq \frac{m^{1/3}}{2Lk}=\frac{1}{2L\eta_0} \leq \frac{1}{2L\eta_t}$ for any $t\geq 0$.
 Due to $0 < \eta_t \leq 1$ and $m\geq (ck)^3$, we have $\alpha_{t+1} = c\eta_t^2 \leq c\eta_t \leq \frac{ck}{m^{1/3}}\leq 1$.
 According to Lemma \ref{lem:C1}, we have
 \begin{align}
  & \frac{1}{\eta_t}\mathbb{E} \|\nabla f_{\mu}(x_{t+1}) - v_{t+1}\|^2 -  \frac{1}{\eta_{t-1}}\mathbb{E} \|\nabla f_{\mu}(x_{t}) - v_t\|^2 \nonumber \\
  & \leq \big(\frac{(1-\alpha_{t+1})^2}{\eta_t} - \frac{1}{\eta_{t-1}}\big)\mathbb{E} \|\nabla f_{\mu}(x_{t}) -v_t\|^2 + 6(1-\alpha_{t+1})^2 dL^2\eta_t \mathbb{E} \| \tilde{x}_{t+1} - x_{t}\|^2 \nonumber \\
  & \quad  + \frac{ 3(1-\alpha_{t+1})^2L^2d^2\mu^2}{\eta_t}+ \frac{2\alpha_{t+1}^2\sigma^2}{\eta_t} \nonumber \\
  & \leq \big(\frac{1\!-\!\alpha_{t+1}}{\eta_t} \!-\! \frac{1}{\eta_{t-1}}\big)\mathbb{E} \|\nabla f_{\mu}(x_{t}) -v_t\|^2 \!+\! 6 dL^2\eta_t \mathbb{E} \| \tilde{x}_{t+1} - x_{t}\|^2 \!+\! \frac{ 3L^2d^2\mu^2}{\eta_t} \!+\! \frac{2\alpha_{t+1}^2\sigma^2}{\eta_t} \nonumber \\
  & = \big(\frac{1}{\eta_t} \!-\! \frac{1}{\eta_{t-1}} - c\eta_t\big)\mathbb{E} \|\nabla f_{\mu}(x_{t}) -v_t\|^2 \!+\! 6 dL^2\eta_t \mathbb{E} \| \tilde{x}_{t+1} - x_{t}\|^2 \!+\! \frac{ 3L^2d^2\mu^2}{\eta_t} \!+\! \frac{2\alpha_{t+1}^2\sigma^2}{\eta_t},
 \end{align}
 where the second inequality is due to $0<\alpha_{t+1}\leq 1$.
 By $\eta_t = \frac{k}{(m+t)^{1/3}}$, we have
 \begin{align}
  \frac{1}{\eta_t} - \frac{1}{\eta_{t-1}} &= \frac{1}{k}\big( (m+t)^{\frac{1}{3}} - (m+t-1)^{\frac{1}{3}}\big) \nonumber \\
  & \leq \frac{1}{3k(m+t-1)^{2/3}} \leq  \frac{1}{3k\big(m/2+t\big)^{2/3}} \nonumber \\
  & \leq \frac{2^{2/3}}{3k(m+t)^{2/3}} = \frac{2^{2/3}}{3k^3}\frac{k^2}{(m+t)^{2/3}}= \frac{2^{2/3}}{3k^3}\eta_t^2 \leq \frac{2}{3k^3}\eta_t,
 \end{align}
 where the first inequality holds by the concavity of function $f(x)=x^{1/3}$, \emph{i.e.}, $(x+y)^{1/3}\leq x^{1/3} + \frac{y}{3x^{2/3}}$; the second inequality is due to $m\geq 2$, and
 the last inequality is due to $0<\eta_t\leq 1$. Let $c \geq \frac{2}{3k^3} + \frac{5}{4} $, we have
 \begin{align}
  & \frac{1}{\eta_t}\mathbb{E} \|\nabla f_{\mu}(x_{t+1}) - v_{t+1}\|^2 - \frac{1}{\eta_{t-1}}\mathbb{E} \|\nabla f_{\mu}(x_{t}) - v_t\|^2 \nonumber \\
  & \leq -\frac{5\eta_t}{4}\mathbb{E} \|\nabla f_{\mu}(x_{t}) -v_t\|^2 \!+\! 6 dL^2\eta_t \mathbb{E} \| \tilde{x}_{t+1} - x_{t}\|^2 \!+\! \frac{ 3L^2d^2\mu^2}{\eta_t} \!+\! \frac{2\alpha_{t+1}^2\sigma^2}{\eta_t}.
 \end{align}

 Next, we define a \emph{Lyapunov} function $R_t = \mathbb{E}\big[f_{\mu}(x_t) + \frac{\gamma}{\eta_{t-1}}\|\nabla f_{\mu}(x_t)-v_t\|^2\big]$ for any $t\geq 1$.
 According to Lemma \ref{lem:B1}, we have
 \begin{align}
  R_{t+1} - R_t & = \mathbb{E}\big[f_{\mu}(x_{t+1}) - f_{\mu}(x_t)\big] + \frac{\gamma}{\eta_{t}}\mathbb{E}\|\nabla f_{\mu}(x_{t+1})-v_{t+1}\|^2  - \frac{\gamma}{\eta_{t-1}}\mathbb{E} \|\nabla f_{\mu}(x_t)-v_t\|^2 \nonumber \\
  & \leq \eta_t\gamma \mathbb{E} \|\nabla f_{\mu}(x_t)-v_t\|^2 - \frac{\eta_t}{2\gamma}\mathbb{E}\|\tilde{x}_{t+1}-x_t\|^2-\frac{5\gamma\eta_t}{4}\mathbb{E} \|\nabla f_{\mu}(x_{t}) -v_t\|^2\nonumber \\
  & \quad + 6 dL^2\eta_t\gamma \mathbb{E} \| \tilde{x}_{t+1} - x_{t}\|^2+ \frac{ 3L^2d^2\mu^2\gamma}{\eta_t} + \frac{2\alpha_{t+1}^2\sigma^2\gamma}{\eta_t} \nonumber \\
  & \leq -\frac{\gamma\eta_t}{4}\mathbb{E} \|\nabla f_{\mu}(x_{t}) -v_t\|^2-\frac{\eta_t}{4\gamma}\mathbb{E} \|\tilde{x}_{t+1}-x_t\|^2 + \frac{ 3L^2d^2\mu^2\gamma}{\eta_t} + \frac{2\alpha_{t+1}^2\sigma^2\gamma}{\eta_t},
 \end{align}
where the last inequality is due to $\gamma \leq \frac{1}{2\sqrt{6d}L}$.
Thus, we obtain
\begin{align} \label{eq:C3}
\frac{\gamma \eta_t}{4}\mathbb{E} \|\nabla f_{\mu}(x_{t}) -v_t\|^2 + \frac{\eta_t}{4\gamma}\mathbb{E}\|\tilde{x}_{t+1}-x_t\|^2 \leq R_t - R_{t+1} + \frac{ 3L^2d^2\mu^2\gamma}{\eta_t} + \frac{2\alpha_{t+1}^2\sigma^2\gamma}{\eta_t}.
\end{align}
Since $\inf_{x\in \mathcal{X}} f(x) =f^*$, we have
$\inf_{x\in \mathcal{X}} f_{\mu}(x) =\inf_{x\in \mathcal{X}} \mathbb{E}_{u\sim U_B}[f(x+\mu u)] =
\inf_{x\in \mathcal{X}} \frac{1}{V}\int_{B}f(x+\mu u)du \geq \frac{1}{V}\int_{B}\inf_{x\in \mathcal{X}}f(x+\mu u)du=f^*$,
where $V$ denotes the volume of the unit ball $B$.

Taking average over $t=1,2,\cdots,T$ on both sides of \eqref{eq:C3}, we have
\begin{align}
 & \frac{1}{T} \sum_{t=1}^T\mathbb{E} \big[\frac{\gamma \eta_t}{4} \|\nabla f_{\mu}(x_{t}) -v_t\|^2 + \frac{\eta_t}{4\gamma}\|\tilde{x}_{t+1}-x_t\|^2 \big] \nonumber \\
 & \leq \frac{f_{\mu}(x_1) - f^*}{T} + \frac{\gamma\|\nabla f_{\mu}(x_1) - v_1\|^2}{T\eta_0}  + \sum_{t=1}^T\frac{3L^2d^2\mu^2\gamma}{T\eta_t} + \sum_{t=1}^T\frac{2\alpha_{t+1}^2\sigma^2\gamma}{T\eta_t} \nonumber \\
 & \leq \frac{f_{\mu}(x_1) - f^*}{T} + \frac{\gamma\sigma^2}{T\eta_0}  + \sum_{t=1}^T\frac{3L^2d^2\mu^2\gamma}{T\eta_t} + \sum_{t=1}^T\frac{2\alpha_{t+1}^2\sigma^2\gamma}{T\eta_t} \nonumber \\
 & = \frac{f_{\mu}(x_1) - f^*}{T} + \frac{\gamma m^{1/3}\sigma^2}{k T}  + \sum_{t=1}^T\frac{3L^2d^2\mu^2\gamma}{T\eta_t} + \sum_{t=1}^T\frac{2c^2\eta_t^3\sigma^2\gamma}{T},
\end{align}
where the second inequality is due to $v_1 = \hat{\nabla}f(x_1,\xi)$ and Assumption \ref{ass:1}.
Since $\eta_t$ is decreasing, i.e., $\eta_T^{-1} \geq \eta_t^{-1}$ for any $0<t<T$, we have
\begin{align}
 & \frac{1}{T} \sum_{t=1}^T\mathbb{E} \big[\frac{1}{4} \|\nabla f_{\mu}(x_{t}) -v_t\|^2 + \frac{1}{4\gamma^2}\|\tilde{x}_{t+1}-x_t\|^2\big] \nonumber \\
 & \leq  \frac{f_{\mu}(x_1) - f^*}{T\eta_T\gamma} + \frac{m^{1/3}\sigma^2}{k T\eta_T}  + \sum_{t=1}^T\frac{3L^2d^2\mu^2}{T\eta_t\eta_T} + \sum_{t=1}^T\frac{2c^2\eta_t^3\sigma^2}{T\eta_T} \nonumber \\
 & \leq \frac{f_{\mu}(x_1) - f^*}{T\eta_T\gamma} + \frac{m^{1/3}\sigma^2}{k T\eta_T}  + \frac{3L^2d^2\mu^2}{T\eta_T} \int^T_1\frac{(m+t)^{1/3}}{k}dt + \frac{2c^2\sigma^2}{T\eta_T} \int^T_1k^3(m+t)^{-1}dt \nonumber \\
 & \leq \frac{f_{\mu}(x_1) - f^*}{T\eta_T\gamma} + \frac{m^{1/3}\sigma^2}{k T\eta_T}  + \frac{9L^2d^2\mu^2}{4kT\eta_T} (m+T)^{4/3} + \frac{2k^3c^2\sigma^2}{T\eta_T} \ln(m+T) \nonumber \\
 & = \frac{f_{\mu}(x_1) - f^*}{T\gamma k}(m+T)^{1/3} + \frac{m^{1/3}\sigma^2}{k^2 T}(m+T)^{1/3}  + \frac{9L^2d^2\mu^2}{4k^2T} (m+T)^{5/3} \nonumber \\
 & \quad + \frac{2k^2c^2\sigma^2}{T}\ln(m+T)(m+T)^{1/3} \nonumber \\
 & \leq \frac{f_{\mu}(x_1) - f^*}{T\gamma k}(m+T)^{1/3} + \frac{m^{1/3}\sigma^2}{k^2 T}(m+T)^{1/3}  + \frac{9L^2}{4k^2T} (m+T)^{1/3} \nonumber \\
 & \quad + \frac{2k^2c^2\sigma^2}{T}\ln(m+T)(m+T)^{1/3},
\end{align}
where the second inequality holds by $\sum_{t=1}^T\frac{1}{\eta_t}dt \leq \int^T_1 \frac{1}{\eta_t}dt = \int^T_1\frac{(m+t)^{1/3}}{k}dt$ and $\sum_{t=1}^T \eta_t^3 dt \leq \int^T_1 \eta_t^3 dt = \int^T_1k^3(m+t)^{-1}$, and the last inequality is due to $0<\mu \leq \frac{1}{d(m+T)^{2/3}}$.
Let $M=\frac{f_{\mu}(x_1) - f^*}{k\gamma} + \frac{m^{1/3}\sigma^2}{k^2} + \frac{9L^2}{4k^2} + 2k^2c^2\sigma^2\ln(m+T)$, we have
\begin{align}
  \frac{1}{T} \sum_{t=1}^T\mathbb{E} \big[\frac{1}{4} \|\nabla f_{\mu}(x_{t}) -v_t\|^2
 + \frac{1}{4\gamma^2} \|\tilde{x}_{t+1}-x_t\|^2 \big] \leq \frac{M}{T}(m+T)^{1/3}.
\end{align}

According to Jensen's inequality, we have
\begin{align}
 & \frac{1}{T} \sum_{t=1}^T\mathbb{E} \big[\frac{1}{2} \|\nabla f_{\mu}(x_{t}) -v_t\|
 + \frac{1}{2\gamma} \|\tilde{x}_{t+1}-x_t\| \big] \nonumber \\
 & \leq \big( \frac{2}{T} \sum_{t=1}^T\mathbb{E} \big[\frac{1}{4} \|\nabla f_{\mu}(x_{t}) -v_t\|^2
 + \frac{1}{4\gamma^2} \|\tilde{x}_{t+1}-x_t\|^2\big]\big)^{1/2} \nonumber \\
 & \leq \frac{\sqrt{2M}}{T^{1/2}}(m+T)^{1/6} \leq \frac{\sqrt{2M}m^{1/6}}{T^{1/2}} + \frac{\sqrt{2M}}{T^{1/3}},
\end{align}
where the last inequality is due to $(a+b)^{1/6} \leq a^{1/6} + b^{1/6}$.
Then we have
\begin{align}
 \frac{1}{T} \sum_{t=1}^T\mathbb{E} \big[ \|\nabla f_{\mu}(x_{t}) -v_t\| + \frac{1}{\gamma} \|\tilde{x}_{t+1}-x_t\| \big] \leq
  \frac{2\sqrt{2M}m^{1/6}}{T^{1/2}} + \frac{2\sqrt{2M}}{T^{1/3}}.
\end{align}

By Lemma \ref{lem:A3}, we have $\|\nabla f_{\mu}(x_{t}) - v_t\|=\|\nabla f_{\mu}(x_{t}) - \nabla f(x_t) + \nabla f(x_t)- v_t\| \geq \|\nabla f(x_t)- v_t\| - \|\nabla f_{\mu}(x_{t}) - \nabla f(x_t)\| \geq \|\nabla f(x_t)- v_t\| - \frac{\mu L d}{2}$.
Thus, we have
\begin{align}
 & \frac{1}{T} \sum_{t=1}^T\mathbb{E} \big[ \|\nabla f(x_{t}) - v_t\| + \frac{1}{\gamma} \|\tilde{x}_{t+1}-x_t\| \big] \nonumber \\
 & \leq \frac{1}{T} \sum_{t=1}^T\mathbb{E} \big[ \|\nabla f_{\mu}(x_{t}) - v_t\| + \frac{\mu L d}{2} + \frac{1}{\gamma} \|\tilde{x}_{t+1}-x_t\| \big] \nonumber \\
 & \leq  \frac{2\sqrt{2M}m^{1/6}}{T^{1/2}} + \frac{2\sqrt{2M}}{T^{1/3}} + \frac{\mu L d}{2} \nonumber \\
 & \leq  \frac{2\sqrt{2M}m^{1/6}}{T^{1/2}} + \frac{2\sqrt{2M}}{T^{1/3}} + \frac{L}{2(m+T)^{2/3}},
\end{align}
where the last inequality is due to $0<\mu \leq \frac{1}{d(m+T)^{2/3}}$. 
Then by using the above inequality \eqref{eq:MG}, we have 
 \begin{align}
  \frac{1}{T}\sum_{t=1}^T\mathbb{E} \|G_{\mathcal{X}}(x_t,\nabla f(x_t),\gamma)\| & \leq \frac{1}{T} \sum_{t=1}^T\mathbb{E} \big[ \|\nabla f(x_{t}) - v_t\| + \frac{1}{\gamma} \|\tilde{x}_{t+1}-x_t\| \big] \nonumber \\
  & \leq  \frac{2\sqrt{2M}m^{1/6}}{T^{1/2}} + \frac{2\sqrt{2M}}{T^{1/3}} + \frac{L}{2(m+T)^{2/3}}.
\end{align}

\end{proof}

\subsection{ Convergence Analysis of Acc-ZOM Algorithm for Unconstrained Mini-Optimization }
\label{Appendix:A2}
In this subsection, we study the convergence properties of our Acc-ZOM algorithm for solving the black-box \textbf{unconstrained}  problem \eqref{eq:1},i.e., $\mathcal{X}= \mathbb{R}^d$.
The following convergence analysis builds on the common metric $\mathbb{E}\|\nabla f(x)\|$ used in \citep{ji2019improved}.

\begin{lemma} \label{lem:B01}
Suppose the sequence $\{x_t\}_{t=1}^T$ be generated from Algorithm \ref{alg:1}.
When $\mathcal{X}=\mathbb{R}^{d}$, given $0<\gamma\leq \frac{1}{2\eta_t L}$,
we have
\begin{align}
  f_{\mu}(x_{t+1}) & \leq f_{\mu}(x_t)
  + \frac{\gamma\eta_t}{2}\|\nabla f_{\mu}(x_t)-v_t\|^2 - \frac{\gamma\eta_t}{2}\|\nabla f_{\mu}(x_t)\|^2
  - \frac{\gamma\eta_t}{4}\|v_t\|^2.
\end{align}
\end{lemma}
\begin{proof}
 According to Assumption \ref{ass:3} and Lemma \ref{lem:A3}, the approximated function $f_{\mu}(x)$ is $L$-smooth.
 Then we have
 \begin{align}
  f_{\mu}(x_{t+1}) &\leq f_{\mu}(x_t) - \gamma\eta_t\langle\nabla f_{\mu}(x_t),v_t\rangle + \frac{\gamma^2\eta_t^2L}{2}\|v_t\|^2   \\
  & = f_{\mu}(x_t) + \frac{\gamma\eta_t}{2}\|\nabla f_{\mu}(x_t)-v_t\|^2 - \frac{\gamma\eta_t}{2}\|\nabla f_{\mu}(x_t)\|^2
  + (\frac{\gamma^2\eta_t^2 L}{2}-\frac{\gamma\eta_t}{2})\|v_t\|^2 \nonumber \\
  & \leq  f_{\mu}(x_t) + \frac{\gamma\eta_t}{2}\|\nabla f_{\mu}(x_t)-v_t\|^2 - \frac{\gamma\eta_t}{2}\|\nabla f_{\mu}(x_t)\|^2 - \frac{\gamma\eta_t}{4}\|v_t\|^2, \nonumber
 \end{align}
 where the last inequality is due to $0< \gamma \leq \frac{1}{2\eta_t L}$.
 Then we have
 \begin{align}
  f_{\mu}(x_{t+1}) & \leq f_{\mu}(x_t)
  + \frac{\gamma\eta_t}{2}\|\nabla f_{\mu}(x_t)-v_t\|^2 - \frac{\gamma\eta_t}{2}\|\nabla f_{\mu}(x_t)\|^2
  - \frac{\gamma\eta_t}{4}\|v_t\|^2.
 \end{align}

\end{proof}

\begin{lemma} \label{lem:C01}
 Suppose the zeroth-order stochastic gradient $\{v_t\}$ be generated from Algorithm \ref{alg:1}, we have
 \begin{align}
 \mathbb{E} \|\nabla f_{\mu}(x_{t+1}) - v_{t+1}\|^2 & \leq (1-\alpha_{t+1})^2 \mathbb{E} \|\nabla f_{\mu}(x_t) -v_t\|^2 + 6(1-\alpha_{t+1})^2 dL^2 \eta_t^2\gamma^2\mathbb{E} \| v_t\|^2 \nonumber \\
 & \quad + 3(1-\alpha_{t+1})^2L^2d^2\mu^2 + 2\alpha_{t+1}^2\sigma^2.
 \end{align}
\end{lemma}
\begin{proof}
 The proof is similar to the proof of Lemma \ref{lem:C1}.
 \end{proof}

 \begin{theorem} \label{th:A01}
  (Restatement of Theorem 3)
 Suppose the sequence $\{x_t\}_{t=1}^T$ be generated from Algorithm \ref{alg:1}.
 When $\mathcal{X}=\mathbb{R}^d$, and let $\eta_t = \frac{k}{(m+t)^{1/3}}$
 for all $t\geq 0$, $0< \gamma \leq \min\big(\frac{m^{1/3}}{2Lk},\frac{1}{2\sqrt{6d}L}\big)$,
 $c\geq \frac{2}{3k^3} + \frac{5}{4}$, $k>0$, $m\geq \max\big(2, k^3, (ck)^3\big)$ and $0<\mu \leq \frac{1}{d(m+T)^{2/3}}$,
 we have
 \begin{align}
  \frac{1}{T} \sum_{t=1}^T\mathbb{E}\|\nabla f(x_t)\|
  \leq  \frac{\sqrt{2M}m^{1/6}}{T^{1/2}} + \frac{\sqrt{2M}}{T^{1/3}} + \frac{L}{2(m+T)^{2/3}},
\end{align}
 where $M=\frac{f_{\mu}(x_1) - f^*}{k\gamma} + \frac{m^{1/3}\sigma^2}{k^2} + \frac{9L^2}{4k^2} + 2k^2c^2\sigma^2\ln(m+T)$.
 \end{theorem}
 \begin{proof}
 This proof is similar to the proof of Theorem \ref{th:A1}. Under the same conditions in Theorem \ref{th:A1},
 by using Lemma \ref{lem:C01} and let $c \geq \frac{2}{3k^3} + \frac{5}{4}$, we have
 \begin{align}
  & \frac{1}{\eta_t}\mathbb{E} \|\nabla f_{\mu}(x_{t+1}) - v_{t+1}\|^2 - \frac{1}{\eta_{t-1}}\mathbb{E} \|\nabla f_{\mu}(x_{t}) - v_t\|^2 \nonumber \\
  & \leq -\frac{5\eta_t}{4}\mathbb{E} \|\nabla f_{\mu}(x_{t}) -v_t\|^2 + 6 dL^2\gamma^2\eta_t \mathbb{E} \| v_{t}\|^2 + \frac{ 3L^2d^2\mu^2}{\eta_t}
  + \frac{2\alpha_{t+1}^2\sigma^2}{\eta_t}.
 \end{align}

 At the same time, we give the  \emph{Lyapunov} function $R_t = \mathbb{E} \big[f_{\mu}(x_t) + \frac{\gamma}{\eta_{t-1}}\|\nabla f_{\mu}(x_t)-v_t\|^2\big]$
 defined in the above Theorem \ref{th:A1}.
 According to Lemma \ref{lem:B01}, we have
 \begin{align}
  R_{t+1} - R_t & = \mathbb{E}\big[f_{\mu}(x_{t+1}) - f_{\mu}(x_t)\big] + \frac{\gamma}{\eta_{t}}\mathbb{E}\|\nabla f_{\mu}(x_{t+1})-v_{t+1}\|^2 - \frac{\gamma}{\eta_{t-1}}\mathbb{E} \|\nabla f_{\mu}(x_t)-v_t\|^2 \nonumber \\
  & \leq \frac{\gamma\eta_t}{2}\mathbb{E}\|\nabla f_{\mu}(x_t)-v_t\|^2 - \frac{\gamma\eta_t}{2}\mathbb{E}\|\nabla f_{\mu}(x_t)\|^2
  - \frac{\gamma\eta_t}{4}\mathbb{E}\|v_t\|^2\nonumber \\
  & \quad -\frac{5\eta_t\gamma}{4}\mathbb{E} \|\nabla f_{\mu}(x_{t}) -v_t\|^2 + 6 dL^2\gamma^3\eta_t \mathbb{E} \| v_{t}\|^2 + \frac{ 3L^2d^2\mu^2\gamma}{\eta_t} + \frac{2\alpha_{t+1}^2\sigma^2\gamma}{\eta_t} \nonumber \\
  & \leq - \frac{\gamma\eta_t}{2}\mathbb{E}\|\nabla f_{\mu}(x_t)\|^2 + \frac{ 3L^2d^2\mu^2\gamma}{\eta_t} + \frac{2\alpha_{t+1}^2\sigma^2\gamma}{\eta_t}  -\big(\frac{\gamma}{4} - 6 dL^2\gamma^3\big)\eta_t\mathbb{E} \| v_{t}\|^2 \nonumber \\
  & \leq - \frac{\gamma\eta_t}{2}\mathbb{E}\|\nabla f_{\mu}(x_t)\|^2 + \frac{ 3L^2d^2\mu^2\gamma}{\eta_t} + \frac{2\alpha_{t+1}^2\sigma^2\gamma}{\eta_t},
 \end{align}
where the last inequality is due to $\gamma \leq \frac{1}{2\sqrt{6d}L}$.
Thus, we can obtain
\begin{align} \label{eq:C03}
\frac{\gamma\eta_t}{2}\mathbb{E}\|\nabla f_{\mu}(x_t)\|^2 \leq R_t - R_{t+1} + \frac{ 3L^2d^2\mu^2\gamma}{\eta_t} + \frac{2\alpha_{t+1}^2\sigma^2\gamma}{\eta_t}.
\end{align}
Since $\inf_{x\in \mathcal{X}} f(x) =f^*$, we have
$\inf_{x\in \mathcal{X}} f_{\mu}(x) =\inf_{x\in \mathcal{X}} \mathbb{E}_{u\sim U_B}[f(x+\mu u)] =
\inf_{x\in \mathcal{X}} \frac{1}{V}\int_{B}f(x+\mu u)du \geq \frac{1}{V}\int_{B}\inf_{x\in \mathcal{X}}f(x+\mu u)du=f^*$,
where $V$ denotes the volume of the unit ball $B$.

Taking average over $t=1,2,\cdots,T$ on both sides of \eqref{eq:C03}, we have
\begin{align}
  &\frac{1}{T} \sum_{t=1}^T\frac{\gamma\eta_t}{2}\mathbb{E} \|\nabla f_{\mu}(x_t)\|^2 \nonumber \\
 & \leq \frac{f_{\mu}(x_1) - f^*}{T} + \frac{\gamma\|\nabla f_{\mu}(x_1) - v_1\|^2}{T\eta_0}  + \sum_{t=1}^T\frac{3L^2d^2\mu^2\gamma}{T\eta_t} + \sum_{t=1}^T\frac{2\alpha_{t+1}^2\sigma^2\gamma}{T\eta_t} \nonumber \\
 & \leq \frac{f_{\mu}(x_1) - f^*}{T} + \frac{\gamma\sigma^2}{T\eta_0}  + \sum_{t=1}^T\frac{3L^2d^2\mu^2\gamma}{T\eta_t} + \sum_{t=1}^T\frac{2\alpha_{t+1}^2\sigma^2\gamma}{T\eta_t} \nonumber \\
 & = \frac{f_{\mu}(x_1) - f^*}{T} + \frac{\gamma m^{1/3}\sigma^2}{k T}  + \sum_{t=1}^T\frac{3L^2d^2\mu^2\gamma}{T\eta_t} + \sum_{t=1}^T\frac{2c^2\eta_t^3\sigma^2\gamma}{T},
\end{align}
where the second inequality is due to $v_1 = \hat{\nabla}f(x_1,\xi)$ and Assumption \ref{ass:1}.
Since $\eta_t$ is decreasing, i.e., $\eta_T^{-1} \geq \eta_t^{-1}$ for any $0<t<T$, we have
\begin{align}
 &\frac{1}{T} \sum_{t=1}^T \frac{1}{2}\mathbb{E} \|\nabla f_{\mu}(x_t)\|^2 \nonumber\\
 & \leq  \frac{f_{\mu}(x_1) - f^*}{T\eta_T\gamma} + \frac{m^{1/3}\sigma^2}{k T\eta_T}  + \sum_{t=1}^T\frac{3L^2d^2\mu^2}{T\eta_t\eta_T} + \sum_{t=1}^T\frac{2c^2\eta_t^3\sigma^2}{T\eta_T} \nonumber \\
 & \leq \frac{f_{\mu}(x_1) - f^*}{T\eta_T\gamma} + \frac{m^{1/3}\sigma^2}{k T\eta_T}  + \frac{3L^2d^2\mu^2}{T\eta_T} \int^T_1\frac{(m+t)^{1/3}}{k}dt \nonumber \\
 & \quad  + \frac{2c^2\sigma^2}{T\eta_T} \int^T_1k^3(m+t)^{-1}dt \nonumber \\
 & \leq \frac{f_{\mu}(x_1) - f^*}{T\eta_T\gamma} + \frac{m^{1/3}\sigma^2}{k T\eta_T}  + \frac{9L^2d^2\mu^2}{4kT\eta_T} (m+T)^{4/3} + \frac{2k^3c^2\sigma^2}{T\eta_T} \ln(m+T) \nonumber \\
 & = \frac{f_{\mu}(x_1) - f^*}{T\gamma k}(m+T)^{1/3} + \frac{m^{1/3}\sigma^2}{k^2 T}(m+T)^{1/3}  + \frac{9L^2d^2\mu^2}{4k^2T} (m+T)^{5/3} \nonumber \\
 & \quad + \frac{2k^2c^2\sigma^2}{T}\ln(m+T)(m+T)^{1/3} \nonumber \\
 & \leq \frac{f_{\mu}(x_1) - f^*}{T\gamma k}(m+T)^{1/3} + \frac{m^{1/3}\sigma^2}{k^2 T}(m+T)^{1/3}  + \frac{9L^2}{4k^2T} (m+T)^{1/3} \nonumber \\
 & \quad + \frac{2k^2c^2\sigma^2}{T}\ln(m+T)(m+T)^{1/3},
\end{align}
where the second inequality holds by $\sum_{t=1}^T\frac{1}{\eta_t}dt \leq \int^T_1 \frac{1}{\eta_t}dt = \int^T_1\frac{(m+t)^{1/3}}{k}dt$ and $\sum_{t=1}^T \eta_t^3 dt \leq \int^T_1 \eta_t^3 dt = \int^T_1k^3(m+t)^{-1}$, and the last inequality is due to $0<\mu \leq \frac{1}{d(m+T)^{2/3}}$.
Let $M=\frac{f_{\mu}(x_1) - f^*}{k\gamma} + \frac{m^{1/3}\sigma^2}{k^2} + \frac{9L^2}{4k^2} + 2k^2c^2\sigma^2\ln(m+T)$, we have
\begin{align}
  \frac{1}{T} \sum_{t=1}^T\mathbb{E} \|\nabla f_{\mu}(x_t)\|^2 \leq \frac{2M}{T}(m+T)^{1/3}.
\end{align}
According to Jensen's inequality, we have
\begin{align}
  \frac{1}{T} \sum_{t=1}^T\mathbb{E} \|\nabla f_{\mu}(x_t)\|
 & \leq \big( \frac{1}{T} \sum_{t=1}^T\mathbb{E} \|\nabla f_{\mu}(x_t)\|^2 \big)^{1/2} \nonumber \\
 & \leq \frac{\sqrt{2M}}{T^{1/2}}(m+T)^{1/6} \leq \frac{\sqrt{2M}m^{1/6}}{T^{1/2}} + \frac{\sqrt{2M}}{T^{1/3}},
\end{align}
where the last inequality is due to $(a+b)^{1/6} \leq a^{1/6} + b^{1/6}$.

By Lemma \ref{lem:A3}, we have $\|\nabla f_{\mu}(x_{t}) \|=\|\nabla f_{\mu}(x_{t}) - \nabla f(x_t) + \nabla f(x_t)\| \geq \|\nabla f(x_t)\|
- \|\nabla f_{\mu}(x_{t}) - \nabla f(x_t)\| \geq \|\nabla f(x_t)\| - \frac{\mu L d}{2}$.
Thus, we have
\begin{align}
 \frac{1}{T} \sum_{t=1}^T\mathbb{E}\|\nabla f(x_t)\|
 & \leq \frac{1}{T} \sum_{t=1}^T\big( \mathbb{E} \|\nabla f_{\mu}(x_{t})\| + \frac{\mu L d}{2} \big) \nonumber \\
 & \leq  \frac{\sqrt{2M}m^{1/6}}{T^{1/2}} + \frac{\sqrt{2M}}{T^{1/3}} + \frac{\mu L d}{2} \nonumber \\
 & \leq  \frac{\sqrt{2M}m^{1/6}}{T^{1/2}} + \frac{\sqrt{2M}}{T^{1/3}} + \frac{L}{2(m+T)^{2/3}},
\end{align}
where the last inequality is due to $0<\mu \leq \frac{1}{d(m+T)^{2/3}}$.

\end{proof}

\subsection{Convergence Analysis of the Acc-ZOMDA Algorithm for Constrained Minimax Optimization}
\label{Appendix:A3}
In the subsection, we study the convergence properties of our  Acc-ZOMDA algorithm for solving the black-box \textbf{constrained}
minimax problem \eqref{eq:2}, i.e., $\mathcal{X}\subset \mathbb{R}^{d_1}$
and $\mathcal{Y} \subset \mathbb{R}^{d_2}$ (or $\mathcal{Y} = \mathbb{R}^{d_2}$), where only the noise function values of $f(x,y)$ can be obtained.
The following convergence analysis builds on a new metric $\mathbb{E}[\mathcal{H}_t]$, where
$\mathcal{H}_t$ is defined in \eqref{eq:14}.

We first let $f_{\mu_1}(x,y) = \mathbb{E}_{u_1\sim U_{B_1}}[f(x+\mu_1u_1,y)]$ and $f_{\mu_2}(x,y) = \mathbb{E}_{u_2\sim U_{B_2}}[f(x,y+\mu_2u_2)]$
denote the smoothing version of $f(x,y)$ w.r.t. $x$ with parameter $\mu_1$ and the smoothing version of $f(x,y)$ w.r.t. $y$
with parameter $\mu_2$, respectively. Here $U_{B_1}$ and
$U_{B_2}$ denote the uniform distributions over the $d_1$-dimensional unit Euclidean ball $B_1$ and
$d_2$-dimensional unit Euclidean ball $B_2$, respectively.
At the same time, let $F_{\mu_1}(x) =\mathbb{E}_{u_1\sim U_{B_1}}[F(x+\mu_1u_1)]$ denote the smoothing approximation of
function $F(x)=\max_{y \in \mathcal{Y}} f(x,y)$.

 \begin{lemma} \label{lem:D1}
 Suppose the sequence $\{x_t,y_t\}_{t=1}^T$ be generated from Algorithm \ref{alg:2}.
 Let $0<\eta_t\leq 1$ and $0< \gamma \leq \frac{1}{2L_g\eta_t}$, we have
 \begin{align}
  F_{\mu_1}(x_{t+1}) - F_{\mu_1}(x_t) & \leq -\frac{\eta_t}{2\gamma}\|\tilde{x}_{t+1}-x_t\|^2 + 6\eta_t\gamma L_f^2\|y^*(x_t)-y_t\|^2 + 2\eta_t\gamma\|\nabla_xf_{\mu_1}(x_t,y_t) -v_t\|^2 \nonumber \\
  & \quad +3\eta_t\gamma\mu_1^2d_1^2L_f^2,
 \end{align}
 where $L_g = L_f + L_f^2/\tau$.
 \end{lemma}
 \begin{proof}
 According to the above Lemma \ref{lem:1} and Lemma \ref{lem:A3},
 the function $F_{\mu_1}(x)$ is $L_g$-smooth.
 By the $L_g$-smoothness of $F_{\mu_1}(x)$, we have
 \begin{align} \label{eq:D1}
  F_{\mu_1}(x_{t+1}) &\leq F_{\mu_1}(x_t) + \langle\nabla F_{\mu_1}(x_t), x_{t+1}-x_t\rangle + \frac{L_g}{2}\|x_{t+1}-x_t\|^2   \\
  & = F_{\mu_1}(x_t) + \eta_t\langle \nabla F_{\mu_1}(x_t),\tilde{x}_{t+1}-x_t\rangle + \frac{L_g\eta_t^2}{2}\|\tilde{x}_{t+1}-x_t\|^2 \nonumber \\
  & = F_{\mu_1}(x_t) + \eta_t\langle \nabla F_{\mu_1}(x_t)-v_t,\tilde{x}_{t+1}-x_t\rangle + \eta_t\langle v_t,\tilde{x}_{t+1}-x_t\rangle + \frac{L_g\eta_t^2}{2}\|\tilde{x}_{t+1}-x_t\|^2. \nonumber
 \end{align}
 By the step 8 of Algorithm \ref{alg:2}, we have $\tilde{x}_{t+1}=\mathcal{P}_{\mathcal{X}}(x_t - \gamma v_t)
 =\arg\min_{x\in\mathcal{X}}\frac{1}{2}\|x- x_t + \gamma v_t\|^2$. Since $\mathcal{X}$ is a convex set
 and the function $\frac{1}{2}\|x- x_t + \gamma v_t\|^2$ is convex, according to Lemma \ref{lem:A1},
 we have
 \begin{align} \label{eq:D2}
  \langle \tilde{x}_{t+1}- x_t + \gamma v_t, x-\tilde{x}_{t+1}\rangle \geq 0, \ \forall x\in \mathcal{X}.
 \end{align}
 In Algorithm \ref{alg:2}, let the initialize solution $x_1 \in \mathcal{X}$, and the sequence $\{x_t\}_{t\geq 1}$ generates as follows:
 \begin{align}
  x_{t+1} = x_{t} + \eta_t(\tilde{x}_{t+1} - x_t) = \eta_t\tilde{x}_{t+1} + (1-\eta_t)x_t,
 \end{align}
 where $0<\eta_t \leq 1$.
 Since $\mathcal{X}$ is convex set and $x_t, \tilde{x}_{t+1} \in \mathcal{X}$, we have $x_{t+1} \in \mathcal{X}$ for any $t \geq 1$.
 Set $x=x_t$ in the inequality \eqref{eq:D2}, we have
 \begin{align} \label{eq:D3}
 \langle v_t, \tilde{x}_{t+1}-x_t\rangle \leq -\frac{1}{\gamma}\|\tilde{x}_{t+1}-x_t\|^2.
 \end{align}

 Next, we decompose the term $\langle \nabla F_{\mu_1}(x_t)-v_t,\tilde{x}_{t+1}-x_t\rangle$ as follows:
 \begin{align}
 &\langle \nabla F_{\mu_1}(x_t)-v_t,\tilde{x}_{t+1}-x_t\rangle \nonumber \\
 & = \underbrace{\langle \nabla F_{\mu_1}(x_t) - \nabla_xf_{\mu_1}(x_t,y_t),\tilde{x}_{t+1}-x_t\rangle}_{=T_1}
 + \underbrace{\langle  \nabla_xf_{\mu_1}(x_t,y_t) -v_t,\tilde{x}_{t+1}-x_t\rangle}_{=T_2}.
 \end{align}
 For the term $T_1$, by the Cauchy-Schwarz inequality and Young's inequality, we have
 \begin{align}
  T_1 &= \langle \nabla F_{\mu_1}(x_t) - \nabla_xf_{\mu_1}(x_t,y_t),\tilde{x}_{t+1}-x_t\rangle \nonumber \\
  & \leq \|\nabla F_{\mu_1}(x_t) - \nabla_xf_{\mu_1}(x_t,y_t)\|\cdot\|\tilde{x}_{t+1}-x_t\| \nonumber \\
  & \leq 2\gamma\|\nabla F_{\mu_1}(x_t) - \nabla_xf_{\mu_1}(x_t,y_t)\|^2 + \frac{1}{8\gamma}\|\tilde{x}_{t+1}-x_t\|^2 \nonumber \\
  & =2\gamma\|\nabla_xf_{\mu_1}(x_t,y^*(x_t)) - \nabla_xf_{\mu_1}(x_t,y_t)\|^2 + \frac{1}{8\gamma}\|\tilde{x}_{t+1}-x_t\|^2 \nonumber \\
  & =2\gamma\|\nabla_xf_{\mu_1}(x_t,y^*(x_t)) -\nabla_xf(x_t,y^*(x_t))+\nabla_xf(x_t,y^*(x_t))-\nabla_xf(x_t,y_t)\nonumber \\
  & \quad +\nabla_xf(x_t,y_t)- \nabla_xf_{\mu_1}(x_t,y_t)\|^2 + \frac{1}{8\gamma}\|\tilde{x}_{t+1}-x_t\|^2 \nonumber \\
  & \leq 6\gamma\|\nabla_xf_{\mu_1}(x_t,y^*(x_t)) -\nabla_xf(x_t,y^*(x_t))\|^2 + 6\gamma\|\nabla_xf(x_t,y^*(x_t))-\nabla_xf(x_t,y_t) \|^2\nonumber \\
  & \quad + 6\gamma\|\nabla_xf(x_t,y_t)- \nabla_xf_{\mu_1}(x_t,y_t)\|^2 + \frac{1}{8\gamma}\|\tilde{x}_{t+1}-x_t\|^2 \nonumber \\
  & \leq 3\gamma\mu_1^2d_1^2L_f^2 + 6\gamma L_f^2\|y^*(x_t)-y_t\|^2 + \frac{1}{8\gamma}\|\tilde{x}_{t+1}-x_t\|^2,
 \end{align}
 where the last inequality holds by Assumption \ref{ass:5}, i.e., implies that the partial gradient $\nabla_x f(x,y)$ is $L_f$-Lipschitz continuous
 and Lemma \ref{lem:A3}, we have
 \begin{align}
  \|\nabla_xf_{\mu_1}(x_t,y^*(x_t)) -\nabla_xf(x_t,y^*(x_t))\| \leq \frac{L_fd_1\mu_1}{2}, \ \|\nabla_xf(x_t,y_t)- \nabla_xf_{\mu_1}(x_t,y_t)\|\leq \frac{L_fd_1\mu_1}{2}, \nonumber
 \end{align}
 and by Assumption \ref{ass:5}, we have
 \begin{align}
 \|\nabla_xf(x_t,y^*(x_t))-\nabla_xf(x_t,y_t) \| \leq \|\nabla f(x_t,y^*(x_t))-\nabla f(x_t,y_t) \|\leq L_f\|y_t-y^*(x_t)\|. \nonumber
 \end{align}
 For the term $T_2$, by the Cauchy-Schwarz inequality and Young's inequality, we have
 \begin{align}
  T_2 & = \langle  \nabla_xf_{\mu_1}(x_t,y_t) -v_t,\tilde{x}_{t+1}-x_t \rangle \nonumber \\
  & \leq \|\nabla_xf_{\mu_1}(x_t,y_t) -v_t\| \cdot \|\tilde{x}_{t+1}-x_t\| \nonumber \\
  & \leq 2\gamma\|\nabla_xf_{\mu_1}(x_t,y_t) -v_t\|^2 + \frac{1}{8\gamma}\|\tilde{x}_{t+1}-x_t\|^2,
 \end{align}
 where the last inequality holds by $\langle a,b\rangle \leq \frac{\lambda}{2}\|a\|^2 + \frac{1}{2\lambda}\|b\|^2$
 with $\lambda=4\gamma$.
 Thus, we have
 \begin{align} \label{eq:D4}
 \langle \nabla F_{\mu_1}(x_t)-v_t,\tilde{x}_{t+1}-x_t\rangle
 & = 3\gamma\mu_1^2d_1^2L_f^2 + 6\gamma L_f^2\|y^*(x_t)-y_t\|^2 + 2\gamma\|\nabla_xf_{\mu_1}(x_t,y_t) -v_t\|^2 \nonumber \\
 & \quad + \frac{1}{4\gamma}\|\tilde{x}_{t+1}-x_t\|^2.
 \end{align}
 Finally, combining the inequalities \eqref{eq:D1}, \eqref{eq:D3} with \eqref{eq:D4}, we have
 \begin{align}
 F_{\mu_1}(x_{t+1}) &\leq F_{\mu_1}(x_t) + 3\eta_t\gamma\mu_1^2d_1^2L_f^2 + 6\eta_t\gamma L_f^2\|y^*(x_t)-y_t\|^2 + 2\eta_t\gamma\|\nabla_xf_{\mu_1}(x_t,y_t) -v_t\|^2 \nonumber \\
 & \quad+ \frac{\eta_t}{4\gamma}\|\tilde{x}_{t+1}-x_t\|^2-\frac{\eta_t}{\gamma}\|\tilde{x}_{t+1}-x_t\|^2 + \frac{L_g\eta_t^2}{2}\|\tilde{x}_{t+1}-x_t\|^2\nonumber \\
 & \leq F_{\mu_1}(x_t) + 3\eta_t\gamma\mu_1^2d_1^2L_f^2 + 6\eta_t\gamma L_f^2\|y^*(x_t)-y_t\|^2
 + 2\eta_t\gamma\|\nabla_xf_{\mu_1}(x_t,y_t) -v_t\|^2 \nonumber \\
 & \quad -\frac{\eta_t}{2\gamma}\|\tilde{x}_{t+1}-x_t\|^2,
 \end{align}
 where the last inequality is due to $0< \gamma \leq \frac{1}{2L_g\eta_t}$.
\end{proof}

\begin{lemma} \label{lem:E1}
Suppose the sequence $\{x_t,y_t\}_{t=1}^T$ be generated from Algorithm \ref{alg:2}.
Under the above assumptions, and set $0< \eta_t\leq 1$
and $0<\lambda\leq \frac{1}{6L_f}$, we have
\begin{align}
\|y_{t+1} - y^*(x_{t+1})\|^2 &\leq (1-\frac{\eta_t\tau\lambda}{4})\|y_t -y^*(x_t)\|^2 -\frac{3\eta_t}{4} \|\tilde{y}_{t+1}-y_t\|^2 \nonumber \\
& \quad + \frac{25\eta_t\lambda}{6\tau}  \|\nabla_y f(x_t,y_t)-w_t\|^2 +  \frac{25\kappa_y^2\eta_t}{6\tau\lambda}\|x_t - \tilde{x}_{t+1}\|^2,
\end{align}
where $\kappa_y = L_f/\tau$.
\end{lemma}
\begin{proof}
According to the assumption \ref{ass:6}, i.e., the function $f(x,y)$ is $\tau$-strongly concave w.r.t $y$,
we have
\begin{align} \label{eq:E1}
 f(x_t,y) & \leq f(x_t,y_t) + \langle\nabla_y f(x_t,y_t), y-y_t\rangle - \frac{\tau}{2}\|y-y_t\|^2 \nonumber \\
 & = f(x_t,y_t) + \langle w_t, y-\tilde{y}_{t+1}\rangle + \langle\nabla_y f(x_t,y_t)-w_t, y-\tilde{y}_{t+1}\rangle \nonumber \\
 & \quad +\langle\nabla_y f(x_t,y_t), \tilde{y}_{t+1}-y_t\rangle- \frac{\tau}{2}\|y-y_t\|^2.
\end{align}
According to the assumption \ref{ass:5}, i.e., the function $f(x,y)$ is $L_f$-smooth, we have
\begin{align} \label{eq:E2}
 -\frac{L_f}{2}\|\tilde{y}_{t+1}-y_t\|^2 \leq f(x_t,\tilde{y}_{t+1}) - f(x_t,y_{t}) -\langle\nabla_y f(x_t,y_t), \tilde{y}_{t+1}-y_t\rangle.
\end{align}
Combining the inequalities \eqref{eq:E1} with \eqref{eq:E2}, we have
\begin{align} \label{eq:E3}
 f(x_t,y) & \leq f(x_t,\tilde{y}_{t+1}) + \langle w_t, y-\tilde{y}_{t+1}\rangle + \langle\nabla_y f(x_t,y_t)-w_t, y-\tilde{y}_{t+1}\rangle \nonumber \\
 & \quad - \frac{\tau}{2}\|y-y_t\|^2 + \frac{L_f}{2}\|\tilde{y}_{t+1}-y_t\|^2.
\end{align}

Next, by the step 10 of Algorithm \ref{alg:2}, we have $\tilde{y}_{t+1}=\mathcal{P}_{\mathcal{Y}}(y_t + \lambda w_t)
 =\arg\min_{y\in\mathcal{Y}}\frac{1}{2}\|y- y_t - \lambda w_t\|^2$. Since $\mathcal{Y}\subset\mathbb{R}^{d_2}$ is a convex set
 and the function $\frac{1}{2}\|y - y_t - \lambda w_t\|^2$ is convex, according to Lemma \ref{lem:A1},
 we have
 \begin{align} \label{eq:E04}
  \langle \tilde{y}_{t+1}- y_t - \lambda w_t, y-\tilde{y}_{t+1}\rangle \geq 0, \ \forall y\in \mathcal{Y}.
 \end{align} 
 When $\mathcal{Y}=\mathbb{R}^{d_2}$, clearly, we still can obtain the above inequality \eqref{eq:E04}. 
Then we obtain
\begin{align} \label{eq:E4}
  \langle w_t, y-\tilde{y}_{t+1}\rangle & \leq \frac{1}{\lambda}\langle \tilde{y}_{t+1}- y_t, y-\tilde{y}_{t+1}\rangle \nonumber \\
  & = \frac{1}{\lambda}\langle \tilde{y}_{t+1}- y_t, y_t-\tilde{y}_{t+1}\rangle + \frac{1}{\lambda}\langle \tilde{y}_{t+1}- y_t, y-y_t\rangle \nonumber \\
  & = -\frac{1}{\lambda}\|\tilde{y}_{t+1}- y_t\|^2 + \frac{1}{\lambda}\langle \tilde{y}_{t+1}- y_t, y-y_t\rangle.
 \end{align}
Combining the inequalities \eqref{eq:E3} with \eqref{eq:E4}, we have
\begin{align}
 f(x_t,y) & \leq f(x_t,\tilde{y}_{t+1}) + \frac{1}{\lambda}\langle \tilde{y}_{t+1}- y_t, y-y_t\rangle + \langle\nabla_y f(x_t,y_t)-w_t, y-\tilde{y}_{t+1}\rangle \nonumber \\
 & \quad -\frac{1}{\lambda}\|\tilde{y}_{t+1}- y_t\|^2- \frac{\tau}{2}\|y-y_t\|^2 + \frac{L_f}{2}\|\tilde{y}_{t+1}-y_t\|^2.
\end{align}
Let $y=y^*(x_t)$ and we obtain
\begin{align}
 f(x_t,y^*(x_t)) & \leq f(x_t,\tilde{y}_{t+1}) + \frac{1}{\lambda}\langle \tilde{y}_{t+1}- y_t, y^*(x_t)-y_t\rangle
 + \langle\nabla_y f(x_t,y_t)-w_t, y^*(x_t)-\tilde{y}_{t+1}\rangle \nonumber \\
 & \quad -\frac{1}{\lambda}\|\tilde{y}_{t+1}- y_t\|^2- \frac{\tau}{2}\|y^*(x_t)-y_t\|^2 + \frac{L_f}{2}\|\tilde{y}_{t+1}-y_t\|^2.
\end{align}
Due to the concavity of $f(\cdot,y)$ and $y^*(x_t) =\arg\max_{y\in \mathcal{Y}} f(x_t,y)$, we have $f(x_t,y^*(x_t)) \geq f(x_t,\tilde{y}_{t+1})$.
Thus, we obtain
\begin{align} \label{eq:E5}
 0 & \leq  \frac{1}{\lambda}\langle \tilde{y}_{t+1}- y_t, y^*(x_t)-y_t\rangle
 + \langle\nabla_y f(x_t,y_t)-w_t, y^*(x_t)-\tilde{y}_{t+1}\rangle \nonumber \\
 & \quad -(\frac{1}{\lambda}-\frac{L_f}{2})\|\tilde{y}_{t+1}- y_t\|^2- \frac{\tau}{2}\|y^*(x_t)-y_t\|^2.
\end{align}

By $y_{t+1} = y_t + \eta_t(\tilde{y}_{t+1}-y_t) $, we have
\begin{align}
 \|y_{t+1}-y^*(x_t)\|^2 & = \|y_t + \eta_t(\tilde{y}_{t+1}-y_t) -y^*(x_t)\|^2 \nonumber \\
 & = \|y_t -y^*(x_t)\|^2 + 2\eta_t\langle \tilde{y}_{t+1}-y_t, y_t -y^*(x_t)\rangle + \eta_t^2\|\tilde{y}_{t+1}-y_t\|^2.
\end{align}
Then we obtain
\begin{align} \label{eq:E6}
 \langle \tilde{y}_{t+1}-y_t, y^*(x_t) - y_t\rangle \leq \frac{1}{2\eta_t}\|y_t -y^*(x_t)\|^2 + \frac{\eta_t}{2}\|\tilde{y}_{t+1}-y_t\|^2 - \frac{1}{2\eta_t}\|y_{t+1}-y^*(x_t)\|^2.
\end{align}
Considering the upper bound of the term $\langle\nabla_y f(x_t,y_t)-w_t, y^*(x_t)-\tilde{y}_{t+1}\rangle$, we have
\begin{align} \label{eq:E7}
 &\langle\nabla_y f(x_t,y_t)-w_t, y^*(x_t)-\tilde{y}_{t+1}\rangle \nonumber \\
 & = \langle\nabla_y f(x_t,y_t)-w_t, y^*(x_t)-y_t\rangle + \langle\nabla_y f(x_t,y_t)-w_t, y_t-\tilde{y}_{t+1}\rangle \nonumber \\
 & \leq \frac{1}{\tau} \|\nabla_y f(x_t,y_t)-w_t\|^2 + \frac{\tau}{4}\|y^*(x_t)-y_t\|^2 + \frac{1}{\tau} \|\nabla_y f(x_t,y_t)-w_t\|^2 + \frac{\tau}{4}\|y_t-\tilde{y}_{t+1}\|^2 \nonumber \\
 & = \frac{2}{\tau} \|\nabla_y f(x_t,y_t)-w_t\|^2 + \frac{\tau}{4}\|y^*(x_t)-y_t\|^2 + \frac{\tau}{4}\|y_t-\tilde{y}_{t+1}\|^2.
\end{align}

Next, combining the inequalities \eqref{eq:E5}, \eqref{eq:E6} with \eqref{eq:E7},
we have
\begin{align}
 &\frac{1}{2\eta_t\lambda}\|y_{t+1}-y^*(x_t)\|^2 \nonumber \\
 & \leq ( \frac{1}{2\eta_t\lambda}-\frac{\tau}{4})\|y_t -y^*(x_t)\|^2 + ( \frac{\eta_t}{2\lambda} + \frac{\tau}{4} + \frac{L_f}{2}-\frac{1}{\lambda}) \|\tilde{y}_{t+1}-y_t\|^2 \nonumber \\
 & \quad + \frac{2}{\tau} \|\nabla_y f(x_t,y_t)-w_t\|^2 \nonumber \\
 & \leq ( \frac{1}{2\eta_t\lambda}-\frac{\tau}{4})\|y_t -y^*(x_t)\|^2 + (\frac{3L_f}{4} -\frac{1}{2\lambda}) \|\tilde{y}_{t+1}-y_t\|^2 + \frac{2}{\tau} \|\nabla_y f(x_t,y_t)-w_t\|^2 \nonumber \\
 & = ( \frac{1}{2\eta_t\lambda}-\frac{\tau}{4})\|y_t -y^*(x_t)\|^2 - \big( \frac{3}{8\lambda} + \frac{1}{8\lambda} -\frac{3L_f}{4}\big) \|\tilde{y}_{t+1}-y_t\|^2 \nonumber \\
 &\quad + \frac{2}{\tau} \|\nabla_y f(x_t,y_t)-w_t\|^2 \nonumber \\
 & \leq  ( \frac{1}{2\eta_t\lambda}-\frac{\tau}{4})\|y_t -y^*(x_t)\|^2 - \frac{3}{8\lambda} \|\tilde{y}_{t+1}-y_t\|^2 + \frac{2}{\tau}\|\nabla_y f(x_t,y_t)-w_t\|^2,
\end{align}
where the second inequality holds by $L_f \geq \tau$ and $0< \eta_t\leq 1$, and the last inequality is due to
$0< \lambda \leq \frac{1}{6L_f}$.
It implies that
\begin{align} \label{eq:E8}
\|y_{t+1}-y^*(x_t)\|^2 \leq ( 1-\frac{\eta_t\tau\lambda}{2})\|y_t -y^*(x_t)\|^2 - \frac{3\eta_t}{4} \|\tilde{y}_{t+1}-y_t\|^2
+ \frac{4\eta_t\lambda}{\tau}\|\nabla_y f(x_t,y_t)-w_t\|^2.
\end{align}

Next, we decompose the term $\|y_{t+1} - y^*(x_{t+1})\|^2$ as follows:
\begin{align} \label{eq:E9}
  &\|y_{t+1} - y^*(x_{t+1})\|^2 \nonumber\\ 
  & = \|y_{t+1} - y^*(x_t) + y^*(x_t) - y^*(x_{t+1})\|^2    \nonumber \\
  & =  \|y_{t+1} - y^*(x_t)\|^2 + 2\langle y_{t+1} - y^*(x_t), y^*(x_t) - y^*(x_{t+1})\rangle  + \|y^*(x_t) - y^*(x_{t+1})\|^2  \nonumber \\
  & \leq (1+\frac{\eta_t\tau\lambda}{4})\|y_{t+1} - y^*(x_t)\|^2  + (1+\frac{4}{\eta_t\tau\lambda})\|y^*(x_t) - y^*(x_{t+1})\|^2 \nonumber \\
  & \leq (1+\frac{\eta_t\tau\lambda}{4})\|y_{t+1} - y^*(x_t)\|^2  + (1+\frac{4}{\eta_t\tau\lambda})\kappa_y^2\|x_t - x_{t+1}\|^2\nonumber \\
  & = (1+\frac{\eta_t\tau\lambda}{4})\|y_{t+1} - y^*(x_t)\|^2  + (1+\frac{4}{\eta_t\tau\lambda})\kappa_y^2\eta_t^2\|x_t - \tilde{x}_{t+1}\|^2,
\end{align}
where the first inequality holds by the Cauchy-Schwarz inequality and Young's inequality, and  the second inequality is due to
Lemma \ref{lem:2}, and  the last equality holds by $x_{t+1}=x_t + \eta_t(\tilde{x}_{t+1}-x_t)$.

Combining the above inequalities \eqref{eq:E8} and \eqref{eq:E9}, we have
\begin{align}
 \|y_{t+1} - y^*(x_{t+1})\|^2 & \leq (1+\frac{\eta_t\tau\lambda}{4})( 1-\frac{\eta_t\tau\lambda}{2})\|y_t -y^*(x_t)\|^2
 - (1+\frac{\eta_t\tau\lambda}{4})\frac{3\eta_t}{4} \|\tilde{y}_{t+1}-y_t\|^2     \nonumber \\
 & \quad + (1+\frac{\eta_t\tau\lambda}{4})\frac{4\eta_t\lambda}{\tau}\|\nabla_y f(x_t,y_t)-w_t\|^2
  + (1+\frac{4}{\eta_t\tau\lambda})\kappa_y^2\eta_t^2\|x_t - \tilde{x}_{t+1}\|^2.
\end{align}
Since $0 < \eta_t \leq 1$, $0< \lambda \leq \frac{1}{6L_f}$ and $L_f\geq \tau$, we have $\lambda \leq \frac{1}{6L_f} \leq \frac{1}{6\tau}$
and $\lambda\eta_t\leq \frac{1}{6\tau }$. Then we obtain
\begin{align}
  (1+\frac{\eta_t\tau\lambda}{4})( 1-\frac{\eta_t\tau\lambda}{2})&= 1-\frac{\eta_t\tau\lambda}{2} +\frac{\eta_t\tau\lambda}{4}
 - \frac{\eta_t^2\tau^2\lambda^2}{8} \leq 1-\frac{\eta_t\tau\lambda}{4}, \nonumber \\
 - (1+\frac{\eta_t\tau\lambda}{4})\frac{3\eta_t}{4} &\leq -\frac{3\eta_t}{4}, \nonumber \\
  (1+\frac{\eta_t\tau\lambda}{4})\frac{4\eta_t\lambda}{\tau} & \leq (1+\frac{1}{24})\frac{4\eta_t\lambda}{\tau}=\frac{25\eta_t\lambda}{6\tau}, \nonumber \\
  (1+\frac{4}{\eta_t\tau\lambda})\kappa_y^2\eta_t^2 & =  \kappa_y^2\eta_t^2 +\frac{4\kappa_y^2\eta_t}{\tau\lambda}
  \leq \frac{\kappa_y^2\eta_t}{6\tau\lambda} +\frac{4\kappa_y^2\eta_t}{\tau\lambda} = \frac{25\kappa_y^2\eta_t}{6\tau\lambda}.
\end{align}
Thus, we have
\begin{align}
     \|y_{t+1} - y^*(x_{t+1})\|^2 &\leq (1-\frac{\eta_t\tau\lambda}{4})\|y_t -y^*(x_t)\|^2 -\frac{3\eta_t}{4} \|\tilde{y}_{t+1}-y_t\|^2 \nonumber \\
     & \quad + \frac{25\eta_t\lambda}{6\tau}  \|\nabla_y f(x_t,y_t)-w_t\|^2 +  \frac{25\kappa_y^2\eta_t}{6\tau\lambda}\|x_t - \tilde{x}_{t+1}\|^2.
\end{align}
\end{proof}

\begin{lemma} \label{lem:F1}
 Suppose the zeroth-order stochastic gradients $\{v_t,w_t\}_{t=1}^T$ be generated from Algorithm \ref{alg:2}, we have
\begin{align} \label{eq:F1}
 &\mathbb{E} \|\nabla_x f_{\mu_1}(x_{t+1},y_{t+1}) - v_{t+1}\|^2 \nonumber\\
 & \leq (1-\alpha_{t+1})^2 \mathbb{E} \|\nabla_x f_{\mu_1}(x_t,y_t) -v_t\|^2
 + \frac{3(1-\alpha_{t+1})^2L^2_f\mu_1^2d_1^2}{b} \nonumber \\
 & \quad + \frac{6d_1L^2_f(1-\alpha_{t+1})^2\eta^2_t}{b}\mathbb{E}\big(\|\tilde{x}_{t+1}-x_t\|^2 + \|\tilde{y}_{t+1}-y_t\|^2\big) + \frac{2\alpha_{t+1}^2\delta^2}{b}.
\end{align}
\begin{align} \label{eq:F2}
 &\mathbb{E} \|\nabla_y f_{\mu_2}(x_{t+1},y_{t+1}) - w_{t+1}\|^2 \nonumber \\
 & \leq (1-\beta_{t+1})^2 \mathbb{E} \|\nabla_y f_{\mu_2}(x_t,y_t) -w_t\|^2
 + \frac{3(1-\beta_{t+1})^2L^2_f\mu_2^2d_2^2}{b} \nonumber \\
 & \quad + \frac{6d_2L^2_f(1-\beta_{t+1})^2\eta^2_t}{b}\mathbb{E}\big(\|\tilde{x}_{t+1}-x_t\|^2 + \|\tilde{y}_{t+1}-y_t\|^2\big) + \frac{2\beta_{t+1}^2\delta^2}{b}.
\end{align}
\end{lemma}
\begin{proof}
 We first prove the inequality \eqref{eq:F1}.
 According to the definition of $v_{t+1}$ in Algorithm \ref{alg:2}, we have
 \begin{align}
  v_{t+1}-v_t & = -\alpha_{t+1}v_t + (1-\alpha_{t+1})\big(\hat{\nabla}_x f(x_{t+1},y_{t+1};\mathcal{B}_{t+1}) - \hat{\nabla}_x f(x_t,y_t;\mathcal{B}_{t+1})\big) \nonumber \\
  & \quad + \alpha_{t+1}\hat{\nabla}_x f(x_{t+1},y_{t+1};\mathcal{B}_{t+1}).
 \end{align}
 Then we have
 \begin{align} \label{eq:F3}
  & \mathbb{E} \|\nabla_x f_{\mu_1}(x_{t+1},y_{t+1}) - v_{t+1}\|^2  \nonumber\\
  & \! = \mathbb{E} \|\nabla_x f_{\mu_1}(x_{t+1},y_{t+1}) - v_t - (v_{t+1}-v_t)\|^2 \nonumber \\
  & \!= \mathbb{E} \|\nabla_x f_{\mu_1}(x_{t+1},y_{t+1}) - v_t + \alpha_{t+1}v_t- \alpha_{t+1}\hat{\nabla}_x f(x_{t+1},y_{t+1};\mathcal{B}_{t+1})
   \nonumber \\
  & \ - (1-\alpha_{t+1})\big(\hat{\nabla}_x f(x_{t+1},y_{t+1};\mathcal{B}_{t+1})- \hat{\nabla}_x f(x_t,y_t;\mathcal{B}_{t+1})\big) \|^2 \nonumber \\
  &\! = \mathbb{E} \|(1-\alpha_{t+1})(\nabla_x f_{\mu_1}(x_t,y_t) -v_t) + \alpha_{t+1}\big(\nabla_x f_{\mu_1}(x_{t+1},y_{t+1})- \hat{\nabla}_x f(x_{t+1},y_{t+1};\mathcal{B}_{t+1})\big) \nonumber \\
  & \ + (1\!-\!\alpha_{t+1})\big(\nabla_x f_{\mu_1}(x_{t+1},y_{t+1})
  \!-\!\nabla_x f_{\mu_1}(x_t,y_t)\!-\!\hat{\nabla}_x f(x_{t+1},y_{t+1};\mathcal{B}_{t+1}) \!+\! \hat{\nabla}_x f(x_t,y_t;\mathcal{B}_{t+1})\big) \|^2 \nonumber \\
  &\! = (1-\alpha_{t+1})^2\mathbb{E}\|\nabla_x f_{\mu_1}(x_t,y_t)-v_t\|^2 + \mathbb{E}\|\alpha_{t+1}\big(\nabla_x f_{\mu_1}(x_{t+1},y_{t+1})- \hat{\nabla}_x f(x_{t+1},y_{t+1};\mathcal{B}_{t+1})\big) \nonumber \\
  & \ + (1\!-\!\alpha_{t+1})\big(\nabla_x f_{\mu_1}(x_{t+1},y_{t+1})
  \!-\!\nabla_x f_{\mu_1}(x_t,y_t)\!-\!\hat{\nabla}_x f(x_{t+1},y_{t+1};\mathcal{B}_{t+1}) \!+\! \hat{\nabla}_x f(x_t,y_t;\mathcal{B}_{t+1})\big) \|^2  \nonumber \\
  &\! \leq (1-\alpha_{t+1})^2\mathbb{E} \|\nabla_x f_{\mu_1}(x_t,y_t) -v_t\|^2 + \frac{2\alpha_{t+1}^2}{b} \mathbb{E} \|\nabla_x f_{\mu_1}(x_{t+1},y_{t+1})- \hat{\nabla}_x f(x_{t+1},y_{t+1};\xi^t_1)\|^2  \nonumber \\
  & \ + \frac{2(1-\alpha_{t+1})^2}{b}\mathbb{E} \|\nabla_x f_{\mu_1}(x_{t+1},y_{t+1}) \!-\! \nabla_x f_{\mu_1}(x_t,y_t)\!-\! \hat{\nabla}_x f(x_{t+1},y_{t+1};\xi_1^t) \!+\! \hat{\nabla}_x f(x_t,y_t;\xi_1^t)\|^2 \nonumber \\
  &\! \leq (1-\alpha_{t+1})^2 \mathbb{E} \|\nabla_x f_{\mu_1}(x_t,y_t) -v_t\|^2 + \frac{2\alpha_{t+1}^2\delta^2}{b} \nonumber \\
  & \ + \frac{2(1-\alpha_{t+1})^2}{b}\underbrace{ \mathbb{E} \| \hat{\nabla}_x f(x_{t+1},y_{t+1};\xi^t_1) -
  \hat{\nabla}_x f(x_t,y_t;\xi^t_1)\|^2}_{=T_1}, 
 \end{align}
 where the fourth equality follows by $\mathbb{E}_{(\hat{U},\mathcal{B}_{t+1})}[\hat{\nabla}_x f(x_{t+1},y_{t+1};\mathcal{B}_{t+1})]=\nabla_x f_{\mu_1}(x_{t+1},y_{t+1})$ and $ \mathbb{E}_{(\hat{U},\mathcal{B}_{t+1})}[\hat{\nabla}_x f(x_{t+1},y_{t+1};\mathcal{B}_{t+1}) - \hat{\nabla}_x f(x_t,y_t;\mathcal{B}_{t+1})]=\nabla_x f_{\mu_1}(x_{t+1},y_{t+1}) - \nabla_x f_{\mu_1}(x_t,y_t)$; the first inequality holds by Young's inequality and the above lemma \ref{lem:A4}; the last inequality is due to the equality $\mathbb{E}\|\zeta-\mathbb{E}[\zeta]\|^2 =\mathbb{E}\|\zeta\|^2 - \|\mathbb{E}[\zeta]\|^2$ and Assumption \ref{ass:2}.

 Next, we consider the upper bound of the above term $T_1$ as follows:
 \begin{align}
 & T_1 = \mathbb{E} \big\| \hat{\nabla}_x f(x_{t+1},y_{t+1};\xi^t_1) - \hat{\nabla}_x f(x_t,y_t;\xi^t_1)\big\|^2  \\
 & = \mathbb{E} \big\| \frac{d_1(f(x_{t+1}+\mu_1u_1,y_{t+1};\xi^t_1)-f(x_{t+1},y_{t+1};\xi^t_1))}{\mu_1}u_1 \nonumber \\
 &\quad - \frac{d_1(f(x_t+\mu_1u_1,y_t;\xi^t_1)-f(x_t,y_t;\xi^t_1))}{\mu_1}u_1 \big\|^2 \nonumber \\
 & =  d_1^2\mathbb{E} \big\| \frac{f(x_{t+1}+\mu_1u_1,y_{t+1};\xi^t_1)-f(x_{t+1},y_{t+1};\xi^t_1)-\big\langle \nabla_x  f(x_{t+1},y_{t+1};\xi^t_1),\mu_1u_1\big\rangle}{\mu_1}u_1 \nonumber \\
 &\qquad \qquad +\big(\big\langle \nabla_x f(x_{t+1},y_{t+1};\xi^t_1),u_1\big\rangle-\big\langle \nabla_x f(x_t,y_t;\xi^t_1),u_1\big\rangle\big)u_1
 \nonumber \\
 &\qquad \qquad - \frac{f(x_t+\mu_1u_1,y_t;\xi^t_1)-f(x_t,y_t;\xi^t_1)-\big\langle \nabla_x  f(x_t,y_t;\xi^t_1),\mu_1u_1\big\rangle}{\mu_1}u_1\big\|^2 \nonumber \\
 & \leq \frac{3L^2_f\mu_1^2d_1^2}{2} + 3d_1^2\mathbb{E}\big\| \big\langle \nabla_x  f(x_{t+1},y_{t+1};\xi^t_1)-\nabla_x  f(x_t,y_t;\xi^t_1),u_1\big\rangle u_1\big\|^2 \nonumber \\
 & = \frac{3L^2_f\mu_1^2d_1^2}{2} + 3d_1^2\mathbb{E} \big\langle \nabla_x  f(x_{t+1},y_{t+1};\xi^t_1)-\nabla_x  f(x_t,y_t;\xi^t_1),u_1\big\rangle^2
 \nonumber \\
 & = \frac{3L^2_f\mu_1^2d_1^2}{2} + 3d_1^2\mathbb{E}\big[ \big(\nabla_x  f(x_{t+1},y_{t+1};\xi^t_1)-\nabla_x  f(x_t,y_t;\xi^t_1)\big)^T(u_1u_1^T)\nonumber \\
 &\quad \cdot \big(\nabla_x  f(x_{t+1},y_{t+1};\xi^t_1)-\nabla_x  f(x_t,y_t;\xi^t_1)\big)\big], \nonumber
 \end{align}
where the above inequality is due to Young's inequality and Assumption \ref{ass:5}, i.e., $f(x,y;\xi)$ is $L_f$-smooth w.r.t $x$, so we have
$f(x_{t+1}+\mu_1u_1,y_{t+1};\xi^t_1)-f(x_{t+1},y_{t+1};\xi^t_1)-\big\langle \nabla_x  f(x_{t+1},y_{t+1};\xi^t_1)
\\ ,\mu_1u_1\big\rangle  \leq \frac{L_f}{2}\|\mu_1u_1\|^2$
and $f(x_t+\mu_1u_1,y_t;\xi^t_1)-f(x_t,y_t;\xi^t_1)-\big\langle \nabla_x  f(x_t,y_t;\xi^t_1),\mu_1u_1\big\rangle \leq \frac{L_f}{2}\|\mu_1u_1\|^2$, and the forth equality holds by $\|u_1\|=1$.

Following the proof of Lemma 5 in \citep{ji2019improved}, we have $u_1^Tu_1 = \frac{1}{d_1}I_{d_1}$.
Thus, we have
\begin{align} \label{eq:F4}
T_1 & \leq \frac{3L^2_f\mu_1^2d_1^2}{2} + 3d_1\mathbb{E} \|\nabla_x  f(x_{t+1},y_{t+1};\xi^t_1)-\nabla_x  f(x_t,y_t;\xi^t_1)\|^2 \nonumber \\
& \leq \frac{3L^2_f\mu_1^2d_1^2}{2} + 3d_1L^2_f\mathbb{E}\big(\|x_{t+1}-x_t\|^2 + \|y_{t+1}-y_t\|^2\big) \nonumber \\
& = \frac{3L^2_f\mu_1^2d_1^2}{2} + 3d_1L^2_f\eta_t^2\mathbb{E}\big(\|\tilde{x}_{t+1}-x_t\|^2 + \|\tilde{y}_{t+1}-y_t\|^2\big),
\end{align}
where the last inequality holds by Assumption \ref{ass:5}.
Plugging the above inequality \eqref{eq:F4} into \eqref{eq:F3}, we obtain
\begin{align}
 &\mathbb{E} \|\nabla_x f_{\mu_1}(x_{t+1},y_{t+1}) - v_{t+1}\|^2\nonumber \\
 & \leq (1-\alpha_{t+1})^2 \mathbb{E} \|\nabla_x f_{\mu_1}(x_t,y_t) -v_t\|^2
 + \frac{3(1-\alpha_{t+1})^2L^2_f\mu_1^2d_1^2}{b} \nonumber \\
 & \quad + \frac{6d_1L^2_f(1-\alpha_{t+1})^2\eta^2_t}{b}\mathbb{E}\big(\|\tilde{x}_{t+1}-x_t\|^2 + \|\tilde{y}_{t+1}-y_t\|^2\big)+ \frac{2\alpha_{t+1}^2\delta^2}{b}. \nonumber
\end{align}

We apply a similar analysis to prove the above inequality \eqref{eq:F2}. We obtain
\begin{align}
 & \mathbb{E} \|\nabla_y f_{\mu_2}(x_{t+1},y_{t+1}) - w_{t+1}\|^2 \nonumber \\
 & \leq (1-\beta_{t+1})^2 \mathbb{E} \|\nabla_y f_{\mu_2}(x_t,y_t) -w_t\|^2
 + \frac{3(1-\beta_{t+1})^2L^2_f\mu_2^2d_2^2}{b} \nonumber \\
 &\quad + \frac{6d_2L^2_f(1-\beta_{t+1})^2\eta^2_t}{b}\mathbb{E}\big(\|\tilde{x}_{t+1}-x_t\|^2 + \|\tilde{y}_{t+1}-y_t\|^2\big)
 + \frac{2\beta_{t+1}^2\delta^2}{b}. \nonumber
\end{align}
 \end{proof}

\begin{theorem} \label{th:A2}
 (Restatement of Theorem 5)
Suppose the sequence $\{x_t,y_t\}_{t=1}^T$ be generated from Algorithm \ref{alg:2}. When $\mathcal{X}\subset \mathbb{R}^{d_1}$, and let $\eta_t = \frac{k}{(m+t)^{1/3}}$
for all $t\geq 0$, $c_1 \geq \frac{2}{3k^3} + \frac{9\tau^2}{4}$ and $c_2 \geq \frac{2}{3k^3} + \frac{625\tilde{d}L^2_f}{3b}$,
$k>0$, $1\leq b\leq \tilde{d}$, $m\geq \max\big( 2, k^3, (c_1k)^3, (c_2k)^3\big)$,  $0<\lambda\leq \min\big(\frac{1}{6L_f},\frac{75\tau}{24}\big)$,
$0< \gamma \leq \min\big( \frac{\lambda\tau}{2L_f}\sqrt{\frac{6b/\tilde{d}}{36 \lambda^2 + 625\kappa_y^2}}, \frac{m^{1/3}}{2L_gk}\big)$, $0<\mu_1\leq \frac{1}{d_1(m+T)^{2/3}}$
and $0<\mu_2\leq \frac{1}{\tilde{d}^{1/2}d_2(m+T)^{2/3}}$,
we have
\begin{align}
 \frac{1}{T} \sum_{t=1}^T \mathbb{E}  \|G_{\mathcal{X}}(x_t,\nabla F(x_t),\gamma)\| & \!\leq\!\frac{1}{T} \sum_{t=1}^T \mathbb{E}\big[ L_f\|y^*(x_t)-y_t\| \!+\! \|\nabla_xf(x_t,y_t) -v_t\| \!+\!\frac{1}{\gamma}\|\tilde{x}_{t+1}-x_t\|\big]
 \nonumber \\
 & \!\leq \! \frac{2\sqrt{3M'}m^{1/6}}{T^{1/2}} + \frac{2\sqrt{3M'}}{T^{1/3}} + \frac{L_f}{2(m+T)^{2/3}},
\end{align}
where $\Delta_1=\|y_1-y^*(x_1)\|^2$ and  $M' =  \frac{F_{\mu_1}(x_1) - F^*}{\gamma k} + \frac{ 25\tilde{d}L^2_f}{k\lambda\tau b}\Delta_1 +  \frac{2m^{1/3}\delta^2}{b\tau^2 k^2} + \frac{36\tau^2L_f^2 + 625L^4_f}{8b\tau^2}(m+T)^{-2/3} + \frac{9L^2_f}{4b\tau^2 k^2} + \frac{2(c_1^2+c_2^2)\delta^2 k^2}{b\tau^2}\ln(m+T)$.
\end{theorem}
\begin{proof}
 Since $\eta_t=\frac{k}{(m+t)^{1/3}}$ on $t$ is decreasing and $m\geq k^3$, we have $\eta_t \leq \eta_0 = \frac{k}{m^{1/3}} \leq 1$ and $\gamma \leq \frac{m^{1/3}}{2L_gk}=\frac{1}{2L_g\eta_0} \leq \frac{1}{2L_g\eta_t}$ for any $t\geq 0$.
 Due to $0 < \eta_t \leq 1$ and $m\geq \max\big( (c_1k)^3, (c_2k)^3 \big)$, we have $\alpha_{t+1} = c_1\eta_t^2 \leq c_1\eta_t \leq \frac{c_1k}{m^{1/3}}\leq 1$ and $\beta_{t+1} = c_2\eta_t^2 \leq c_2\eta_t \leq \frac{c_2k}{m^{1/3}}\leq 1$.
 According to Lemma \ref{lem:F1}, we have
 \begin{align}
  & \frac{1}{\eta_t}\mathbb{E} \|\nabla_x f_{\mu_1}(x_{t+1},y_{t+1}) - v_{t+1}\|^2 - \frac{1}{\eta_{t-1}}\mathbb{E} \|\nabla_x f_{\mu_1}(x_t,y_t) - v_t\|^2  \\
  & \leq \big(\frac{(1-\alpha_{t+1})^2}{\eta_t} - \frac{1}{\eta_{t-1}}\big)\mathbb{E} \|\nabla_x f_{\mu_1}(x_t,y_t) -v_t\|^2 + \frac{3(1-\alpha_{t+1})^2L^2_f\mu_1^2d_1^2}{b\eta_t} + \frac{2\alpha_{t+1}^2\delta^2}{b\eta_t}\nonumber \\
  & \quad + \frac{6d_1L^2_f(1-\alpha_{t+1})^2\eta_{t}}{b}\mathbb{E}\big(\|\tilde{x}_{t+1}-x_t\|^2 + \|\tilde{y}_{t+1}-y_t\|^2\big) \nonumber \\
  & \leq \big(\frac{1\!-\!\alpha_{t+1}}{\eta_t} \!-\! \frac{1}{\eta_{t-1}}\big)\mathbb{E} \|\nabla_x f_{\mu_1}(x_t,y_t) -v_t\|^2 \!+\! \frac{6d_1L^2_f\eta_{t}}{b}\mathbb{E}\big(\|\tilde{x}_{t+1}-x_t\|^2 + \|\tilde{y}_{t+1}-y_t\|^2\big) \nonumber \\
  & \quad + \frac{3L^2_f\mu_1^2d_1^2}{b\eta_t} \!+\! \frac{2\alpha_{t+1}^2\delta^2}{b\eta_t}\nonumber \\
  & = \big(\frac{1}{\eta_t} \!-\! \frac{1}{\eta_{t-1}} - c_1\eta_t\big)\mathbb{E} \|\nabla_x f_{\mu_1}(x_t,y_t) -v_t\|^2 \!+\! \frac{6d_1L^2_f\eta_{t}}{b}\mathbb{E}\big(\|\tilde{x}_{t+1}\!-\!x_t\|^2 + \|\tilde{y}_{t+1}\!-\!y_t\|^2\big) \nonumber \\
  & \quad + \frac{3L^2_f\mu_1^2d_1^2}{b\eta_t} \!+\! \frac{2\alpha_{t+1}^2\delta^2}{b\eta_t}, \nonumber
 \end{align}
 where the second inequality is due to $0<\alpha_{t+1}\leq 1$.
 By a similar way, we obtain
 \begin{align}
  & \frac{1}{\eta_t}\mathbb{E} \|\nabla_y f_{\mu_2}(x_{t+1},y_{t+1}) - w_{t+1}\|^2
  -  \frac{1}{\eta_{t-1}}\mathbb{E} \|\nabla_y f_{\mu_2}(x_t,y_t) - w_t\|^2  \\
  & \leq \big(\frac{1}{\eta_t} \!-\! \frac{1}{\eta_{t-1}} - c_2\eta_t\big)\mathbb{E} \|\nabla_y f_{\mu_2}(x_t,y_t) -w_t\|^2 \!+\! \frac{6d_2L^2_f\eta_{t}}{b}\mathbb{E}\big(\|\tilde{x}_{t+1}-x_t\|^2 + \|\tilde{y}_{t+1}-y_t\|^2\big) \nonumber \\
  & \quad + \frac{3L^2_f\mu_2^2d_2^2}{b\eta_t} \!+\! \frac{2\beta_{t+1}^2\delta^2}{b\eta_t}. \nonumber
 \end{align}
By $\eta_t = \frac{k}{(m+t)^{1/3}}$, we have
 \begin{align}
  \frac{1}{\eta_t} - \frac{1}{\eta_{t-1}} &= \frac{1}{k}\big( (m+t)^{\frac{1}{3}} - (m+t-1)^{\frac{1}{3}}\big) \nonumber \\
  & \leq \frac{1}{3k(m+t-1)^{2/3}} \leq \frac{1}{3k\big(m/2+t\big)^{2/3}} \nonumber \\
  & \leq \frac{2^{2/3}}{3k(m+t)^{2/3}} = \frac{2^{2/3}}{3k^3}\frac{k^2}{(m/2+t)^{2/3}}= \frac{2^{2/3}}{3k^3}\eta_t^2 \leq \frac{2}{3k^3}\eta_t,
 \end{align}
 where the first inequality holds by the concavity of function $f(x)=x^{1/3}$, \emph{i.e.}, $(x+y)^{1/3}\leq x^{1/3} + \frac{y}{3x^{2/3}}$; the second inequality is due to $m\geq 2$, and
 the last inequality is due to $0<\eta_t\leq 1$.
 Let $c_1 \geq \frac{2}{3k^3} + \frac{9\tau^2}{4}$, we have
 \begin{align} \label{eq:G1}
  & \frac{1}{\eta_t}\mathbb{E} \|\nabla_x f_{\mu_1}(x_{t+1},y_{t+1}) - v_{t+1}\|^2 - \frac{1}{\eta_{t-1}}\mathbb{E} \|\nabla_x f_{\mu_1}(x_t,y_t) - v_t\|^2  \\
  & \leq -\frac{9\tau^2}{4}\eta_t\mathbb{E} \|\nabla_x f_{\mu_1}(x_t,y_t) -v_t\|^2 + \frac{6d_1L^2_f\eta_{t}}{b}\mathbb{E}\big(\|\tilde{x}_{t+1}-x_t\|^2 + \|\tilde{y}_{t+1}-y_t\|^2\big) \nonumber\\
  &\quad + \frac{3L^2_f\mu_1^2d_1^2}{b\eta_t} + \frac{2\alpha_{t+1}^2\delta^2}{b\eta_t}. \nonumber
 \end{align}
 Let $c_2 \geq \frac{2}{3k^3} + \frac{625\tilde{d}L^2_f}{3b}$ with $\tilde{d}=d_1+d_2$, we have
  \begin{align} \label{eq:G2}
  & \frac{1}{\eta_t}\mathbb{E} \|\nabla_y f_{\mu_2}(x_{t+1},y_{t+1}) - w_{t+1}\|^2
  -  \frac{1}{\eta_{t-1}}\mathbb{E} \|\nabla_y f_{\mu_2}(x_t,y_t) - w_t\|^2  \\
  & \leq - \frac{625\tilde{d}L^2_f}{3b}\eta_t\mathbb{E} \|\nabla_y f_{\mu_2}(x_t,y_t) -w_t\|^2 + \frac{6d_2L^2_f\eta_{t}}{b}\mathbb{E}\big(\|\tilde{x}_{t+1}-x_t\|^2 + \|\tilde{y}_{t+1}-y_t\|^2\big) \nonumber \\
  &\quad + \frac{3L^2_f\mu_2^2d_2^2}{b\eta_t} + \frac{2\beta_{t+1}^2\delta^2}{b\eta_t}. \nonumber
 \end{align}
 According to Lemma \ref{lem:E1}, we have
 \begin{align} 
     &\|y_{t+1} - y^*(x_{t+1})\|^2 \nonumber\\
     & \leq (1-\frac{\eta_t\tau\lambda}{4})\|y_t -y^*(x_t)\|^2 - \frac{3\eta_t}{4} \|\tilde{y}_{t+1}-y_t\|^2
     + \frac{25\eta_t\lambda}{6\tau}  \|\nabla_y f(x_t,y_t)-w_t\|^2 \nonumber \\
     & \quad + \frac{25\kappa_y^2\eta_t}{6\tau\lambda}\|x_t -\tilde{x}_{t+1}\|^2 \nonumber \\
     & = (1-\frac{\eta_t\tau\lambda}{4})\|y_t -y^*(x_t)\|^2 -\frac{3\eta_t}{4} \|\tilde{y}_{t+1}-y_t\|^2 + \frac{25\kappa_y^2\eta_t}{6\tau\lambda}\|x_t - \tilde{x}_{t+1}\|^2 \nonumber \\
     & \quad + \frac{25\eta_t\lambda}{6\tau}  \|\nabla_y f(x_t,y_t) -\nabla_y f_{\mu_2}(x_t,y_t) + \nabla_y f_{\mu_2}(x_t,y_t) -w_t\|^2  \nonumber \\
     & \leq (1-\frac{\eta_t\tau\lambda}{4})\|y_t -y^*(x_t)\|^2 -\frac{3\eta_t}{4} \|\tilde{y}_{t+1}-y_t\|^2 +  \frac{25\kappa_y^2\eta_t}{6\tau\lambda}\|x_t - \tilde{x}_{t+1}\|^2\nonumber \\
     & \quad + \frac{25\lambda\mu^2_2L^2_fd^2_2\eta_t}{12\tau} + \frac{25\eta_t\lambda}{3\tau}\|\nabla_y f_{\mu_2}(x_t,y_t) -w_t\|^2, \nonumber
\end{align}
 where the last inequality is due to Young's inequality and Lemma \ref{lem:A3}.
 Thus, we have
\begin{align} \label{eq:G3}
 &\|y_{t+1} - y^*(x_{t+1})\|^2 - \|y_t -y^*(x_t)\|^2 \nonumber \\
 & \leq -\frac{\eta_t\tau\lambda}{4}\|y_t -y^*(x_t)\|^2 -\frac{3\eta_t}{4} \|\tilde{y}_{t+1}-y_t\|^2 +  \frac{25\kappa_y^2\eta_t}{6\tau\lambda}\|x_t - \tilde{x}_{t+1}\|^2\nonumber \\
 & \quad + \frac{25\lambda\mu^2_2L^2_fd^2_2\eta_t}{12\tau} + \frac{25\eta_t\lambda}{3\tau}\|\nabla_y f_{\mu_2}(x_t,y_t) -w_t\|^2.
\end{align}

 Next, we define a \emph{Lyapunov} function (i.e., potential function) $\Phi_t$, for any $t\geq 1$
 \begin{align}
 \Phi_t & =  \mathbb{E}\big[ F_{\mu_1}(x_t) + \frac{ 25\gamma\tilde{d}L^2_f}{\lambda\tau b}\|y_t-y^*(x_t)\|^2 + \frac{\gamma}{\tau^2\eta_{t-1}}\|\nabla_x f_{\mu_1}(x_t,y_t)-v_t\|^2 \nonumber \\
 & \quad + \frac{\gamma}{\tau^2\eta_{t-1}}\|\nabla_y f_{\mu_2}(x_t,y_t)-w_t\|^2 \big]. \nonumber
 \end{align}
 By using Lemma \ref{lem:D1}, we have
 \begin{align}
 &\Phi_{t+1} - \Phi_t \nonumber \\
 &\! = \mathbb{E}\big[F_{\mu_1}(x_{t+1}) - F_{\mu_1}(x_t)\big] + \frac{25\tilde{d}L^2_f\gamma}{\lambda\tau}
 \big( \mathbb{E}\|y_{t+1}-y^*(x_{t+1})\|^2 - \mathbb{E}\|y_t-y^*(x_t)\|^2 \big) \nonumber \\
 & \ + \frac{\gamma}{\tau^2}\big(\frac{1}{\eta_t}\mathbb{E}\|\nabla_x f_{\mu_1}(x_{t+1},y_{t+1})-v_{t+1}\|^2 - \frac{1}{\eta_{t-1}}\mathbb{E}\|\nabla_x f_{\mu_1}(x_t,y_t)-v_t\|^2 \big) \nonumber \\
 & \ + \frac{\gamma}{\tau^2}\big(\frac{1}{\eta_t}\mathbb{E}\|\nabla_y f_{\mu_2}(x_{t+1},y_{t+1})-w_{t+1}\|^2 - \frac{1}{\eta_{t-1}}\mathbb{E}\|\nabla_y f_{\mu_2}(x_t,y_t)-w_t\|^2 \big) \nonumber \\
 & \! \leq -\frac{\eta_t}{2\gamma}\mathbb{E}\|\tilde{x}_{t+1}-x_t\|^2 + 6\eta_t\gamma L_f^2\mathbb{E}\|y^*(x_t)-y_t\|^2 + 2\eta_t\gamma\mathbb{E}\|\nabla_xf_{\mu_1}(x_t,y_t) -v_t\|^2 +3\eta_t\gamma\mu_1^2d_1^2L_f^2  \nonumber \\
 & \ + \frac{25\tilde{d}L^2_f\gamma}{b\lambda\tau} \big(\!-\!\frac{\eta_t\tau\lambda}{4}\mathbb{E}\|y_t \!-\!y^*(x_t)\|^2 \!-\!\frac{3\eta_t}{4} \mathbb{E}\|\tilde{y}_{t+1}\!-\!y_t\|^2 \!+\! \frac{25\eta_t\lambda}{3\tau} \mathbb{E} \|\nabla_y f_{\mu_2}(x_t,y_t)\!-\!w_t\|^2  \nonumber \\
 & \ + \frac{25\lambda\mu^2_2L^2_fd^2_2\eta_t}{12\tau} + \frac{25\kappa_y^2\eta_t}{6\tau\lambda}\mathbb{E}\|x_t - \tilde{x}_{t+1}\|^2 \big) \!-\! \frac{9\gamma\eta_t}{4} \mathbb{E} \|\nabla_x f_{\mu_1}(x_t,y_t) -v_t\|^2  \nonumber \\
 & \ + \frac{6d_1L^2_f\eta_{t}\gamma}{b\tau^2}\big(\mathbb{E}\|\tilde{x}_{t+1}-x_t\|^2+ \mathbb{E}\|\tilde{y}_{t+1}-y_t\|^2\big) + \frac{6d_2L^2_f\eta_{t}\gamma}{b\tau^2}\big(\mathbb{E}\|\tilde{x}_{t+1}-x_t\|^2
 + \mathbb{E}\|\tilde{y}_{t+1}-y_t\|^2\big)  \nonumber \\
 & \ -\frac{625\tilde{d}L^2_f\gamma\eta_t}{3b\tau^2}\mathbb{E}\|\nabla_y f_{\mu_2}(x_t,y_t) -w_t\|^2 + \frac{3L^2_f\mu_1^2d_1^2\gamma}{b\eta_t\tau^2} \!+\! \frac{2\alpha_{t+1}^2\delta^2\gamma}{b\eta_t\tau^2}  + \frac{3L^2_f\mu_2^2d_2^2\gamma}{b\eta_t\tau^2} \!+\! \frac{2\beta_{t+1}^2\delta^2\gamma}{b\eta_t\tau^2} \nonumber \\
 & \!\leq -\frac{\gamma L_f^2\eta_t}{4}\mathbb{E}\|y^*(x_t)-y_t\|^2 \!-\! \frac{\gamma\eta_t}{4}\mathbb{E}\|\nabla_xf_{\mu_1}(x_t,y_t) -v_t\|^2 \!+\! 3\mu_1^2d_1^2L_f^2\eta_t\gamma \!+\! \frac{625\tilde{d}d^2_2L^4_f\mu^2_2\eta_t\gamma}{12b\tau^2} \nonumber \\
 & \ + \frac{3L^2_f\mu_1^2d_1^2\gamma}{b\eta_t\tau^2} + \frac{2\alpha_{t+1}^2\delta^2\gamma}{b\eta_t\tau^2} + \frac{3L^2_f\mu_2^2d_2^2\gamma}{b\eta_t\tau^2} \!+\! \frac{2\beta_{t+1}^2\delta^2\gamma}{b\eta_t\tau^2} \nonumber \\
 & \ - \big(\frac{75\tilde{d}L^2_f\gamma}{4b\lambda\tau} - \frac{6\tilde{d}L^2_f\gamma}{b\tau^2}\big)\eta_t\mathbb{E}\|\tilde{y}_{t+1}-y_t\|^2 -\big(\frac{1}{2\gamma} -\frac{6\tilde{d}L^2_f\gamma}{b\tau^2} - \frac{625\tilde{d}L^2_f\kappa_y^2\gamma}{6b\lambda^2\tau^2}\big)\eta_t\mathbb{E}\|\tilde{x}_{t+1}-x_t\|^2 \nonumber \\
 & \!\leq -\frac{\gamma L_f^2\eta_t}{4}\mathbb{E}\|y^*(x_t)-y_t\|^2 -\frac{\gamma\eta_t}{4}\mathbb{E}\|\nabla_xf_{\mu_1}(x_t,y_t) -v_t\|^2 -\frac{\eta_t}{4\gamma}\mathbb{E}\|\tilde{x}_{t+1}-x_t\|^2 + 3\mu_1^2d_1^2L_f^2\eta_t\gamma  \nonumber \\
 & \ + \frac{625\tilde{d}d^2_2L^4_f\mu^2_2\eta_t\gamma}{12b\tau^2}+ \frac{3L^2_f\mu_1^2d_1^2\gamma}{b\eta_t\tau^2} \!+\! \frac{2\alpha_{t+1}^2\delta^2\gamma}{b\eta_t\tau^2} + \frac{3L^2_f\mu_2^2d_2^2\gamma}{b\eta_t\tau^2} \!+\! \frac{2\beta_{t+1}^2\delta^2\gamma}{b\eta_t\tau^2},
 \end{align}
 where the first inequality holds by combining the above inequalities \eqref{eq:G1}, \eqref{eq:G2} and \eqref{eq:G3},
  and the second inequality is due to $1\leq b\leq \tilde{d}$
 and the last inequality is due to $0< \gamma \leq \frac{\lambda\tau^2}{2L_f}\sqrt{\frac{6b/\tilde{d}}{36 \lambda^2 + 625\kappa_y^2}}$
 and $\lambda\leq \frac{75\tau}{24}$.
 Thus, we have
 \begin{align} \label{eq:G4}
 & \frac{ L_f^2\eta_t}{4}\mathbb{E}\|y^*(x_t)-y_t\|^2 + \frac{\eta_t}{4}\mathbb{E}\|\nabla_xf_{\mu_1}(x_t,y_t) -v_t\|^2
 + \frac{\eta_t}{4\gamma^2}\mathbb{E}\|\tilde{x}_{t+1}-x_t\|^2  \\
 & \leq \frac{\Phi_t - \Phi_{t+1}}{\gamma} + 3\mu_1^2d_1^2L_f^2\eta_t + \frac{625\tilde{d}d^2_2L^4_f\mu^2_2\eta_t}{12b\tau^2} + \frac{3L^2_f\mu_1^2d_1^2}{b\eta_t\tau^2} \!+\! \frac{2\alpha_{t+1}^2\delta^2}{b\eta_t\tau^2} \!+\! \frac{3L^2_f\mu_2^2d_2^2}{b\eta_t\tau^2} + \frac{2\beta_{t+1}^2\delta^2}{b\eta_t\tau^2}. \nonumber
 \end{align}
Since $\inf_{x\in \mathcal{X}} F(x) =F^*$, we have
$\inf_{x\in \mathcal{X}} F_{\mu_1}(x) =\inf_{x\in \mathcal{X}} \mathbb{E}_{u_1\sim U_B}[F(x+\mu_1 u_1)] =
\inf_{x\in \mathcal{X}} \frac{1}{V}\int_{B} \\ F(x+\mu_1 u_1)du_1 \geq \frac{1}{V}\int_{B}\inf_{x\in \mathcal{X}}F(x+\mu_1 u_1)du_1=F^*$,
where $V$ denotes volume of the unit ball $B$.

Taking average over $t=1,2,\cdots,T$ on both sides of \eqref{eq:G4}, we have:
\begin{align}
 & \frac{1}{T} \sum_{t=1}^T \big( \frac{ L_f^2\eta_t}{4}\mathbb{E}\|y^*(x_t)-y_t\|^2 + \frac{\eta_t}{4}\mathbb{E}\|\nabla_xf_{\mu_1}(x_t,y_t) -v_t\|^2 +\frac{\eta_t}{4\gamma^2}\mathbb{E}\|\tilde{x}_{t+1}-x_t\|^2 \big) \nonumber \\
 & \leq  \sum_{t=1}^T \frac{\Phi_t - \Phi_{t+1}}{T\gamma} + \frac{1}{T}\sum_{t=1}^T\big( 3\mu_1^2d_1^2L_f^2\eta_t + \frac{625\tilde{d}d^2_2L^4_f\mu^2_2\eta_t}{12b\tau^2} + \frac{3L^2_f\mu_1^2d_1^2}{b\eta_t\tau^2} \nonumber \\
 & \quad + \frac{2\alpha_{t+1}^2\delta^2}{b\eta_t\tau^2} + \frac{3L^2_f\mu_2^2d_2^2}{b\eta_t\tau^2} + \frac{2\beta_{t+1}^2\delta^2}{b\eta_t\tau^2}\big). \nonumber
\end{align}
Let $\Delta_1=\|y_1-y^*(x_1)\|^2$, we have
\begin{align} \label{eq:G5}
 \Phi_1 &= F_{\mu_1}(x_1) + \frac{ 25\gamma\tilde{d}L^2_f}{\lambda\tau b}\|y_1-y^*(x_1)\|^2 +\frac{\gamma}{\eta_0\tau^2}\mathbb{E}\|\nabla_x f_{\mu_1}(x_1,y_1)-v_1\|^2 \nonumber \\
 &\quad + \frac{\gamma}{\eta_{0}\tau^2}\mathbb{E}\|\nabla_y f_{\mu_2}(x_1,y_1)-w_1\|^2 \nonumber \\
 & = F_{\mu_1}(x_1)+ \frac{ 25\gamma\tilde{d}L^2_f}{\lambda\tau b}\|y_1-y^*(x_1)\|^2+\frac{\gamma}{\eta_0\tau^2}\mathbb{E}\|\nabla_x f_{\mu_1}(x_1,y_1)-\hat{\nabla}_x f(x_1,y_1;\mathcal{B}_1)\|^2 \nonumber \\
 & \quad + \frac{\gamma}{\eta_{0}\tau^2}\mathbb{E}\|\nabla_y f_{\mu_2}(x_1,y_1)-\hat{\nabla}_y f(x_1,y_1;\mathcal{B}_1)\|^2 \nonumber \\
 & \leq F_{\mu_1}(x_1) + \frac{ 25\gamma\tilde{d}L^2_f}{\lambda\tau b}\Delta_1 +\frac{2\gamma\delta^2}{b\eta_0\tau^2},
\end{align}
where the last inequality holds by Assumption \ref{ass:2}.
Since $\eta_t$ is decreasing, i.e., $\eta_T^{-1} \geq \eta_t^{-1}$ for any $0\leq t\leq T$, we have
 \begin{align}
 & \frac{1}{T} \sum_{t=1}^T \big( \frac{L_f^2}{4}\mathbb{E}\|y^*(x_t)-y_t\|^2 + \frac{1}{4}\mathbb{E}\|\nabla_xf_{\mu_1}(x_t,y_t) -v_t\|^2 +\frac{1}{4\gamma^2}\mathbb{E}\|\tilde{x}_{t+1}-x_t\|^2 \big) \nonumber \\
 & \leq  \frac{1}{T\gamma\eta_T} \sum_{t=1}^T\big(\Phi_t - \Phi_{t+1}\big)+ \frac{1}{T\eta_T}\sum_{t=1}^T\big( 3\mu_1^2d_1^2L_f^2\eta_t + \frac{625\tilde{d}d^2_2L^4_f\mu^2_2\eta_t}{12b\tau^2} + \frac{3L^2_f\mu_1^2d_1^2}{b\eta_t\tau^2} \!+\! \frac{2\alpha_{t+1}^2\delta^2}{b\eta_t\tau^2} \nonumber \\
 & \quad  + \frac{3L^2_f\mu_2^2d_2^2}{b\eta_t\tau^2} + \frac{2\beta_{t+1}^2\delta^2}{b\eta_t\tau^2}\big) \nonumber \\
 & \leq \frac{1}{T\gamma\eta_T} \big( F_{\mu_1}(x_1) - F^* + \frac{ 25\gamma\tilde{d}L^2_f}{\lambda\tau b}\Delta_1+\frac{2\delta^2\gamma}{b\tau^2\eta_0} \big)
 + \frac{1}{T\eta_T}\sum_{t=1}^T\big( 3\mu_1^2d_1^2L_f^2\eta_t
 + \frac{625\tilde{d}d^2_2L^4_f\mu^2_2\eta_t}{12b\tau^2}
 \nonumber \\
 & \quad + \frac{3L^2_f\mu_1^2d_1^2}{b\eta_t\tau^2}+ \frac{2\alpha_{t+1}^2\delta^2}{b\eta_t\tau^2}+ \frac{3L^2_f\mu_2^2d_2^2}{b\eta_t\tau^2} \!+\! \frac{2\beta_{t+1}^2\delta^2}{b\eta_t\tau^2}\big)
 \nonumber \\
 & = \frac{F_{\mu_1}(x_1) - F^*}{T\gamma\eta_T} + \frac{ 25\tilde{d}L^2_f}{T\eta_T\lambda\tau b}\Delta_1+\frac{2\delta^2}{Tb\tau^2\eta_T\eta_0} + \frac{36\tau^2\mu_1^2d_1^2L_f^2 + 625\tilde{d}d^2_2L^4_f\mu^2_2}{12b\tau^2T\eta_T}\sum_{t=1}^T\eta_t  \nonumber \\
 & \quad + \frac{3L^2_f\big(\mu_1^2 d_1^2 + \mu_2^2d_2^2\big)}{Tb\tau^2\eta_T}\sum_{t=1}^T\frac{1}{\eta_t}+ \frac{2(c_1^2+c_2^2)\delta^2}{Tb\tau^2\eta_T}\sum_{t=1}^T\eta_t^3
 \nonumber \\
 & \leq \frac{F_{\mu_1}(x_1) - F^*}{T\gamma\eta_T} + \frac{ 25\tilde{d}L^2_f}{T\eta_T\lambda\tau b}\Delta_1 +\frac{2\delta^2}{Tb\tau^2\eta_T\eta_0} + \frac{36\tau^2\mu_1^2d_1^2L_f^2 + 625\tilde{d}d^2_2L^4_f\mu^2_2}{12b\tau^2T\eta_T}\int^T_1\frac{k}{(m+t)^{1/3}}dt  \nonumber \\
 &\quad + \frac{3L^2_f\big(\mu_1^2 d_1^2 + \mu_2^2d_2^2\big)}{Tb\tau^2\eta_T}\int^T_1\frac{(m+t)^{1/3}}{k}dt + \frac{2(c_1^2+c_2^2)\delta^2}{Tb\tau^2\eta_T}\int^T_1\frac{k^3}{m+t} dt\nonumber \\
 & \leq \frac{F_{\mu_1}(x_1) - F^*}{T\gamma\eta_T} + \frac{ 25\tilde{d}L^2_f}{T\eta_T\lambda\tau b}\Delta_1 +\frac{2\delta^2}{Tb\tau^2\eta_T\eta_0} + \frac{36\tau^2\mu_1^2d_1^2L_f^2k + 625\tilde{d}d^2_2L^4_f\mu^2_2k}{8b\tau^2T\eta_T}(m+T)^{2/3} \nonumber \\
 &\quad + \frac{9L^2_f\big(\mu_1^2 d_1^2 + \mu_2^2d_2^2\big)}{4Tb\tau^2\eta_T k}(m+T)^{4/3}  + \frac{2(c_1^2+c_2^2)\delta^2 k^3}{Tb\tau^2\eta_T}\ln(m+T) \nonumber \\
 & \leq \frac{F_{\mu_1}(x_1) - F^*}{T\gamma\eta_T}+ \frac{ 25\tilde{d}L^2_f}{T\eta_T\lambda\tau b}\Delta_1 +\frac{2\delta^2}{Tb\tau^2\eta_T\eta_0} + \frac{36\tau^2L_f^2k + 625L^4_fk}{8b\tau^2T\eta_T}(m+T)^{-2/3} + \frac{9L^2_f}{4Tb\tau^2\eta_T k} \nonumber \\
 & \quad + \frac{2(c_1^2+c_2^2)\delta^2 k^3}{Tb\tau^2\eta_T}\ln(m+T) \nonumber \\
 & = \big(\frac{F_{\mu_1}(x_1) - F^*}{T\gamma k} + \frac{ 25\tilde{d}L^2_f}{T k\lambda\tau b}\Delta_1 +\frac{2m^{1/3}\delta^2}{Tb\tau^2 k^2}\big)(m+T)^{1/3} + \frac{36\tau^2L_f^2 + 625L^4_f}{8b\tau^2T}(m+T)^{-1/3}  \nonumber \\
& \quad + \frac{9L^2_f}{4Tb\tau^2 k^2}(m+T)^{1/3} + \frac{2(c_1^2+c_2^2)\delta^2 k^2}{Tb\tau^2}\ln(m+T)(m+T)^{1/3},
\end{align}
where the second inequality holds by the above inequality \eqref{eq:G5}, and the last inequality is due to $0<\mu_1\leq \frac{1}{d_1(m+T)^{2/3}}$
and $0<\mu_2\leq \frac{1}{\tilde{d}^{1/2}d_2(m+T)^{2/3}}$. Let $M' =  \frac{F_{\mu_1}(x_1) - F^*}{\gamma k} + \frac{ 25\tilde{d}L^2_f}{k\lambda\tau b}\Delta_1+ \frac{2m^{1/3}\delta^2}{b\tau^2 k^2} + \frac{36\tau^2L_f^2 + 625L^4_f}{8b\tau^2}(m+T)^{-2/3} + \frac{9L^2_f}{4b\tau^2 k^2} + \frac{2(c_1^2+c_2^2)\delta^2 k^2}{b\tau^2}\ln(m+T)$,
we have
\begin{align}
 \frac{1}{T} \sum_{t=1}^T \big( \frac{ L_f^2}{4}\mathbb{E}\|y^*(x_t)\!-\!y_t\|^2 \!+\! \frac{1}{4}\mathbb{E}\|\nabla_xf_{\mu_1}(x_t,y_t) \!-\!v_t\|^2 \!+\! \frac{1}{4\gamma^2}\mathbb{E}\|\tilde{x}_{t+1}\!-\!x_t\|^2 \big) \leq \frac{M'}{T}(m+T)^{1/3}.
\end{align}
According to Jensen's inequality, we have
\begin{align}
 &  \frac{1}{T} \sum_{t=1}^T \big( \frac{ L_f}{2}\mathbb{E}\|y^*(x_t)-y_t\| + \frac{1}{2}\mathbb{E}\|\nabla_xf_{\mu_1}(x_t,y_t) -v_t\| +\frac{1}{2\gamma}\mathbb{E}\|\tilde{x}_{t+1}-x_t\| \big) \nonumber \\
 & \leq \big( \frac{3}{T} \sum_{t=1}^T \big( \frac{ L_f^2}{4}\mathbb{E}\|y^*(x_t)-y_t\|^2 + \frac{1}{4}\mathbb{E}\|\nabla_xf_{\mu_1}(x_t,y_t) -v_t\|^2 +\frac{1}{4\gamma^2}\mathbb{E}\|\tilde{x}_{t+1}-x_t\|^2 \big)\big)^{1/2} \nonumber \\
 & \leq \frac{\sqrt{3M'}}{T^{1/2}}(m+T)^{1/6} \leq \frac{\sqrt{3M'}m^{1/6}}{T^{1/2}} + \frac{\sqrt{3M'}}{T^{1/3}},
\end{align}
where the last inequality is due to $(a+b)^{1/6} \leq a^{1/6} + b^{1/6}$.
Thus we obtain
\begin{align}
 & \frac{1}{T} \sum_{t=1}^T \mathbb{E}  \big[ L_f\|y^*(x_t)-y_t\| + \|\nabla_xf_{\mu_1}(x_t,y_t) -v_t\| +\frac{1}{\gamma}\|\tilde{x}_{t+1}-x_t\| \big] \nonumber \\
 & \leq \frac{2\sqrt{3M'}m^{1/6}}{T^{1/2}} + \frac{2\sqrt{3M'}}{T^{1/3}}. \nonumber
\end{align}

According to Lemma \ref{lem:A3}, we have $\|\nabla_x f_{\mu_1}(x_t,y_t) - v_t\|=\|\nabla_x f_{\mu_1}(x_t,y_t) - \nabla_x f(x_t,y_t)
+ \nabla_x f(x_t,y_t)- v_t\| \geq \|\nabla_x f(x_t,y_t)- v_t\| - \|\nabla_x f_{\mu_1}(x_t,y_t) - \nabla_x f(x_t,y_t)\| \geq \|\nabla_x f(x_t,y_t)- v_t\|
- \frac{\mu_1 L_f d_1}{2}$.
Thus, we have
\begin{align}
 & \frac{1}{T} \sum_{t=1}^T \mathbb{E} \big[L_f\|y^*(x_t)-y_t\| + \|\nabla_xf(x_t,y_t) -v_t\| +\frac{1}{\gamma}\|\tilde{x}_{t+1}-x_t\|\big] \nonumber \\
 & \leq \frac{1}{T} \sum_{t=1}^T \mathbb{E}  \big[L_f\|y^*(x_t)-y_t\| + \|\nabla_xf_{\mu_1}(x_t,y_t) -v_t\| +\frac{\mu_1 L_f d_1}{2} +\frac{1}{\gamma}\|\tilde{x}_{t+1}-x_t\|\big] \nonumber \\
 & \leq  \frac{2\sqrt{3M'}m^{1/6}}{T^{1/2}} + \frac{2\sqrt{3M'}}{T^{1/3}} + \frac{\mu_1 L_f d_1}{2} \nonumber \\
 & \leq  \frac{2\sqrt{3M'}m^{1/6}}{T^{1/2}} + \frac{2\sqrt{3M'}}{T^{1/3}} + \frac{L_f}{2(m+T)^{2/3}},
\end{align}
where the last inequality is due to $0<\mu_1 \leq \frac{1}{d_1(m+T)^{2/3}}$. Then by using the above inequality \eqref{eq:MH}, we have 
\begin{align}
 \frac{1}{T} \sum_{t=1}^T \mathbb{E}  \|G_{\mathcal{X}}(x_t,\nabla F(x_t),\gamma)\| & \leq\frac{1}{T} \sum_{t=1}^T \mathbb{E}\big[ L_f\|y^*(x_t)\!-\!y_t\| \!+\! \|\nabla_xf(x_t,y_t) \!-\!v_t\| \!+\!\frac{1}{\gamma}\|\tilde{x}_{t+1}\!-\!x_t\|\big]
 \nonumber \\
 &\leq  \frac{2\sqrt{3M'}m^{1/6}}{T^{1/2}} + \frac{2\sqrt{3M'}}{T^{1/3}} + \frac{L_f}{2(m+T)^{2/3}}.
\end{align}

\end{proof}

\subsection{ Convergence Analysis of Acc-ZOMDA Algorithm for Unconstrained Minimax Optimization }
\label{Appendix:A4}
In this subsection, we study the convergence properties of our  Acc-ZOMDA algorithm for solving the black-box \textbf{unconstrained}
minimax problem \eqref{eq:2}, i.e., $\mathcal{X}= \mathbb{R}^{d_1}$
and $\mathcal{Y} = \mathbb{R}^{d_2}$ (or $\mathcal{Y} \subset \mathbb{R}^{d_2}$). The following convergence analysis builds on the common convergence metric
$\mathbb{E}\|\nabla F(x_t)\|$ used in \citep{lin2019gradient},
where $F(x)=\max_{y\in \mathcal{Y}}f(x,y)$.

\begin{lemma} \label{lem:D01}
Suppose the sequence $\{x_t,y_t\}_{t=1}^T$ be generated from Algorithm \ref{alg:2}.
When $\mathcal{X}=\mathbb{R}^{d_1}$, given $0<\gamma\leq \frac{1}{2\eta_t L_g}$,
we have
\begin{align}
  F_{\mu_1}(x_{t+1}) & \leq F_{\mu_1}(x_t)
  + 3\eta_t\gamma L_f^2 \|y_t - y^*(x_t)\|^2 + \frac{3\eta_t\gamma L_f^2d_1^2\mu_1^2}{2}\nonumber \\
  & \quad + \gamma\eta_t\|\nabla_x f_{\mu_1}(x_t,y_t)-v_t\|^2 - \frac{\gamma\eta_t}{2}\|\nabla F_{\mu_1}(x_t)\|^2 - \frac{\gamma\eta_t}{4}\|v_t\|^2.
\end{align}
\end{lemma}
\begin{proof}
 According to Lemma \ref{lem:1} and Lemma \ref{lem:A3}, the approximated function $F_{\mu_1}(x)$ has $L_g$-Lipschitz continuous gradient.
 Then we have
 \begin{align}
  &F_{\mu_1}(x_{t+1}) \nonumber \\
  &\leq F_{\mu_1}(x_t) - \gamma\eta_t\langle\nabla F_{\mu_1}(x_t),v_t\rangle + \frac{\gamma^2\eta_t^2L_g}{2}\|v_t\|^2   \\
  & = F_{\mu_1}(x_t) + \frac{\gamma\eta_t}{2}\|\nabla F_{\mu_1}(x_t)-v_t\|^2 - \frac{\gamma\eta_t}{2}\|\nabla F_{\mu_1}(x_t)\|^2
  + (\frac{\gamma^2\eta_t^2 L_g}{2}-\frac{\gamma\eta_t}{2})\|v_t\|^2 \nonumber \\
  & = F_{\mu_1}(x_t) + \frac{\gamma\eta_t}{2}\|\nabla F_{\mu_1}(x_t)-\nabla_x f_{\mu_1}(x_t,y_t) + \nabla_x f_{\mu_1}(x_t,y_t)-v_t\|^2
  - \frac{\gamma\eta_t}{2}\|\nabla F_{\mu_1}(x_t)\|^2 \nonumber \\
  & \quad + (\frac{\gamma^2\eta_t^2 L_g}{2}-\frac{\gamma\eta_t}{2})\|v_t\|^2 \nonumber \\
  & \leq F_{\mu_1}(x_t) + \gamma\eta_t\|\nabla F_{\mu_1}(x_t)-\nabla_x f_{\mu_1}(x_t,y_t)\|^2 + \gamma\eta_t\|\nabla_x f_{\mu_1}(x_t,y_t)-v_t\|^2
  \nonumber \\
  & \quad - \frac{\gamma\eta_t}{2}\|\nabla F_{\mu_1}(x_t)\|^2  + (\frac{\gamma^2\eta_t^2 L_g}{2}-\frac{\gamma\eta_t}{2})\|v_t\|^2 \nonumber \\
  & \leq  F_{\mu_1}(x_t) + \gamma\eta_t\|\nabla F_{\mu_1}(x_t)-\nabla_x f_{\mu_1}(x_t,y_t)\|^2 + \gamma\eta_t\|\nabla_x f_{\mu_1}(x_t,y_t)-v_t\|^2 \nonumber \\
  & \quad - \frac{\gamma\eta_t}{2}\|\nabla F_{\mu_1}(x_t)\|^2 - \frac{\gamma\eta_t}{4}\|v_t\|^2, \nonumber
 \end{align}
 where the last inequality is due to $0< \gamma \leq \frac{1}{2\eta_t L}$.

 Considering an upper bound of $\|\nabla F_{\mu_1}(x_t)-\nabla_x f_{\mu_1}(x_t,y_t)\|^2$, we have
 \begin{align}
  & \|\nabla F_{\mu_1}(x_t)-\nabla_x f_{\mu_1}(x_t,y_t)\|^2 \nonumber \\
  & = \|\nabla_x f_{\mu_1}(x_t,y^*(x_t)) - \nabla_x f_{\mu_1}(x_t,y_t)\|^2 \nonumber \\
  & = \|\nabla_x f_{\mu_1}(x_t,y^*(x_t)) - \nabla_xf(x_t,y^*(x_t)) + \nabla_xf(x_t,y^*(x_t)) - \nabla_xf(x_t,y_t)\nonumber \\
  & \quad + \nabla_xf(x_t,y_t) - \nabla_x f_{\mu_1}(x_t,y_t)\|^2 \nonumber \\
  & \leq 3 \|\nabla_x f_{\mu_1}(x_t,y^*(x_t)) - \nabla_xf(x_t,y^*(x_t))\|^2 + 3 \|\nabla_xf(x_t,y^*(x_t)) - \nabla_xf(x_t,y_t)\|^2 \nonumber \\
  & \quad + 3 \|\nabla_xf(x_t,y_t) - \nabla_x f_{\mu_1}(x_t,y_t)\|^2 \nonumber \\
  & \leq \frac{3L_f^2d_1^2\mu_1^2}{2} + 3L_f^2 \|y_t - y^*(x_t)\|^2,
 \end{align}
 the last inequality holds by Assumption \ref{ass:5}
 and Lemma \ref{lem:A3}, i.e., we have
 \begin{align}
  \|\nabla_xf_{\mu_1}(x_t,y^*(x_t)) -\nabla_xf(x_t,y^*(x_t))\| \leq \frac{L_fd_1\mu_1}{2},
  \ \|\nabla_xf(x_t,y_t)- \nabla_xf_{\mu_1}(x_t,y_t)\|\leq \frac{L_fd_1\mu_1}{2}, \nonumber
 \end{align}
 and
 \begin{align}
 \|\nabla_xf(x_t,y^*(x_t))-\nabla_xf(x_t,y_t) \| \leq \|\nabla f(x_t,y^*(x_t))-\nabla f(x_t,y_t) \|\leq L_f\|y_t-y^*(x_t)\|. \nonumber
 \end{align}
 Then we have
 \begin{align}
  F_{\mu_1}(x_{t+1}) & \leq F_{\mu_1}(x_t)
  + \frac{3\eta_t\gamma L_f^2d_1^2\mu_1^2}{2} + 3\eta_t\gamma L_f^2 \|y_t - y^*(x_t)\|^2\nonumber \\
  & \quad + \gamma\eta_t\|\nabla_x f_{\mu_1}(x_t,y_t)-v_t\|^2 - \frac{\gamma\eta_t}{2}\|\nabla F_{\mu_1}(x_t)\|^2
  - \frac{\gamma\eta_t}{4}\|v_t\|^2.
 \end{align}

\end{proof}

\begin{lemma} \label{lem:E01}
Suppose the sequence $\{x_t,y_t\}_{t=1}^T$ be generated from Algorithm \ref{alg:2}. Under the above assumptions, and set $0< \eta_t\leq 1$
and $0<\lambda\leq \frac{1}{6L_f}$, we have
\begin{align}
\|y_{t+1} - y^*(x_{t+1})\|^2 &\leq (1-\frac{\eta_t\tau\lambda}{4})\|y_t -y^*(x_t)\|^2 -\frac{3\eta_t}{4} \|\tilde{y}_{t+1}-y_t\|^2 \nonumber \\
& \quad + \frac{25\eta_t\lambda}{6\tau}  \|\nabla_y f(x_t,y_t)-w_t\|^2 +  \frac{25\kappa_y^2\gamma^2\eta_t}{6\tau\lambda}\|v_t\|^2,
\end{align}
where $\kappa_y = L_f/\tau$.
\end{lemma}
\begin{proof}
This proof is similar to the proof of Lemma \ref{lem:E1}.
\end{proof}

\begin{lemma} \label{lem:F01}
 Suppose the zeroth-order stochastic gradients $\{v_t,w_t\}_{t=1}^T$ be generated from Algorithm \ref{alg:2}, we have
\begin{align}
&\mathbb{E} \|\nabla_x f_{\mu_1}(x_{t+1},y_{t+1}) - v_{t+1}\|^2 \nonumber \\
& \leq (1-\alpha_{t+1})^2 \mathbb{E} \|\nabla_x f_{\mu_1}(x_t,y_t) -v_t\|^2
 + \frac{3(1-\alpha_{t+1})^2L^2_f\mu_1^2d_1^2}{b} \nonumber \\
 & \quad +\frac{6d_1L^2_f(1-\alpha_{t+1})^2\eta^2_t}{b}\big(\gamma^2\mathbb{E}\|v_t\|^2
 + \mathbb{E}\|\tilde{y}_{t+1}-y_t\|^2\big) + \frac{2\alpha_{t+1}^2\delta^2}{b}.
\end{align}
\begin{align}
&\mathbb{E} \|\nabla_y f_{\mu_2}(x_{t+1},y_{t+1}) - w_{t+1}\|^2 \nonumber \\
&\leq (1-\beta_{t+1})^2 \mathbb{E} \|\nabla_y f_{\mu_2}(x_t,y_t) -w_t\|^2
 + \frac{3(1-\beta_{t+1})^2L^2_f\mu_2^2d_2^2}{b} \nonumber \\
 & \quad +\frac{6d_2L^2_f(1-\beta_{t+1})^2\eta^2_t}{b}\big(\gamma^2\mathbb{E}\|v_t\|^2
 + \mathbb{E}\|\tilde{y}_{t+1}-y_t\|^2\big) + \frac{2\beta_{t+1}^2\delta^2}{b}.
\end{align}
\end{lemma}
\begin{proof}
This proof is similar to the proof of Lemma \ref{lem:F1}.
\end{proof}

\begin{theorem} \label{th:A02}
 (Restatement of Theorem 7)
Suppose the sequence $\{x_t,y_t\}_{t=1}^T$ be generated from Algorithm \ref{alg:2}. When $\mathcal{X}=\mathbb{R}^{d_1}$, 
and let $\eta_t = \frac{k}{(m+t)^{1/3}}$
for all $t\geq 0$, $c_1 \geq \frac{2}{3k^3} + \frac{9\tau^2}{4}$ and $c_2 \geq \frac{2}{3k^3} + \frac{625\tilde{d}L^2_f}{3b}$,
$k>0$, $1\leq b\leq \tilde{d}$, $m\geq \max\big( 2, k^3, (c_1k)^3, (c_2k)^3\big)$,  $0<\lambda\leq \min\big(\frac{1}{6L_f},\frac{75\tau}{24}\big)$,
$0< \gamma \leq \min\big( \frac{\lambda\tau}{2L_f}\sqrt{\frac{6b/\tilde{d}}{36 \lambda^2 + 625\kappa_y^2}}, \frac{m^{1/3}}{2L_gk}\big)$, $0<\mu_1\leq \frac{1}{d_1(m+T)^{2/3}}$
and $0<\mu_2\leq \frac{1}{\tilde{d}^{1/2}d_2(m+T)^{2/3}}$,
we have
\begin{align}
 \frac{1}{T} \sum_{t=1}^T \mathbb{E}\|\nabla F(x_t)\| \leq  \frac{\sqrt{2M'}m^{1/6}}{T^{1/2}} + \frac{\sqrt{2M'}}{T^{1/3}} + \frac{L_f}{2(m+T)^{2/3}},
\end{align}
where $\Delta_1=\|y_1-y^*(x_1)\|^2$ and $M' = \frac{F_{\mu_1}(x_1) - F^*}{\gamma k} + \frac{ 25\tilde{d}L^2_f}{k\lambda\tau b}\Delta_1 +\frac{2m^{1/3}\delta^2}{b\tau^2 k^2} + \frac{36\tau^2L_f^2 + 625L^4_f}{8b\tau^2}(m+T)^{-2/3} + \frac{9L^2_f}{4b\tau^2 k^2} + \frac{2(c_1^2+c_2^2)\delta^2 k^2}{b\tau^2}\ln(m+T)$.
\end{theorem}
\begin{proof}
This proof is similar to the proof of Theorem \ref{th:A2}.
Following the above proof of Theorem \ref{th:A2},
let $c_1 \geq \frac{2}{3k^3} + \frac{9\tau^2}{4}$, we have
 \begin{align} \label{eq:G01}
  & \frac{1}{\eta_t}\mathbb{E} \|\nabla_x f_{\mu_1}(x_{t+1},y_{t+1}) - v_{t+1}\|^2 - \frac{1}{\eta_{t-1}}\mathbb{E} \|\nabla_x f_{\mu_1}(x_t,y_t) - v_t\|^2  \\
  & \!\leq\! -\frac{9\tau^2}{4}\eta_t\mathbb{E} \|\nabla_x f_{\mu_1}(x_t,y_t) -v_t\|^2 \!+\! \frac{6d_1L^2_f\eta_{t}}{b}\big(\gamma^2\mathbb{E}\|v_t\|^2 + \mathbb{E}\|\tilde{y}_{t+1}-y_t\|^2\big) + \frac{3L^2_f\mu_1^2d_1^2}{b\eta_t} \!+\! \frac{2\alpha_{t+1}^2\delta^2}{b\eta_t}. \nonumber
 \end{align}
 Similarly, let $c_2 \geq \frac{2}{3k^3} + \frac{625\tilde{d}L^2_f}{3b}$ with $\tilde{d}=d_1+d_2$, we  also have
  \begin{align} \label{eq:G02}
  & \frac{1}{\eta_t}\mathbb{E} \|\nabla_y f_{\mu_2}(x_{t+1},y_{t+1}) - w_{t+1}\|^2
  -  \frac{1}{\eta_{t-1}}\mathbb{E} \|\nabla_y f_{\mu_2}(x_t,y_t) - w_t\|^2  \\
  & \! \leq \! - \frac{625\tilde{d}L^2_f}{3b}\eta_t\mathbb{E} \|\nabla_y f_{\mu_2}(x_t,y_t) \!-\! w_t\|^2 \!+\! \frac{6d_2L^2_f\eta_{t}}{b}\big(\gamma^2\mathbb{E}\|v_t\|^2 \!+\! \mathbb{E}\|\tilde{y}_{t+1}\!-\!y_t\|^2\big) \!+\! \frac{3L^2_f\mu_2^2d_2^2}{b\eta_t} \!+\! \frac{2\beta_{t+1}^2\delta^2}{b\eta_t}. \nonumber
 \end{align}
 According to Lemma \ref{lem:E01}, we have
 \begin{align}
     \|y_{t+1} - y^*(x_{t+1})\|^2 & \leq (1-\frac{\eta_t\tau\lambda}{4})\|y_t \!-\!y^*(x_t)\|^2 \!-\! \frac{3\eta_t}{4} \|\tilde{y}_{t+1}\!-\!y_t\|^2
     \! + \! \frac{25\eta_t\lambda}{6\tau}  \|\nabla_y f(x_t,y_t)\!-\!w_t\|^2 \nonumber \\
     & \quad + \frac{25\kappa_y^2\gamma^2\eta_t}{6\tau\lambda}\|v_t\|^2 \nonumber \\
     & = (1-\frac{\eta_t\tau\lambda}{4})\|y_t -y^*(x_t)\|^2 -\frac{3\eta_t}{4} \|\tilde{y}_{t+1}-y_t\|^2  + \frac{25\kappa_y^2\gamma^2\eta_t}{6\tau\lambda}\|v_t\|^2  \nonumber \\
     & \quad + \frac{25\eta_t\lambda}{6\tau}  \|\nabla_y f(x_t,y_t) -\nabla_y f_{\mu_2}(x_t,y_t) + \nabla_y f_{\mu_2}(x_t,y_t) -w_t\|^2  \nonumber \\
     & \leq (1-\frac{\eta_t\tau\lambda}{4})\|y_t -y^*(x_t)\|^2 -\frac{3\eta_t}{4} \|\tilde{y}_{t+1}-y_t\|^2  + \frac{25\kappa_y^2\gamma^2\eta_t}{6\tau\lambda}\|v_t\|^2 \nonumber \\
     & \quad + \frac{25\lambda\mu^2_2L^2_fd^2_2\eta_t}{12\tau} + \frac{25\eta_t\lambda}{3\tau}\|\nabla_y f_{\mu_2}(x_t,y_t) -w_t\|^2,
\end{align}
where the last inequality is due to Young's inequality and Lemma \ref{lem:A3}.
Thus, we have
\begin{align} \label{eq:G03}
 & \|y_{t+1} - y^*(x_{t+1})\|^2 - \|y_t -y^*(x_t)\|^2 \nonumber \\
 & \leq -\frac{\eta_t\tau\lambda}{4}\|y_t -y^*(x_t)\|^2 -\frac{3\eta_t}{4} \|\tilde{y}_{t+1}-y_t\|^2  + \frac{25\kappa_y^2\gamma^2\eta_t}{6\tau\lambda}\|v_t\|^2 \nonumber \\
 & \quad + \frac{25\lambda\mu^2_2L^2_fd^2_2\eta_t}{12\tau} + \frac{25\eta_t\lambda}{3\tau}\|\nabla_y f_{\mu_2}(x_t,y_t) -w_t\|^2.
\end{align}

At the same time, we give the  \emph{Lyapunov} function $\Phi_t$ defined in the above proof of Theorem \ref{th:A2}, 
 \begin{align}
 \Phi_t & =  \mathbb{E} \big[ F_{\mu_1}(x_t) + \frac{ 25\gamma\tilde{d}L^2_f}{\lambda\tau b}\|y_t-y^*(x_t)\|^2 + \frac{\gamma}{\tau^2\eta_{t-1}}\|\nabla_x f_{\mu_1}(x_t,y_t)-v_t\|^2 \nonumber \\
 &\quad + \frac{\gamma}{\tau^2\eta_{t-1}}\|\nabla_y f_{\mu_2}(x_t,y_t)-w_t\|^2 \big] . \nonumber
 \end{align}
By using Lemma \ref{lem:D01}, we have
\begin{align}
 &\Phi_{t+1} - \Phi_t \nonumber \\
 & = \mathbb{E}\big[F_{\mu_1}(x_{t+1}) - F_{\mu_1}(x_t)\big] + \frac{25\tilde{d}L^2_f\gamma}{\lambda\tau}
 \big( \mathbb{E}\|y_{t+1}-y^*(x_{t+1})\|^2 - \mathbb{E}\|y_t-y^*(x_t)\|^2 \big) \nonumber \\
 & \quad +\frac{\gamma}{\tau^2}\big(\frac{1}{\eta_t}\mathbb{E} \|\nabla_x f_{\mu_1}(x_{t+1},y_{t+1})-v_{t+1}\|^2 - \frac{1}{\eta_{t-1}}\mathbb{E} \|\nabla_x f_{\mu_1}(x_t,y_t)-v_t\|^2 \big) \nonumber \\
 & \quad + \frac{\gamma}{\tau^2}\big(\frac{1}{\eta_t}\mathbb{E} \|\nabla_y f_{\mu_2}(x_{t+1},y_{t+1})-w_{t+1}\|^2 - \frac{1}{\eta_{t-1}}\mathbb{E} \|\nabla_y f_{\mu_2}(x_t,y_t)-w_t\|^2 \big) \nonumber \\
 & \leq  3\eta_t\gamma L_f^2 \mathbb{E}\|y_t - y^*(x_t)\|^2 \!+\! \frac{3\eta_t\gamma L_f^2d_1^2\mu_1^2}{2} \!+\! \gamma\eta_t\mathbb{E} \|\nabla_x f_{\mu_1}(x_t,y_t)-v_t\|^2 \!-\! \frac{\gamma\eta_t}{2}\mathbb{E}\|\nabla F_{\mu_1}(x_t)\|^2   \nonumber \\
 & \quad - \frac{\gamma\eta_t}{4}\mathbb{E}\|v_t\|^2 \!+\! \frac{25\tilde{d}L^2_f\gamma}{b\lambda\tau} \big(\! -\!\frac{\eta_t\tau\lambda}{4}\mathbb{E}\|y_t -y^*(x_t)\|^2 \!-\!\frac{3\eta_t}{4} \mathbb{E}\|\tilde{y}_{t+1}-y_t\|^2+ \frac{25\kappa_y^2\gamma^2\eta_t}{6\tau\lambda}\mathbb{E}\|v_t\|^2 \nonumber \\
 & \quad + \frac{25\eta_t\lambda}{3\tau}  \mathbb{E} \|\nabla_y f_{\mu_2}(x_t,y_t)-w_t\|^2 + \frac{25\lambda\mu^2_2L^2_fd^2_2\eta_t}{12\tau} \big) - \frac{9\gamma\eta_t}{4} \mathbb{E} \|\nabla_x f_{\mu_1}(x_t,y_t) -v_t\|^2  \nonumber \\
 & \quad + \frac{6d_1L^2_f\eta_{t}\gamma}{b\tau^2}\big(\gamma^2\mathbb{E}\|v_t\|^2  + \mathbb{E}\|\tilde{y}_{t+1}-y_t\|^2\big)-\frac{625\tilde{d}L^2_f\gamma}{3b\tau^2} \eta_t\mathbb{E}\|\nabla_y f_{\mu_2}(x_t,y_t) -w_t\|^2  \nonumber \\
 & \quad + \frac{6d_2L^2_f\eta_{t}\gamma}{b\tau^2}\big(\gamma^2\mathbb{E}\|v_t\|^2
 + \mathbb{E}\|\tilde{y}_{t+1}-y_t\|^2\big) + \frac{3L^2_f\mu_1^2d_1^2\gamma}{b\eta_t\tau^2} \!+\! \frac{2\alpha_{t+1}^2\delta^2\gamma}{b\eta_t\tau^2}+ \frac{3L^2_f\mu_2^2d_2^2\gamma}{b\eta_t\tau^2} \!+\! \frac{2\beta_{t+1}^2\delta^2\gamma}{b\eta_t\tau^2} \nonumber \\
 & \leq -\frac{13\gamma L_f^2\eta_t}{4}\mathbb{E}\|y_t - y^*(x_t)\|^2 - \frac{5\gamma\eta_t}{4}\mathbb{E} \|\nabla_xf_{\mu_1}(x_t,y_t) -v_t\|^2 - \frac{\gamma\eta_t}{2}\mathbb{E} \|\nabla F_{\mu_1}(x_t)\|^2 \nonumber \\
 & \quad + \frac{3\mu_1^2d_1^2L_f^2\eta_t\gamma}{2} \!+\! \frac{625\tilde{d}d^2_2L^4_f\mu^2_2\eta_t\gamma}{12b\tau^2} + \frac{3L^2_f\mu_1^2d_1^2\gamma}{b\eta_t\tau^2} + \frac{2\alpha_{t+1}^2\delta^2\gamma}{b\eta_t\tau^2} + \frac{3L^2_f\mu_2^2d_2^2\gamma}{b\eta_t\tau^2} \!+\! \frac{2\beta_{t+1}^2\delta^2\gamma}{b\eta_t\tau^2} \nonumber \\
 & \quad - \big(\frac{75\tilde{d}L^2_f\gamma}{4b\lambda\tau} - \frac{6\tilde{d}L^2_f\gamma}{b\tau^2}\big)\eta_t\mathbb{E}\|\tilde{y}_{t+1}-y_t\|^2 -\big(\frac{\gamma}{4} -\frac{6\tilde{d}L^2_f\gamma^3}{b\tau^2} - \frac{625\tilde{d}L^2_f\kappa_y^2\gamma^3}{6b\lambda^2\tau^2}\big)\eta_t\mathbb{E}\|v_t\|^2 \nonumber \\
 & \leq - \frac{\gamma\eta_t}{2}\mathbb{E} \|\nabla F_{\mu_1}(x_t)\|^2 + 3\mu_1^2d_1^2L_f^2\eta_t\gamma+ \frac{625\tilde{d}d^2_2L^4_f\mu^2_2\eta_t\gamma}{12b\tau^2}  \nonumber \\
 & \quad + \frac{3L^2_f\mu_1^2d_1^2\gamma}{b\eta_t\tau^2} \!+\! \frac{2\alpha_{t+1}^2\delta^2\gamma}{b\eta_t\tau^2} + \frac{3L^2_f\mu_2^2d_2^2\gamma}{b\eta_t\tau^2} \!+\! \frac{2\beta_{t+1}^2\delta^2\gamma}{b\eta_t\tau^2},
\end{align}
where the first inequality holds by combining the above inequalities \eqref{eq:G01}, \eqref{eq:G02} and \eqref{eq:G03},
  and the second inequality is due to $1\leq b\leq \tilde{d}$
 and the last inequality is due to $0< \gamma \leq \frac{\lambda\tau^2}{2L_f}\sqrt{\frac{6b/\tilde{d}}{36 \lambda^2 + 625\kappa_y^2}}$
 and $\lambda\leq \frac{75\tau}{24}$.
Thus, we have
 \begin{align} \label{eq:G04}
 \frac{\eta_t}{2}\mathbb{E} \|\nabla F_{\mu_1}(x_t)\|^2 & \leq \frac{\Phi_t - \Phi_{t+1}}{\gamma} + 3\mu_1^2d_1^2L_f^2\eta_t + \frac{625\tilde{d}d^2_2L^4_f\mu^2_2\eta_t}{12b\tau^2}\nonumber \\
 &\quad + \frac{3L^2_f\mu_1^2d_1^2}{b\eta_t\tau^2} \!+\! \frac{2\alpha_{t+1}^2\delta^2}{b\eta_t\tau^2} + \frac{3L^2_f\mu_2^2d_2^2}{b\eta_t\tau^2} \!+\! \frac{2\beta_{t+1}^2\delta^2}{b\eta_t\tau^2}.
 \end{align}
Since $\inf_{x\in \mathcal{X}} F(x) =F^*$, we have
$\inf_{x\in \mathcal{X}} F_{\mu_1}(x) =\inf_{x\in \mathcal{X}} \mathbb{E}_{u_1\sim U_B}[F(x+\mu_1 u_1)] =
\inf_{x\in \mathcal{X}} \frac{1}{V}\int_{B} \\ F(x+\mu_1 u_1)du_1 \geq \frac{1}{V}\int_{B}\inf_{x\in \mathcal{X}}F(x+\mu_1 u_1)du_1=F^*$,
where $V$ denotes the volume of the unit ball $B$.
Let $\Delta_1=\|y_1-y^*(x_1)\|^2$, we have
\begin{align} \label{eq:G05}
 \Phi_1 &= F_{\mu_1}(x_1)+ \frac{ 25\gamma\tilde{d}L^2_f}{\lambda\tau b}\|y_1-y^*(x_1)\|^2+\frac{\gamma}{\eta_0\tau^2}\mathbb{E}\|\nabla_x f_{\mu_1}(x_1,y_1)-v_1\|^2 \nonumber \\
 &\quad + \frac{\gamma}{\eta_{0}\tau^2}\mathbb{E}\|\nabla_y f_{\mu_2}(x_1,y_1)-w_1\|^2 \nonumber \\
 & = F_{\mu_1}(x_1)+ \frac{ 25\gamma\tilde{d}L^2_f}{\lambda\tau b}\|y_1-y^*(x_1)\|^2+\frac{\gamma}{\eta_0\tau^2}\mathbb{E}\|\nabla_x f_{\mu_1}(x_1,y_1)-\hat{\nabla}_x f(x_1,y_1;\mathcal{B}_1)\|^2 \nonumber \\
 &\quad + \frac{\gamma}{\eta_{0}\tau^2}\mathbb{E}\|\nabla_y f_{\mu_2}(x_1,y_1)-\hat{\nabla}_y f(x_1,y_1;\mathcal{B}_1)\|^2 \nonumber \\
 & \leq F_{\mu_1}(x_1)+ \frac{ 25\gamma\tilde{d}L^2_f}{\lambda\tau b}\Delta_1  +\frac{2\gamma\delta^2}{b\eta_0\tau^2},
\end{align}
where the last inequality holds by Assumption \ref{ass:2}.

Taking average over $t=1,2,\cdots,T$ on both sides of \eqref{eq:G04} and by using  $\eta_T^{-1} \geq \eta_t^{-1}$
for any $0\leq t\leq T$, we have
 \begin{align}
 & \frac{1}{T}\sum_{t=1}^T \frac{1}{2}\mathbb{E}\|\nabla F_{\mu_1}(x_t)\|^2 \nonumber \\
 & \leq  \frac{1}{T\gamma\eta_T} \sum_{t=1}^T\big(\Phi_t - \Phi_{t+1}\big)  +  \frac{1}{T\eta_T}\sum_{t=1}^T\big( 3\mu_1^2d_1^2L_f^2\eta_t + \frac{625\tilde{d}d^2_2L^4_f\mu^2_2\eta_t}{12b\tau^2} + \frac{3L^2_f\mu_1^2d_1^2}{b\eta_t\tau^2} \!+\! \frac{2\alpha_{t+1}^2\delta^2}{b\eta_t\tau^2} \nonumber \\
 & \quad  + \frac{3L^2_f\mu_2^2d_2^2}{b\eta_t\tau^2} + \frac{2\beta_{t+1}^2\delta^2}{b\eta_t\tau^2}\big) \nonumber \\
 & \leq \frac{F_{\mu_1}(x_1) - F^*}{T\gamma\eta_T} + \frac{ 25\tilde{d}L^2_f}{T\eta_T\lambda\tau b}\Delta_1 +\frac{2\delta^2}{Tb\tau^2\eta_T\eta_0} + \frac{36\tau^2\mu_1^2d_1^2L_f^2 + 625\tilde{d}d^2_2L^4_f\mu^2_2}{12b\tau^2T\eta_T}\sum_{t=1}^T\eta_t  \nonumber \\
 & \quad + \frac{3L^2_f\big(\mu_1^2 d_1^2 + \mu_2^2d_2^2\big)}{Tb\tau^2\eta_T}\sum_{t=1}^T\frac{1}{\eta_t}+ \frac{2(c_1^2+c_2^2)\delta^2}{Tb\tau^2\eta_T}\sum_{t=1}^T\eta_t^3
 \nonumber \\
 & \leq \frac{F_{\mu_1}(x_1) - F^*}{T\gamma\eta_T}+ \frac{ 25\tilde{d}L^2_f}{T\eta_T\lambda\tau b}\Delta_1 +\frac{2\delta^2}{Tb\tau^2\eta_T\eta_0} + \frac{36\tau^2\mu_1^2d_1^2L_f^2 + 625\tilde{d}d^2_2L^4_f\mu^2_2}{12b\tau^2T\eta_T}\int^T_1\frac{k}{(m+t)^{1/3}}dt  \nonumber \\
 &\quad + \frac{3L^2_f\big(\mu_1^2 d_1^2 + \mu_2^2d_2^2\big)}{Tb\tau^2\eta_T}\int^T_1\frac{(m+t)^{1/3}}{k}dt + \frac{2(c_1^2+c_2^2)\delta^2}{Tb\tau^2\eta_T}\int^T_1\frac{k^3}{m+t} dt\nonumber \\
 & \leq \frac{F_{\mu_1}(x_1) - F^*}{T\gamma\eta_T} + \frac{ 25\tilde{d}L^2_f}{T\eta_T\lambda\tau b}\Delta_1 + \frac{2\delta^2}{Tb\tau^2\eta_T\eta_0} + \frac{36\tau^2L_f^2k + 625L^4_fk}{8b\tau^2T\eta_T}(m+T)^{-2/3} + \frac{9L^2_f}{4Tb\tau^2\eta_T k} \nonumber \\
 & \quad + \frac{2(c_1^2+c_2^2)\delta^2 k^3}{Tb\tau^2\eta_T}\ln(m+T) \nonumber \\
 & = \big(\frac{F_{\mu_1}(x_1) - F^*}{T\gamma k}+ \frac{ 25\tilde{d}L^2_f}{Tk\lambda\tau b}\Delta_1 +\frac{2m^{1/3}\delta^2}{Tb\tau^2 k^2}\big)(m+T)^{1/3} + \frac{36\tau^2L_f^2 + 625L^4_f}{8b\tau^2T}(m+T)^{-1/3}  \nonumber \\
& \quad + \frac{9L^2_f}{4Tb\tau^2 k^2}(m+T)^{1/3} + \frac{2(c_1^2+c_2^2)\delta^2 k^2}{Tb\tau^2}\ln(m+T)(m+T)^{1/3},
\end{align}
where the second inequality holds by the above inequality \eqref{eq:G05}, and the last inequality is due to $0<\mu_1\leq \frac{1}{d_1(m+T)^{2/3}}$
and $0<\mu_2\leq \frac{1}{\tilde{d}^{1/2}d_2(m+T)^{2/3}}$. Let $M' = \frac{F_{\mu_1}(x_1) - F^*}{\gamma k} + \frac{ 25\tilde{d}L^2_f}{k\lambda\tau b}\Delta_1 +\frac{2m^{1/3}\delta^2}{b\tau^2 k^2} + \frac{36\tau^2L_f^2 + 625L^4_f}{8b\tau^2}(m+T)^{-2/3} + \frac{9L^2_f}{4b\tau^2 k^2} + \frac{2(c_1^2+c_2^2)\delta^2 k^2}{b\tau^2}\ln(m+T)$,
we have
\begin{align}
 \frac{1}{T} \sum_{t=1}^T \mathbb{E}\|\nabla F_{\mu_1}(x_t)\|^2  \leq \frac{2M'}{T}(m+T)^{1/3}.
\end{align}
According to Jensen's inequality, we have
\begin{align}
  \frac{1}{T} \sum_{t=1}^T \mathbb{E}\|\nabla F_{\mu_1}(x_t)\|
 & \leq \big( \frac{1}{T} \sum_{t=1}^T \mathbb{E} \|\nabla F_{\mu_1}(x_t)\|^2 \big)^{1/2} \nonumber \\
 & \leq \frac{\sqrt{2M'}}{T^{1/2}}(m+T)^{1/6} \leq \frac{\sqrt{2M'}m^{1/6}}{T^{1/2}} + \frac{\sqrt{2M'}}{T^{1/3}},
\end{align}
where the last inequality is due to $(a+b)^{1/6} \leq a^{1/6} + b^{1/6}$.
According to Lemma \ref{lem:A3}, we have $\|\nabla F_{\mu_1}(x_t)\|=\|\nabla F_{\mu_1}(x_t) - \nabla F(x_t)
+ \nabla F(x_t)\| \geq \|\nabla F(x_t)\| - \|\nabla F_{\mu_1}(x_t) - \nabla F(x_t)\| \geq \|\nabla F(x_t)\|
- \frac{\mu_1 L_f d_1}{2}$.
Thus, we have
\begin{align}
 \frac{1}{T} \sum_{t=1}^T \mathbb{E}\|\nabla F(x_t)\|
 & \leq \frac{1}{T} \sum_{t=1}^T  \big( \mathbb{E}\|\nabla F_{\mu_1}(x_t)\| + \frac{\mu_1 L_f d_1}{2} \big) \nonumber \\
 & \leq  \frac{\sqrt{2M'}m^{1/6}}{T^{1/2}} + \frac{\sqrt{2M'}}{T^{1/3}} + \frac{\mu_1 L_f d_1}{2} \nonumber \\
 & \leq  \frac{\sqrt{2M'}m^{1/6}}{T^{1/2}} + \frac{\sqrt{2M'}}{T^{1/3}} + \frac{L_f}{2(m+T)^{2/3}},
\end{align}
where the last inequality is due to $0<\mu_1 \leq \frac{1}{d_1(m+T)^{2/3}}$.

\end{proof}

\subsection{ Convergence Analysis of Acc-MDA Algorithm for Constrained Minimax Optimization}
\label{Appendix:A5}
In this subsection, we study the convergence properties of our  Acc-MDA algorithm for solving the  \textbf{constrained}
minimax problem \eqref{eq:2}, i.e., $\mathcal{X}\subset \mathbb{R}^{d_1}$
and $\mathcal{Y} \subset \mathbb{R}^{d_2}$ (or $\mathcal{Y} = \mathbb{R}^{d_2}$),
where the noise stochastic gradients of function $f(x,y)$ can be obtained. The following convergence analysis builds on a new metric $\mathbb{E}[\mathcal{H}_t]$, where $\mathcal{H}_t$ is defined in \eqref{eq:14}.

\begin{lemma} \label{lem:D3}
 Suppose the sequence $\{x_t,y_t\}_{t=1}^T$ be generated from Algorithm \ref{alg:4}.
 Let $0<\eta_t\leq 1$ and $0< \gamma \leq \frac{1}{2L_g\eta_t}$,
 we have
 \begin{align}
  F(x_{t+1}) - F(x_t) & \leq -\frac{\eta_t}{2\gamma}\|\tilde{x}_{t+1}-x_t\|^2 + 2\eta_t\gamma L_f^2\|y^*(x_t)-y_t\|^2 + 2\eta_t\gamma\|\nabla_xf(x_t,y_t) -v_t\|^2,
 \end{align}
 where $L_g = L_f + L_f^2/\tau$.
 \end{lemma}
 \begin{proof}
 This proof is similar to the proof of Lemma \ref{lem:D1}.
 According to Lemma \ref{lem:1}, the function $F(x)$ has $L_g$-Lipschitz continuous gradient.
 Then we have
 \begin{align} \label{eq:J1}
  F(x_{t+1}) &\leq F(x_t) + \langle\nabla F(x_t), x_{t+1}-x_t\rangle + \frac{L_g}{2}\|x_{t+1}-x_t\|^2   \\
  & = F(x_t) + \eta_t\langle \nabla F(x_t),\tilde{x}_{t+1}-x_t\rangle + \frac{L_g\eta_t^2}{2}\|\tilde{x}_{t+1}-x_t\|^2 \nonumber \\
  & = F(x_t) + \eta_t\langle \nabla F(x_t)-v_t,\tilde{x}_{t+1}-x_t\rangle + \eta_t\langle v_t,\tilde{x}_{t+1}-x_t\rangle + \frac{L_g\eta_t^2}{2}\|\tilde{x}_{t+1}-x_t\|^2. \nonumber
 \end{align}
 By the step 8 of Algorithm \ref{alg:4}, we have $\tilde{x}_{t+1}=\mathcal{P}_{\mathcal{X}}(x_t - \gamma v_t)
 =\arg\min_{x\in\mathcal{X}}\frac{1}{2}\|x- x_t + \gamma v_t\|^2$. Since $\mathcal{X}$ is a convex set
 and the function $\frac{1}{2}\|x- x_t + \gamma v_t\|^2$ is convex, according to Lemma \ref{lem:A1},
 we have
 \begin{align} \label{eq:J2}
  \langle \tilde{x}_{t+1}- x_t + \gamma v_t, x-\tilde{x}_{t+1}\rangle \geq 0, \ \forall x\in \mathcal{X}.
 \end{align}
 In Algorithm \ref{alg:4}, let the initialize solution $x_1 \in \mathcal{X}$, and the sequence $\{x_t\}_{t\geq 1}$ generates as follows:
 \begin{align}
  x_{t+1} = x_{t} + \eta_t(\tilde{x}_{t+1} - x_t) = \eta_t\tilde{x}_{t+1} + (1-\eta_t)x_t,
 \end{align}
 where $0<\eta_t \leq 1$. Since $\mathcal{X}$ is convex set and $x_t, \tilde{x}_{t+1} \in \mathcal{X}$, we have $x_{t+1} \in \mathcal{X}$ for any $t >0$.
 Set $x=x_t$ in the inequality \eqref{eq:J2}, we have
 \begin{align} \label{eq:J3}
 \langle v_t, \tilde{x}_{t+1}-x_t\rangle \leq -\frac{1}{\gamma}\|\tilde{x}_{t+1}-x_t\|^2.
 \end{align}

 Next, we decompose the term $\langle \nabla F(x_t)-v_t,\tilde{x}_{t+1}-x_t\rangle$ as follows:
 \begin{align}
 &\langle \nabla F(x_t)-v_t,\tilde{x}_{t+1}-x_t\rangle \nonumber \\
 & = \underbrace{\langle \nabla F(x_t) - \nabla_xf(x_t,y_t),\tilde{x}_{t+1}-x_t\rangle}_{=T_1}
 + \underbrace{\langle  \nabla_xf(x_t,y_t) -v_t,\tilde{x}_{t+1}-x_t\rangle}_{=T_2}.
 \end{align}
 For the term $T_1$, by the Cauchy-Schwarz inequality and Young's inequality, we have
 \begin{align}
  T_1 &= \langle \nabla F(x_t) - \nabla_xf(x_t,y_t),\tilde{x}_{t+1}-x_t\rangle \nonumber \\
  & \leq \|\nabla F(x_t) - \nabla_xf(x_t,y_t)\|\cdot\|\tilde{x}_{t+1}-x_t\| \nonumber \\
  & \leq 2\gamma\|\nabla F(x_t) - \nabla_xf(x_t,y_t)\|^2 + \frac{1}{8\gamma}\|\tilde{x}_{t+1}-x_t\|^2 \nonumber \\
  & =2\gamma\|\nabla_xf(x_t,y^*(x_t)) - \nabla_xf(x_t,y_t)\|^2 + \frac{1}{8\gamma}\|\tilde{x}_{t+1}-x_t\|^2 \nonumber \\
  & \leq 2\gamma\|\nabla f(x_t,y^*(x_t)) - \nabla f(x_t,y_t)\|^2 + \frac{1}{8\gamma}\|\tilde{x}_{t+1}-x_t\|^2\nonumber \\
  & \leq 2\gamma L^2_f\|y^*(x_t)-y_t\|^2 + \frac{1}{8\gamma}\|\tilde{x}_{t+1}-x_t\|^2,
 \end{align}
 where the last inequality holds by Assumption \ref{ass:5}.

 For the term $T_2$, by the Cauchy-Schwarz inequality and Young's inequality, we have
 \begin{align}
  T_2 & = \langle  \nabla_xf(x_t,y_t) -v_t,\tilde{x}_{t+1}-x_t \rangle \nonumber \\
  & \leq \|\nabla_xf(x_t,y_t) -v_t\| \cdot \|\tilde{x}_{t+1}-x_t\| \nonumber \\
  & \leq 2\gamma\|\nabla_xf(x_t,y_t) -v_t\|^2 + \frac{1}{8\gamma}\|\tilde{x}_{t+1}-x_t\|^2,
 \end{align}
 where the last inequality holds by $\langle a,b\rangle \leq \frac{\lambda}{2}\|a\|^2 + \frac{1}{2\lambda}\|b\|^2$
 with $\lambda=4\gamma$.
 Thus, we have
 \begin{align} \label{eq:J4}
 \langle \nabla F(x_t)-v_t,\tilde{x}_{t+1}-x_t\rangle
 & = 2\gamma L_f^2\|y^*(x_t)-y_t\|^2 + 2\gamma\|\nabla_xf(x_t,y_t) -v_t\|^2 + \frac{1}{4\gamma}\|\tilde{x}_{t+1}-x_t\|^2.
 \end{align}
 Finally, combining the inequalities \eqref{eq:J1}, \eqref{eq:J3} with \eqref{eq:J4}, we have
 \begin{align}
 F(x_{t+1}) &\leq F(x_t) + 2\eta_t\gamma L_f^2\|y^*(x_t)-y_t\|^2 + 2\eta_t\gamma\|\nabla_xf(x_t,y_t) -v_t\|^2 + \frac{\eta_t}{4\gamma}\|\tilde{x}_{t+1}-x_t\|^2 \nonumber \\
 & \quad-\frac{\eta_t}{\gamma}\|\tilde{x}_{t+1}-x_t\|^2 + \frac{L_g\eta_t^2}{2}\|\tilde{x}_{t+1}-x_t\|^2\nonumber \\
 & \leq F(x_t) + 2\eta_t\gamma L_f^2\|y^*(x_t)-y_t\|^2
 + 2\eta_t\gamma\|\nabla_xf(x_t,y_t) -v_t\|^2 -\frac{\eta_t}{2\gamma}\|\tilde{x}_{t+1}-x_t\|^2,
 \end{align}
 where the last inequality is due to $0< \gamma \leq \frac{1}{2L_g\eta_t}$.

\end{proof}

\begin{lemma} \label{lem:E3}
Suppose the sequence $\{x_t,y_t\}_{t=1}^T$ be generated from Algorithm \ref{alg:4}. Under the above assumptions, and set $0< \eta_t\leq 1$
and $\lambda\leq \frac{1}{6L_f}$, we have
\begin{align}
     \|y_{t+1} - y^*(x_{t+1})\|^2 &\leq (1-\frac{\eta_t\tau\lambda}{4})\|y_t -y^*(x_t)\|^2 -\frac{3\eta_t}{4} \|\tilde{y}_{t+1}-y_t\|^2 \nonumber \\
     & \quad + \frac{25\eta_t\lambda}{6\tau}  \|\nabla_y f(x_t,y_t)-w_t\|^2 +  \frac{25\kappa_y^2\eta_t}{6\tau\lambda}\|x_t - \tilde{x}_{t+1}\|^2,
\end{align}
where $\kappa_y = L_f/\tau$.
\end{lemma}
\begin{proof}
 This proof is the same to the proof of Lemma \ref{lem:E1}.
\end{proof}

\begin{lemma} \label{lem:F3}
 Suppose the stochastic gradients $\{v_t,w_t\}_{t=1}^T$ be generated from Algorithm \ref{alg:4}, we have
\begin{align}
\mathbb{E} \|\nabla_x f(x_{t+1},y_{t+1}) - v_{t+1}\|^2 & \leq (1-\alpha_{t+1})^2 \mathbb{E} \|\nabla_x f(x_t,y_t) -v_t\|^2
 + \frac{2\alpha_{t+1}^2\delta^2}{b} \nonumber \\
 & + \frac{2(1-\alpha_{t+1})^2L^2_f\eta^2_t}{b}\big(\mathbb{E}\|\tilde{x}_{t+1}-x_t\|^2 + \mathbb{E}\|\tilde{y}_{t+1}-y_t\|^2\big).
\end{align}
\begin{align}
\mathbb{E} \|\nabla_y f(x_{t+1},y_{t+1}) - w_{t+1}\|^2 & \leq (1-\beta_{t+1})^2 \mathbb{E} \|\nabla_y f(x_t,y_t) -w_t\|^2 + \frac{2\beta_{t+1}^2\delta^2}{b}
 \nonumber \\
 & + \frac{2(1-\beta_{t+1})^2L^2_f\eta^2_t}{b}\big(\mathbb{E} \|\tilde{x}_{t+1}-x_t\|^2 + \mathbb{E} \|\tilde{y}_{t+1}-y_t\|^2\big).
\end{align}
\end{lemma}
\begin{proof}
 This proof is the same to the proof of Lemma \ref{lem:F1}.
 According to the definition of $w_{t+1}$ in Algorithm \ref{alg:4}, we have
 \begin{align}
  w_{t+1}-w_t & = -\beta_{t+1}w_t + (1-\beta_{t+1})\big(\nabla_y f(x_{t+1},y_{t+1};\mathcal{B}_{t+1}) - \nabla_y f(x_t,y_t;\mathcal{B}_{t+1})\big) \nonumber \\
  & \quad + \beta_{t+1}\nabla_y f(x_{t+1},y_{t+1};\mathcal{B}_{t+1}). \nonumber
 \end{align}
 Then we have
 \begin{align}
  & \mathbb{E} \|\nabla_y f(x_{t+1},y_{t+1}) - v_{t+1}\|^2 \nonumber \\
  &\! = \mathbb{E} \|\nabla_y f(x_{t+1},y_{t+1}) - v_t - (v_{t+1}-v_t)\|^2 \\
  & \!= \mathbb{E} \|\nabla_y f(x_{t+1},y_{t+1}) - v_t + \beta_{t+1}v_t- \beta_{t+1}\nabla_y f(x_{t+1},y_{t+1};\mathcal{B}_{t+1})
   \nonumber \\
  & \ - (1-\beta_{t+1})(\nabla_y f(x_{t+1},y_{t+1};\mathcal{B}_{t+1})- \nabla_y f(x_t,y_t;\mathcal{B}_{t+1})) \|^2 \nonumber \\
  & \!= \mathbb{E} \|(1-\beta_{t+1})(\nabla_y f(x_t,y_t) -v_t) + \beta_{t+1}\big(\nabla_y f(x_{t+1},y_{t+1})- \nabla_y f(x_{t+1},y_{t+1};\mathcal{B}_{t+1})\big)\nonumber \\
  & \ + (1-\beta_{t+1})\big(\nabla_y f(x_{t+1},y_{t+1})
  -\nabla_y f(x_t,y_t)-\nabla_y f(x_{t+1},y_{t+1};\mathcal{B}_{t+1}) + \nabla_y f(x_t,y_t;\mathcal{B}_{t+1})\big)\|^2 \nonumber \\
  &\! = (1-\beta_{t+1})^2\mathbb{E}\|\nabla_y f(x_t,y_t)-v_t\|^2 +  \mathbb{E} \|\beta_{t+1}\big(\nabla_y f(x_{t+1},y_{t+1})- \nabla_y f(x_{t+1},y_{t+1};\mathcal{B}_{t+1})\big)\nonumber \\
  & \ + (1-\beta_{t+1})\big(\nabla_y f(x_{t+1},y_{t+1})
  -\nabla_y f(x_t,y_t)-\nabla_y f(x_{t+1},y_{t+1};\mathcal{B}_{t+1}) + \nabla_y f(x_t,y_t;\mathcal{B}_{t+1})\big)\|^2 \nonumber \\
  & \!\leq (1-\beta_{t+1})^2\mathbb{E} \|\nabla_y f(x_t,y_t) -v_t\|^2 + 2\beta_{t+1}^2 \mathbb{E} \|\nabla_y f(x_{t+1},y_{t+1})- \nabla_y f(x_{t+1},y_{t+1};\mathcal{B}_{t+1})\|^2  \nonumber \\
  & \ + 2(1-\beta_{t+1})^2\mathbb{E} \|\nabla_y f(x_{t+1},y_{t+1}) \!-\! \nabla_y f(x_t,y_t)\!-\! \nabla_y f(x_{t+1},y_{t+1};\mathcal{B}_{t+1}) \!+\! \nabla_y f(x_t,y_t;\mathcal{B}_{t+1})\|^2 \nonumber \\
  & \!\leq (1-\beta_{t+1})^2 \mathbb{E} \|\nabla_y f(x_t,y_t) -v_t\|^2 + \frac{2(1-\beta_{t+1})^2}{b} \mathbb{E} \| \nabla_y f(x_{t+1},y_{t+1};\mathcal{B}_{t+1})  \nonumber \\
  & \ -
  \nabla_y f(x_t,y_t;\mathcal{B}_{t+1})\|^2+ \frac{2\beta_{t+1}^2\delta^2}{b} \nonumber \\
  &\! \leq (1-\beta_{t+1})^2 \mathbb{E} \|\nabla_y f(x_t,y_t) -v_t\|^2 + \frac{2(1-\beta_{t+1})^2L^2_f \eta_t^2}{b}\big( \mathbb{E}\|\tilde{x}_{t+1} -
  x_t\|^2 + \mathbb{E}\|\tilde{y}_{t+1} -y_t\|^2 \big) \nonumber\\
  & \ + \frac{2\beta_{t+1}^2\delta^2}{b},  \nonumber
 \end{align}
 where the fourth equality follows by $\mathbb{E}_{\mathcal{B}_{t+1}}[\nabla_y f(x_{t+1},y_{t+1};\mathcal{B}_{t+1})]=\nabla_y f(x_{t+1},y_{t+1})$ and $ \mathbb{E}_{\mathcal{B}_{t+1}}[\nabla_y f(x_{t+1},y_{t+1};\mathcal{B}_{t+1}) - \nabla_y f(x_t,y_t;\mathcal{B}_{t+1})]=\nabla_y f(x_{t+1},y_{t+1}) - \nabla_y f(x_t,y_t)$; the second inequality is due to Lemma \ref{lem:A4} and Assumption \ref{ass:2}; the last inequality holds by Assumption \ref{ass:5}.
 Similarly, we can obtain
 \begin{align}
  \mathbb{E} \|\nabla_x f(x_t,y_t) - v_t\|^2 & \leq (1-\alpha_t)^2 \mathbb{E} \|\nabla_x f(x_{t-1},y_{t-1}) -v_{t-1}\|^2
 + \frac{2\alpha_t^2\delta^2}{b} \nonumber \\
 & \quad + \frac{2(1-\alpha_t)^2L^2_f\eta^2_{t-1}}{b}\big(\mathbb{E}\|\tilde{x}_t-x_{t-1}\|^2 + \mathbb{E}\|\tilde{y}_t-y_{t-1}\|^2\big).
 \end{align}

\end{proof}

\begin{theorem} \label{th:A4}
 (Restatement of Theorem 9)
Suppose the sequence $\{x_t,y_t\}_{t=1}^T$ be generated from Algorithm \ref{alg:4}. When $\mathcal{X}\subset \mathbb{R}^{d_1}$, and let $\eta_t = \frac{k}{(m+t)^{1/3}}$
for all $t\geq 0$, $c_1 \geq \frac{2}{3k^3} + \frac{9\tau^2}{4}$
and $c_2 \geq \frac{2}{3k^3} + \frac{75L^2_f}{2}$, $k>0$, $m\geq \max\big(2, k^3, (c_1k)^3, (c_2k)^3\big)$,  $0<\lambda\leq \min\big(\frac{1}{6L_f},
\frac{27b\tau}{16}\big)$ and
$0< \gamma \leq \min\big( \frac{\lambda\tau}{2L_f}\sqrt{\frac{2b}{8\lambda^2 + 75\kappa_y^2b}}, \frac{m^{1/3}}{2L_gk}\big)$,
we have
\begin{align}
 \frac{1}{T} \sum_{t=1}^T \mathbb{E}  \|G_{\mathcal{X}}(x_t,\nabla F(x_t),\gamma)\| & \!\leq \! \frac{1}{T} \sum_{t=1}^T \mathbb{E}\big[L_f\|y^*(x_t)-y_t\| \!+\! \|\nabla_xf(x_t,y_t) -v_t\| \!+\! \frac{1}{\gamma}\|\tilde{x}_{t+1}-x_t\| \big] \nonumber \\
 & \!\leq\! \frac{2\sqrt{3M''}m^{1/6}}{T^{1/2}} + \frac{2\sqrt{3M''}}{T^{1/3}},
\end{align}
where $\Delta_1=\|y_1-y^*(x_1)\|^2$ and  $M'' =  \frac{F(x_1) - F^*}{\gamma k} +  \frac{9L^2_f\Delta_1}{ k\lambda\tau}+\frac{2m^{1/3}\delta^2}{b \tau^2k^2} + \frac{2(c_1^2+c_2^2)\delta^2 k^2}{b\tau^2}\ln(m+T)$.
\end{theorem}
\begin{proof}
 Since $\eta_t=\frac{k}{(m+t)^{1/3}}$ on $t$ is decreasing and $m\geq k^3$, we have $\eta_t \leq \eta_0 = \frac{k}{m^{1/3}} \leq 1$ and $\gamma \leq \frac{m^{1/3}}{2L_gk}=\frac{1}{2L_g\eta_0} \leq \frac{1}{2L_g\eta_t}$ for any $t\geq 0$.
 Due to $0 < \eta_t \leq 1$ and $m\geq \max\big( (c_1k)^3, (c_2k)^3 \big)$, we have $\alpha_{t+1} = c_1\eta_t^2 \leq c_1\eta_t \leq \frac{c_1k}{m^{1/3}}\leq 1$ and $\beta_{t+1} = c_2\eta_t^2 \leq c_2\eta_t \leq \frac{c_2k}{m^{1/3}}\leq 1$.
 According to Lemma \ref{lem:F3}, we have
 \begin{align}
  & \frac{1}{\eta_t}\mathbb{E} \|\nabla_x f(x_{t+1},y_{t+1}) - v_{t+1}\|^2 - \frac{1}{\eta_{t-1}}\mathbb{E} \|\nabla_x f(x_t,y_t) - v_t\|^2  \\
  & \leq \big(\frac{(1\!-\!\alpha_{t+1})^2}{\eta_t} \!-\! \frac{1}{\eta_{t-1}}\big)\mathbb{E} \|\nabla_x f(x_t,y_t) -v_t\|^2 +\! \frac{2L^2_f(1\!-\!\alpha_{t+1})^2\eta_{t}}{b}\mathbb{E}\big(\|\tilde{x}_{t+1}-x_t\|^2 \!+\! \|\tilde{y}_{t+1}-y_t\|^2\big) \nonumber \\
  & \quad  + \frac{2\alpha_{t+1}^2\delta^2}{b\eta_t}\nonumber \\
  & \leq \big(\frac{1\!-\!\alpha_{t+1}}{\eta_t} \!-\! \frac{1}{\eta_{t-1}}\big)\mathbb{E} \|\nabla_x f(x_t,y_t) -v_t\|^2 \!+\! \frac{2 L^2_f\eta_{t}}{b}\mathbb{E}\big(\|\tilde{x}_{t+1}-x_t\|^2 \!+\! \|\tilde{y}_{t+1}-y_t\|^2\big) \!+\! \frac{2\alpha_{t+1}^2\delta^2}{b\eta_t}\nonumber \\
  & = \big(\frac{1}{\eta_t} \!-\! \frac{1}{\eta_{t-1}} - c_1\eta_t\big)\mathbb{E} \|\nabla_x f(x_t,y_t) -v_t\|^2 \!+\! \frac{2L^2_f\eta_{t}}{b}\mathbb{E}\big(\|\tilde{x}_{t+1}-x_t\|^2 \!+\! \|\tilde{y}_{t+1}-y_t\|^2\big) \!+\! \frac{2\alpha_{t+1}^2\delta^2}{b\eta_t}, \nonumber
 \end{align}
 where the second inequality is due to $0<\alpha_{t+1}\leq 1$.
 Similarly, according to Lemma \ref{lem:F3}, we can obtain
 \begin{align}
  & \frac{1}{\eta_t}\mathbb{E} \|\nabla_y f(x_{t+1},y_{t+1}) - w_{t+1}\|^2
  -  \frac{1}{\eta_{t-1}}\mathbb{E} \|\nabla_y f(x_t,y_t) - w_t\|^2  \\
  & \leq \big(\frac{1}{\eta_t} \!-\! \frac{1}{\eta_{t-1}} - c_2\eta_t\big)\mathbb{E} \|\nabla_y f(x_t,y_t) -w_t\|^2 \!+\! \frac{2L^2_f\eta_{t}}{b}\mathbb{E}\big(\|\tilde{x}_{t+1}-x_t\|^2 \!+\! \|\tilde{y}_{t+1}-y_t\|^2\big)  \!+\! \frac{2\beta_{t+1}^2\delta^2}{b\eta_t}. \nonumber
 \end{align}
 By $\eta_t = \frac{k}{(m+t)^{1/3}}$, we have
 \begin{align}
  \frac{1}{\eta_t} - \frac{1}{\eta_{t-1}} &= \frac{1}{k}\big( (m+t)^{\frac{1}{3}} - (m+t-1)^{\frac{1}{3}}\big) \nonumber \\
  & \leq \frac{1}{3k(m+t-1)^{2/3}} \leq \frac{1}{3k\big(m/2+t\big)^{2/3}} \nonumber \\
  & \leq \frac{2^{2/3}}{3k(m+t)^{2/3}} = \frac{2^{2/3}}{3k^3}\frac{k^2}{(m+t)^{2/3}}= \frac{2^{2/3}}{3k^3}\eta_t^2 \leq \frac{2}{3k^3}\eta_t,
 \end{align}
 where the first inequality holds by the concavity of function $f(x)=x^{1/3}$, \emph{i.e.}, $(x+y)^{1/3}\leq x^{1/3} + \frac{y}{3x^{2/3}}$; the second inequality is due to $m\geq 2$,  and
 the last inequality is due to $0<\eta_t\leq 1$.
 Let $c_1 \geq \frac{2}{3k^3} + \frac{9\tau^2}{4}$, we have
 \begin{align} \label{eq:K1}
  & \frac{1}{\eta_t}\mathbb{E} \|\nabla_x f(x_{t+1},y_{t+1}) - v_{t+1}\|^2 - \frac{1}{\eta_{t-1}}\mathbb{E} \|\nabla_x f(x_t,y_t) - v_t\|^2  \\
  & \leq -\frac{9\tau^2\eta_t}{4}\mathbb{E} \|\nabla_x f(x_t,y_t) -v_t\|^2 + \frac{2L^2_f\eta_{t}}{b}\mathbb{E}\big(\|\tilde{x}_{t+1}-x_t\|^2 + \|\tilde{y}_{t+1}-y_t\|^2\big) + \frac{2\alpha_{t+1}^2\delta^2}{b\eta_t}. \nonumber
 \end{align}
 Let $c_2 \geq \frac{2}{3k^3} + \frac{75L^2_f}{2}$, we have
  \begin{align} \label{eq:K2}
  & \frac{1}{\eta_t}\mathbb{E} \|\nabla_y f(x_{t+1},y_{t+1}) - w_{t+1}\|^2
  -  \frac{1}{\eta_{t-1}}\mathbb{E} \|\nabla_y f(x_t,y_t) - w_t\|^2  \\
  & \leq - \frac{75L^2_f\eta_t}{2\tau^2}\mathbb{E} \|\nabla_y f(x_t,y_t) -w_t\|^2 \!+\! \frac{2L^2_f\eta_{t}}{b}\mathbb{E}\big(\|\tilde{x}_{t+1}-x_t\|^2 + \|\tilde{y}_{t+1}-y_t\|^2\big) + \frac{2\beta_{t+1}^2\delta^2}{b\eta_t}. \nonumber
 \end{align}

 Next, we define a \emph{Lyapunov} function, for any $t\geq 1$
 \begin{align}
 \Omega_t & = \mathbb{E}\big[ F(x_t) + \frac{9L^2_f\gamma}{\lambda\tau}\|y_t-y^*(x_t)\|^2 + \frac{\gamma}{\tau^2\eta_{t-1}}\|\nabla_x f(x_t,y_t)-v_t\|^2 \nonumber \\
 & \quad + \frac{\gamma}{\tau^2\eta_{t-1}}\|\nabla_y f(x_t,y_t)-w_t\|^2 \big]. \nonumber
 \end{align}
 Then we have
 \begin{align}
 & \Omega_{t+1} - \Omega_t \nonumber \\
 & = \mathbb{E}\big[F(x_{t+1}) - F(x_t)\big] + \frac{ 9L^2_f\gamma}{\lambda\tau}
 \big( \mathbb{E}\|y_{t+1}-y^*(x_{t+1})\|^2 - \mathbb{E}\|y_t-y^*(x_t)\|^2 \big)
  \nonumber \\
 & \quad + \frac{\gamma}{\tau^2} \big(\frac{1}{\eta_t}\mathbb{E}\|\nabla_x f(x_{t+1},y_{t+1})-v_{t+1}\|^2- \frac{1}{\eta_{t-1}}\mathbb{E}\|\nabla_x f(x_t,y_t)-v_t\|^2  \nonumber \\
 &\quad 
 + \frac{1}{\eta_t}\mathbb{E}\|\nabla_y f(x_{t+1},y_{t+1})-w_{t+1}\|^2- \frac{1}{\eta_{t-1}}\mathbb{E}\|\nabla_y f(x_t,y_t)-w_t\|^2 \big) \nonumber \\
 & \leq -\frac{\eta_t}{2\gamma}\mathbb{E}\|\tilde{x}_{t+1}-x_t\|^2 + 2\eta_t\gamma L_f^2\mathbb{E}\|y^*(x_t)-y_t\|^2 + 2\eta_t\gamma\mathbb{E} \|\nabla_x f(x_t,y_t) -v_t\|^2 \nonumber \\
 & \quad + \frac{9L^2_f\gamma}{\lambda\tau} \big( \!-\!\frac{\eta_t\tau\lambda}{4}\mathbb{E}\|y_t - y^*(x_t)\|^2 -\frac{3\eta_t}{4} \mathbb{E}\|\tilde{y}_{t+1}-y_t\|^2 + \frac{25\eta_t\lambda}{6\tau} \mathbb{E}\|\nabla_y f(x_t,y_t)-w_t\|^2  \nonumber \\
 & \quad + \frac{25\kappa_y^2\eta_t}{6\tau\lambda}\mathbb{E}\|x_t - \tilde{x}_{t+1}\|^2 \big) -\frac{9\gamma\eta_t}{4} \mathbb{E} \|\nabla_x f(x_t,y_t) -v_t\|^2 - \frac{2L^2_f\eta_{t}\gamma}{b\tau^2}\big(\mathbb{E}\|\tilde{x}_{t+1}-x_t\|^2  \nonumber \\
 & \quad + \mathbb{E}\|\tilde{y}_{t+1}-y_t\|^2\big)+ \frac{2\alpha_{t+1}^2\delta^2\gamma}{b\tau^2\eta_t}-\frac{75L^2_f\gamma}{2\tau^2} \eta_t \mathbb{E}\|\nabla_y f(x_t,y_t) -w_t\|^2  \nonumber \\
 &\quad + \frac{2L^2_f\eta_{t}\gamma}{b\tau^2}\big(\mathbb{E}\|\tilde{x}_{t+1}-x_t\|^2 + \mathbb{E}\|\tilde{y}_{t+1}-y_t\|^2\big) + \frac{2\beta_{t+1}^2\delta^2\gamma}{b\tau^2\eta_t} \nonumber \\
 & \leq -\frac{L_f^2\eta_t\gamma }{4}\mathbb{E}\|y^*(x_t)-y_t\|^2 - \frac{\gamma\eta_t}{4}\mathbb{E}\|\nabla_xf(x_t,y_t) -v_t\|^2 + \frac{2\alpha_{t+1}^2\delta^2\gamma}{b\tau^2\eta_t} + \frac{2\beta_{t+1}^2\delta^2\gamma}{b\tau^2\eta_t}\nonumber \\
 & \quad - \big( \frac{27L^2_f\gamma}{4\lambda\tau} - \frac{4L^2_f\gamma}{b\tau^2}\big)\eta_t \mathbb{E}\|\tilde{y}_{t+1}-y_t\|^2-\big(\frac{1}{2\gamma} -\frac{4L^2_f\gamma}{b\tau^2} -\frac{75L^2_f\kappa_y^2\gamma}{2\lambda^2\tau^2}\big)\eta_t\mathbb{E}\|\tilde{x}_{t+1}-x_t\|^2 \nonumber \\
 & \leq -\frac{L_f^2\eta_t\gamma}{4}\mathbb{E}\|y^*(x_t)-y_t\|^2 -\frac{\gamma\eta_t}{4}\mathbb{E}\|\nabla_xf(x_t,y_t) -v_t\|^2 -\frac{\eta_t}{4\gamma}\mathbb{E}\|\tilde{x}_{t+1}-x_t\|^2 \nonumber \\
 & \quad + \frac{2\alpha_{t+1}^2\delta^2\gamma}{b\tau^2\eta_t} + \frac{2\beta_{t+1}^2\delta^2\gamma}{b\tau^2\eta_t},
 \end{align}
 where the first inequality holds by Lemmas \ref{lem:D3}, \ref{lem:E3} and the above inequalities \eqref{eq:K1}, \eqref{eq:K2}, and
 the last inequality is due to $0< \gamma \leq \frac{\lambda\tau}{2L_f}\sqrt{\frac{2b}{8\lambda^2 + 75\kappa_y^2b}}$ and $\lambda\leq \frac{27b\tau}{16}$.
 Then we have
 \begin{align} \label{eq:K4}
 & \frac{L_f^2\eta_t}{4}\mathbb{E}\|y^*(x_t)-y_t\|^2 + \frac{\eta_t}{4}\mathbb{E}\|\nabla_xf(x_t,y_t) -v_t\|^2 +\frac{\eta_t}{4\gamma^2}\mathbb{E}\|\tilde{x}_{t+1}-x_t\|^2 \nonumber \\
 & \leq \frac{\Omega_t - \Omega_{t+1}}{\gamma}+ \frac{2\alpha_{t+1}^2\delta^2}{b\tau^2\eta_t} + \frac{2\beta_{t+1}^2\delta^2}{b\tau^2\eta_t}.
 \end{align}

Taking average over $t=1,2,\cdots,T$ on both sides of \eqref{eq:K4}, we have
\begin{align}
 & \frac{1}{T} \sum_{t=1}^T \big( \frac{L_f^2\eta_t}{4}\mathbb{E}\|y^*(x_t)-y_t\|^2 + \frac{\eta_t}{4}\mathbb{E}\|\nabla_xf(x_t,y_t) -v_t\|^2 +\frac{\eta_t}{4\gamma^2}\mathbb{E}\|\tilde{x}_{t+1}-x_t\|^2 \big) \nonumber \\
 & \leq  \sum_{t=1}^T \frac{\Omega_t - \Omega_{t+1}}{T\gamma} + \frac{1}{T}\sum_{t=1}^T\big( \frac{2\alpha_{t+1}^2\delta^2}{b\tau^2\eta_t} + \frac{2\beta_{t+1}^2\delta^2}{b\tau^2\eta_t}\big). \nonumber
\end{align}
Let $\Delta_1=\|y_1-y^*(x_1)\|^2$, we have
\begin{align} \label{eq:K5}
 \Omega_1 &= F(x_1)\!+\! \frac{9L^2_f\gamma}{\lambda\tau}\|y_1-y^*(x_1)\|^2\!+\!\frac{\gamma}{\tau^2\eta_0}\mathbb{E}\|\nabla_x f(x_1,y_1)-v_1\|^2 \!+\! \frac{\gamma}{\tau^2\eta_{0}}\mathbb{E}\|\nabla_y f(x_1,y_1)-w_1\|^2 \nonumber \\
 & = F(x_1) + \frac{9L^2_f\gamma}{\lambda\tau}\|y_1-y^*(x_1)\|^2+ \frac{\gamma}{\tau^2\eta_0}\mathbb{E}\|\nabla_x f(x_1,y_1)- \frac{1}{b}\sum_{i=1}^b\hat{\nabla}_x f(x_1,y_1;\xi^1_i)\|^2 \nonumber \\
 & \quad + \frac{\gamma}{\tau^2\eta_0}\mathbb{E}\|\nabla_y f(x_1,y_1)- \frac{1}{b}\sum_{i=1}^b\hat{\nabla}_y f(x_1,y_1;\xi^1_i)\|^2 \nonumber \\
 & \leq F(x_1)+ \frac{9L^2_f\gamma}{\lambda\tau}\Delta_1 +\frac{2\gamma\delta^2}{b\tau^2\eta_0},
\end{align}
where the last inequality holds by Assumption \ref{ass:2}.
Since $\eta_t$ is decreasing, i.e., $\eta_T^{-1} \geq \eta_t^{-1}$ for any $0\leq t\leq T$, we have
 \begin{align}
 & \frac{1}{T} \sum_{t=1}^T \mathbb{E}\big( \frac{L_f^2}{4}\|y^*(x_t)-y_t\|^2 + \frac{1}{4}\|\nabla_xf(x_t,y_t) -v_t\|^2 +\frac{1}{4\gamma^2}\|\tilde{x}_{t+1}-x_t\|^2 \big) \nonumber \\
 & \leq  \frac{1}{T\gamma\eta_T} \sum_{t=1}^T\big(\Omega_t - \Omega_{t+1}\big)+ \frac{1}{T\eta_T}\sum_{t=1}^T\big(  \frac{2\alpha_{t+1}^2\delta^2}{b\tau^2\eta_t} +  \frac{2\beta_{t+1}^2\delta^2}{b\tau^2\eta_t}\big) \nonumber \\
 & \leq \frac{1}{T\gamma\eta_T} \big( F(x_1) - F^* + \frac{9L^2_f\gamma}{\lambda\tau}\Delta_1+\frac{2\gamma\delta^2}{b\tau^2\eta_0} \big) + \frac{1}{T\eta_T}\sum_{t=1}^T\big( \frac{2\alpha_{t+1}^2\delta^2}{b\tau^2\eta_t}  + \frac{2\beta_{t+1}^2\delta^2}{b\tau^2\eta_t}\big)  \nonumber \\
 & = \frac{F(x_1) - F^*}{T\gamma\eta_T} + \frac{9L^2_f}{T\eta_T\lambda\tau}\Delta_1 +\frac{2\delta^2}{Tb\tau^2\eta_T\eta_0}
 + \frac{2(c_1^2+c_2^2)\delta^2}{Tb\tau^2\eta_t}\sum_{t=1}^T\eta_t^3 \nonumber \\
 & \leq \frac{F(x_1) - F^*}{T\gamma\eta_T}+ \frac{9L^2_f}{T\eta_T\lambda\tau}\Delta_1 +\frac{2\delta^2}{Tb\tau^2\eta_T\eta_0} + \frac{2(c_1^2+c_2^2)\delta^2}{Tb\tau^2\eta_T}\int^T_1\frac{k^3}{m+t} dt\nonumber \\
 & \leq \frac{F(x_1) - F^*}{T\gamma\eta_T}+ \frac{9L^2_f}{T\eta_T\lambda\tau}\Delta_1 +\frac{2\delta^2}{Tb\tau^2\eta_T\eta_0}
 + \frac{2(c_1^2+c_2^2)\delta^2 k^3}{Tb\tau^2\eta_T}\ln(m+T) \nonumber \\
 & = \bigg(\frac{F(x_1) - F^*}{T\gamma k} + \frac{9L^2_f}{T k\lambda\tau}\Delta_1+\frac{2m^{1/3}\delta^2}{Tb \tau^2k^2}
  + \frac{2(c_1^2+c_2^2)\delta^2 k^2}{Tb\tau^2}\ln(m+T)\bigg)(m+T)^{1/3},
\end{align}
where the second inequality holds by the above inequality \eqref{eq:K5}.
Let $M'' =  \frac{F(x_1) - F^*}{\gamma k} + \frac{9L^2_f\Delta_1}{ k\lambda\tau} +\frac{2m^{1/3}\delta^2}{b \tau^2k^2} + \frac{2(c_1^2+c_2^2)\delta^2 k^2}{b\tau^2}\ln(m+T)$,
we have
\begin{align}
 \frac{1}{T} \sum_{t=1}^T \mathbb{E}\big[ \frac{L_f^2}{4}\|y^*(x_t)-y_t\|^2 + \frac{1}{4}\|\nabla_xf(x_t,y_t) -v_t\|^2 +\frac{1}{4\gamma^2}\|\tilde{x}_{t+1}-x_t\|^2 \big]  \leq \frac{M''}{T}(m+T)^{1/3}. \nonumber
\end{align}
According to Jensen's inequality, we have
\begin{align}
 &  \frac{1}{T} \sum_{t=1}^T \mathbb{E}\big[ \frac{L_f}{2}\|y^*(x_t)-y_t\| + \frac{1}{2}\|\nabla_xf(x_t,y_t) -v_t\| +\frac{1}{2\gamma}\|\tilde{x}_{t+1}-x_t\| \big] \nonumber \\
 & \leq \big( \frac{3}{T} \sum_{t=1}^T \mathbb{E}\big[ \frac{L_f^2}{4}\|y^*(x_t)-y_t\|^2 + \frac{1}{4}\|\nabla_xf(x_t,y_t) -v_t\|^2 +\frac{1}{4\gamma^2}\|\tilde{x}_{t+1}-x_t\|^2 \big] \big)^{1/2} \nonumber \\
 & \leq \frac{\sqrt{3M''}}{T^{1/2}}(m+T)^{1/6} \leq \frac{\sqrt{3M''}m^{1/6}}{T^{1/2}} + \frac{\sqrt{3M''}}{T^{1/3}},
\end{align}
where the last inequality is due to $(a+b)^{1/6} \leq a^{1/6} + b^{1/6}$ for all $a,b>0$.
Thus we obtain
\begin{align}
 \frac{1}{T} \sum_{t=1}^T \mathbb{E}\big[ L_f\|y^*(x_t)-y_t\| + \|\nabla_xf(x_t,y_t) -v_t\| +\frac{1}{\gamma}\|\tilde{x}_{t+1}-x_t\| \big] \leq \frac{2\sqrt{3M''}m^{1/6}}{T^{1/2}} + \frac{2\sqrt{3M''}}{T^{1/3}}. \nonumber
\end{align}
Then by using the above inequality \eqref{eq:MH}, we have
\begin{align}
 \frac{1}{T} \sum_{t=1}^T \mathbb{E}  \|G_{\mathcal{X}}(x_t,\nabla F(x_t),\gamma)\| & \! \leq \frac{1}{T} \sum_{t=1}^T \mathbb{E}\big[L_f\|y^*(x_t)-y_t\| \!+\! \|\nabla_xf(x_t,y_t) -v_t\| \!+\!\frac{1}{\gamma}\|\tilde{x}_{t+1}-x_t\| \big] \nonumber \\
 & \! \leq \frac{2\sqrt{3M''}m^{1/6}}{T^{1/2}} + \frac{2\sqrt{3M''}}{T^{1/3}}.
\end{align}

\end{proof}

\subsection{ Convergence Analysis of Acc-MDA Algorithm for Unconstrained Minimax Optimization }
\label{Appendix:A6}
In this subsection, we study the convergence properties of our  Acc-MDA algorithm for solving the \textbf{unconstrained}
minimax problem \eqref{eq:2}, i.e., $\mathcal{X} = \mathbb{R}^{d_1}$
and $\mathcal{Y} = \mathbb{R}^{d_2}$ (or $\mathcal{Y} \subset \mathbb{R}^{d_2}$). The following convergence analysis builds on the common convergence metric
$\mathbb{E}\|\nabla F(x_t)\|$ used in \citep{lin2019gradient},
where $F(x)=\max_{y\in \mathcal{Y}}f(x,y)$.

\begin{lemma} \label{lem:D03}
Suppose the sequence $\{x_t,y_t\}_{t=1}^T$ be generated from Algorithm \ref{alg:4}.
When $\mathcal{X}=\mathbb{R}^{d_1}$, given $0<\gamma\leq \frac{1}{2\eta_t L_g}$,
we have
\begin{align}
  F(x_{t+1}) &\leq F(x_t)
  + \eta_t\gamma L_f^2 \|y_t - y^*(x_t)\|^2  + \gamma\eta_t\|\nabla_x f(x_t,y_t)-v_t\|^2  \nonumber \\
  & \quad - \frac{\gamma\eta_t}{2}\|\nabla F(x_t)\|^2 - \frac{\gamma\eta_t}{4}\|v_t\|^2.
\end{align}
\end{lemma}

\begin{proof}
This proof is similar to the proof of Lemma \ref{lem:D01}.
 According to Lemma \ref{lem:1}, the approximated function $F(x)$ has $L_g$-Lipschitz continuous gradient.
 Then we have
 \begin{align}
  F(x_{t+1}) &\leq F(x_t) - \gamma\eta_t\langle\nabla F(x_t),v_t\rangle + \frac{\gamma^2\eta_t^2L_g}{2}\|v_t\|^2   \\
  & = F(x_t) + \frac{\gamma\eta_t}{2}\|\nabla F(x_t)-v_t\|^2 - \frac{\gamma\eta_t}{2}\|\nabla F(x_t)\|^2
  + (\frac{\gamma^2\eta_t^2 L_g}{2}-\frac{\gamma\eta_t}{2})\|v_t\|^2 \nonumber \\
  & = F(x_t) + \frac{\gamma\eta_t}{2}\|\nabla F(x_t)-\nabla_x f(x_t,y_t) + \nabla_x f(x_t,y_t)-v_t\|^2
  - \frac{\gamma\eta_t}{2}\|\nabla F(x_t)\|^2 \nonumber \\
  & \quad + (\frac{\gamma^2\eta_t^2 L_g}{2}-\frac{\gamma\eta_t}{2})\|v_t\|^2 \nonumber \\
  & \leq F(x_t) + \gamma\eta_t\|\nabla F(x_t)-\nabla_x f(x_t,y_t)\|^2 + \gamma\eta_t\|\nabla_x f(x_t,y_t)-v_t\|^2
  \nonumber \\
  & \quad - \frac{\gamma\eta_t}{2}\|\nabla F(x_t)\|^2  + (\frac{\gamma^2\eta_t^2 L_g}{2}-\frac{\gamma\eta_t}{2})\|v_t\|^2 \nonumber \\
  & \leq  F(x_t) + \gamma\eta_t\|\nabla F(x_t)-\nabla_x f(x_t,y_t)\|^2 + \gamma\eta_t\|\nabla_x f(x_t,y_t)-v_t\|^2 \nonumber \\
  & \quad - \frac{\gamma\eta_t}{2}\|\nabla F(x_t)\|^2 - \frac{\gamma\eta_t}{4}\|v_t\|^2, \nonumber
 \end{align}
 where the last inequality is due to $0< \gamma \leq \frac{1}{2\eta_t L}$.

 Considering an upper bound of $\|\nabla F(x_t)-\nabla_x f(x_t,y_t)\|^2$, we have
 \begin{align}
   \|\nabla F(x_t)-\nabla_x f(x_t,y_t)\|^2
   = \|\nabla_x f(x_t,y^*(x_t)) - \nabla_x f(x_t,y_t)\|^2  \leq L_f^2 \|y_t - y^*(x_t)\|^2,
 \end{align}
 the last inequality holds by Assumption \ref{ass:5}.
 Then we have
 \begin{align}
  F(x_{t+1}) & \leq F(x_t) +  L_f^2 \|y_t - y^*(x_t)\|^2+ \gamma\eta_t\|\nabla_x f(x_t,y_t)-v_t\|^2\nonumber \\
  & \quad - \frac{\gamma\eta_t}{2}\|\nabla F(x_t)\|^2 - \frac{\gamma\eta_t}{4}\|v_t\|^2.
 \end{align}

\end{proof}

\begin{lemma} \label{lem:E03}
Suppose the sequence $\{x_t,y_t\}_{t=1}^T$ be generated from Algorithm \ref{alg:4}. Under the above assumptions, and set $0< \eta_t\leq 1$
and $\lambda\leq \frac{1}{6L_f}$, we have
\begin{align}
     \|y_{t+1} - y^*(x_{t+1})\|^2 &\leq (1-\frac{\eta_t\tau\lambda}{4})\|y_t -y^*(x_t)\|^2 -\frac{3\eta_t}{4} \|\tilde{y}_{t+1}-y_t\|^2 \nonumber \\
     & \quad + \frac{25\eta_t\lambda}{6\tau}  \|\nabla_y f(x_t,y_t)-w_t\|^2 +  \frac{25\kappa_y^2\gamma^2\eta_t}{6\tau\lambda}\|v_t\|^2,
\end{align}
where $\kappa_y = L_f/\tau$.
\end{lemma}
\begin{proof}
 This proof is the same to the proof of Lemma \ref{lem:E1}.
\end{proof}

\begin{lemma} \label{lem:F03}
 Suppose the stochastic gradients $\{v_t,w_t\}_{t=1}^T$ be generated from Algorithm \ref{alg:4}, we have
\begin{align}
 \mathbb{E} \|\nabla_x f(x_{t+1},y_{t+1}) - v_{t+1}\|^2 & \leq (1-\alpha_{t+1})^2 \mathbb{E} \|\nabla_x f(x_t,y_t) -v_t\|^2
 + \frac{2\alpha_{t+1}^2\delta^2}{b} \nonumber \\
 & \quad + \frac{2(1-\alpha_{t+1})^2L^2_f\eta^2_t}{b}\big(\gamma^2\mathbb{E}\|v_t\|^2 + \mathbb{E}\|\tilde{y}_{t+1}-y_t\|^2\big).
\end{align}
\begin{align}
 \mathbb{E} \|\nabla_y f(x_{t+1},y_{t+1}) - w_{t+1}\|^2 & \leq (1-\beta_{t+1})^2 \mathbb{E} \|\nabla_y f(x_t,y_t) -w_t\|^2 + \frac{2\beta_{t+1}^2\delta^2}{b}
 \nonumber \\
 & \quad + \frac{2(1-\beta_{t+1})^2L^2_f\eta^2_t}{b}\big(\gamma^2\mathbb{E} \|v_t\|^2 + \mathbb{E} \|\tilde{y}_{t+1}-y_t\|^2\big).
\end{align}
\end{lemma}
\begin{proof}
 This proof is the same to the proof of Lemma \ref{lem:F3}.
\end{proof}

\begin{theorem} \label{th:A04}
 (Restatement of Theorem 12)
Suppose the sequence $\{x_t,y_t\}_{t=1}^T$ be generated from Algorithm \ref{alg:4}. When $\mathcal{X}=\mathbb{R}^{d_1}$,
and let $\eta_t = \frac{k}{(m+t)^{1/3}}$
for all $t\geq 0$, $c_1 \geq \frac{2}{3k^3} + \frac{9\tau^2}{4}$
and $c_2 \geq \frac{2}{3k^3} + \frac{75L^2_f}{2}$, $k>0$, $m\geq \max\big(2, k^3, (c_1k)^3, (c_2k)^3\big)$,  $0<\lambda\leq \min\big(\frac{1}{6L_f},
\frac{27b\tau}{16}\big)$ and
$0< \gamma \leq \min\big( \frac{\lambda\tau}{2L_f}\sqrt{\frac{2b}{8\lambda^2 + 75\kappa_y^2b}}, \frac{m^{1/3}}{2L_gk}\big)$,
we have
\begin{align}
 \frac{1}{T} \sum_{t=1}^T \mathbb{E}\|\nabla F(x_t)\|
  \leq \frac{\sqrt{2M''}m^{1/6}}{T^{1/2}} + \frac{\sqrt{2M''}}{T^{1/3}},
\end{align}
where $\Delta_1=\|y_1-y^*(x_1)\|^2$ and $M'' =  \frac{F(x_1) - F^*}{\gamma k} + \frac{9L^2_f\Delta_1}{k\lambda\tau} +\frac{2m^{1/3}\delta^2}{b\tau^2k^2} + \frac{2(c_1^2+c_2^2)\delta^2 k^2}{b\tau^2}\ln(m+T)$.
\end{theorem}
\begin{proof}
This proof is the similar to the proof of Theorem \ref{th:A4}.
As in the proof of Theorem \ref{th:A4},
 let $c_1 \geq \frac{2}{3k^3} + \frac{9\tau^2}{4}$, we have
 \begin{align} \label{eq:K01}
  & \frac{1}{\eta_t}\mathbb{E} \|\nabla_x f(x_{t+1},y_{t+1}) - v_{t+1}\|^2 - \frac{1}{\eta_{t-1}}\mathbb{E} \|\nabla_x f(x_t,y_t) - v_t\|^2  \\
  & \leq -\frac{9}{4}\eta_t\mathbb{E} \|\nabla_x f(x_t,y_t) -v_t\|^2 + \frac{2L^2_f\eta_{t}}{b}\mathbb{E}\big(\gamma^2\|v_t\|^2 + \|\tilde{y}_{t+1}-y_t\|^2\big) + \frac{2\alpha_{t+1}^2\delta^2}{b\eta_t}. \nonumber
 \end{align}
 Let $c_2 \geq \frac{2}{3k^3} + \frac{75L^2_f}{2}$, we have
  \begin{align} \label{eq:K02}
  & \frac{1}{\eta_t}\mathbb{E} \|\nabla_y f(x_{t+1},y_{t+1}) - w_{t+1}\|^2
  -  \frac{1}{\eta_{t-1}}\mathbb{E} \|\nabla_y f(x_t,y_t) - w_t\|^2  \\
  & \leq - \frac{75L^2_f}{2\tau^2}\eta_t\mathbb{E} \|\nabla_y f(x_t,y_t) -w_t\|^2 + \frac{2L^2_f\eta_{t}}{b}\mathbb{E}\big(\gamma^2\|v_t\|^2 + \|\tilde{y}_{t+1}-y_t\|^2\big) + \frac{2\beta_{t+1}^2\delta^2}{b\eta_t}. \nonumber
 \end{align}
 According to Lemma \ref{lem:E03}, we have
 \begin{align} \label{eq:K03}
  \|y_{t+1} - y^*(x_{t+1})\|^2 - \|y_t -y^*(x_t)\|^2 & \leq -\frac{\eta_t\tau\lambda}{4}\|y_t -y^*(x_t)\|^2
  -\frac{3\eta_t}{4} \|\tilde{y}_{t+1}-y_t\|^2 \nonumber \\
     & \quad + \frac{25\eta_t\lambda}{6\tau}  \|\nabla_y f(x_t,y_t)-w_t\|^2 +  \frac{25\kappa_y^2\gamma^2\eta_t}{6\tau\lambda}\|v_t\|^2.
  \end{align}

 At the same time, we give the \emph{Lyapunov} function $\Omega_t$ defined in the proof of the Theorem \ref{th:A4}, 
 \begin{align}
 \Omega_t & = \mathbb{E}\big[ F(x_t) +  \frac{9L^2_f\gamma}{\lambda\tau}\|y_t-y^*(x_t)\|^2 +  \frac{\gamma}{\tau^2\eta_{t-1}}\|\nabla_x f(x_t,y_t)-v_t\|^2 \nonumber \\
 & \quad + \frac{\gamma}{\tau^2\eta_{t-1}}\|\nabla_y f(x_t,y_t)-w_t\|^2 \big]. \nonumber
 \end{align}
 By using Lemma \ref{lem:D03}, we have
 \begin{align}
 & \Omega_{t+1} - \Omega_t \nonumber \\
 &\! = \mathbb{E}\big[F(x_{t+1}) - F(x_t)\big] + \frac{ 9L^2_f\gamma}{\lambda\tau}
 \big( \mathbb{E}\|y_{t+1}-y^*(x_{t+1})\|^2 - \mathbb{E}\|y_t-y^*(x_t)\|^2 \big)
  \nonumber \\
 & \ + \frac{\gamma}{\tau^2} \big(\frac{1}{\eta_t}\mathbb{E}\|\nabla_x f(x_{t+1},y_{t+1})-v_{t+1}\|^2- \frac{1}{\eta_{t-1}}\mathbb{E}\|\nabla_x f(x_t,y_t)-v_t\|^2  \nonumber \\
 &\ +  \frac{1}{\eta_t}\mathbb{E}\|\nabla_y f(x_{t+1},y_{t+1})-w_{t+1}\|^2- \frac{1}{\eta_{t-1}}\mathbb{E}\|\nabla_y f(x_t,y_t)-w_t\|^2 \big) \nonumber \\
 & \! \leq \eta_t\gamma L_f^2 \mathbb{E}\|y_t - y^*(x_t)\|^2  + \gamma\eta_t\mathbb{E}\|\nabla_x f(x_t,y_t)-v_t\|^2  - \frac{\gamma\eta_t}{2}\mathbb{E}\|\nabla F(x_t)\|^2 - \frac{\gamma\eta_t}{4}\mathbb{E}\|v_t\|^2 \nonumber \\
 & \ +\! \frac{9L^2_f\gamma}{\lambda\tau} \big( -\frac{\eta_t\tau\lambda}{4}\mathbb{E}\|y_t -y^*(x_t)\|^2 -\frac{3\eta_t}{4} \mathbb{E}\|\tilde{y}_{t+1}-y_t\|^2 + \frac{25\eta_t\lambda}{6\tau} \mathbb{E}\|\nabla_y f(x_t,y_t)-w_t\|^2  \nonumber \\
 & \ +\! \frac{25\kappa_y^2\gamma^2\eta_t}{6\tau\lambda}\mathbb{E}\|v_t\|^2 \big) -\frac{9\gamma\eta_t}{4} \mathbb{E} \|\nabla_x f(x_t,y_t) -v_t\|^2 \!+\! \frac{2L^2_f\eta_{t}\gamma}{b\tau^2}\big(\gamma^2\mathbb{E}\|v_t\|^2 + \mathbb{E}\|\tilde{y}_{t+1}-y_t\|^2\big) \nonumber \\
 & \ +\! \frac{2\alpha_{t+1}^2\delta^2\gamma}{b\tau^2\eta_t}\!-\!\frac{75L^2_f\gamma}{2\tau^2} \eta_t\mathbb{E}\|\nabla_y f(x_t,y_t) \!-\!w_t\|^2 \!+\! \frac{2L^2_f\eta_{t}\gamma}{b\tau^2}\big(\gamma^2\mathbb{E}\|v_t\|^2 \!+\! \mathbb{E}\|\tilde{y}_{t+1}\!-\!y_t\|^2\big) \!+\! \frac{2\beta_{t+1}^2\delta^2\gamma}{b\tau^2\eta_t} \nonumber \\
 & \!\leq -\frac{5L_f^2\eta_t\gamma }{4}\mathbb{E}\|y^*(x_t)-y_t\|^2 - \frac{5\gamma\eta_t}{4}\mathbb{E}\|\nabla_xf(x_t,y_t) -v_t\|^2 - \frac{\gamma\eta_t}{2}\mathbb{E}\|\nabla F(x_t)\|^2 + \frac{2\alpha_{t+1}^2\delta^2\gamma}{b\tau^2\eta_t} \nonumber \\
 & \ + \frac{2\beta_{t+1}^2\delta^2\gamma}{b\tau^2\eta_t}- \big( \frac{27L^2_f\gamma}{4\lambda\tau} - \frac{4L^2_f\gamma}{b\tau^2}\big)\eta_t \mathbb{E}\|\tilde{y}_{t+1}-y_t\|^2-\big(\frac{\gamma}{4} -\frac{4L^2_f\gamma^3}{b\tau^2} -\frac{75L^2_f\kappa_y^2\gamma^3}{2\lambda^2\tau^2}\big)\eta_t\mathbb{E}\|v_t\|^2 \nonumber \\
 & \! \leq - \frac{\gamma\eta_t}{2}\mathbb{E}\|\nabla F(x_t)\|^2 + \frac{2\alpha_{t+1}^2\delta^2\gamma}{b\tau^2\eta_t} + \frac{2\beta_{t+1}^2\delta^2\gamma}{b\tau^2\eta_t},
 \end{align}
 where the first inequality holds by combining the above inequalities \eqref{eq:K01}, \eqref{eq:K02} and \eqref{eq:K03},
 and the last inequality is due to $0< \gamma \leq \frac{\lambda\tau}{2L_f}\sqrt{\frac{2b}{8\lambda^2 + 75\kappa_y^2b}}$
 and $\lambda\leq \frac{27b\tau}{16}$.
 Then we have
 \begin{align} \label{eq:K04}
  \frac{\eta_t}{2}\mathbb{E}\|\nabla F(x_t)\|^2
 \leq \frac{\Omega_t - \Omega_{t+1}}{\gamma}+ \frac{2\alpha_{t+1}^2\delta^2}{b\tau^2\eta_t} + \frac{2\beta_{t+1}^2\delta^2}{b\tau^2\eta_t}.
 \end{align}
Let $\Delta_1=\|y_1-y^*(x_1)\|^2$, we have
\begin{align} \label{eq:K05}
 \Omega_1 &= F(x_1)+ \frac{9L^2_f\gamma}{\lambda\tau}\|y_1-y^*(x_1)\|^2+\frac{\gamma}{\tau^2\eta_0}\mathbb{E}\|\nabla_x f(x_1,y_1)-v_1\|^2 + \frac{\gamma}{\tau^2\eta_{0}}\mathbb{E}\|\nabla_y f(x_1,y_1)-w_1\|^2 \nonumber \\
 & = F(x_1) +  \frac{9L^2_f\gamma}{\lambda\tau}\|y_1-y^*(x_1)\|^2 +\frac{\gamma}{\tau^2\eta_0}\mathbb{E}\|\nabla_x f(x_1,y_1)- \frac{1}{b}\sum_{i=1}^b\hat{\nabla}_x f(x_1,y_1;\xi^1_i)\|^2 \nonumber \\
 & \quad + \frac{\gamma}{\tau^2\eta_0}\mathbb{E}\|\nabla_y f(x_1,y_1)- \frac{1}{b}\sum_{i=1}^b\hat{\nabla}_y f(x_1,y_1;\xi^1_i)\|^2 \nonumber \\
 & \leq F(x_1)+  \frac{9L^2_f\gamma}{\lambda\tau}\Delta_1 +\frac{2\gamma\delta^2}{b\tau^2\eta_0},
\end{align}
where the last inequality holds by Assumption \ref{ass:2}.

Taking average over $t=1,2,\cdots,T$ on both sides of \eqref{eq:K04} and due to $\eta_T^{-1} \geq \eta_t^{-1}$
for any $0\leq t\leq T$,
we have
 \begin{align}
 & \frac{1}{T} \sum_{t=1}^T \frac{1}{2}\mathbb{E}\|\nabla F(x_t)\|^2 \nonumber \\
 & \leq  \frac{1}{T\gamma\eta_T} \sum_{t=1}^T\big(\Omega_t - \Omega_{t+1}\big)+ \frac{1}{T\eta_T}\sum_{t=1}^T\big(  \frac{2\alpha_{t+1}^2\delta^2}{b\tau^2\eta_t} +  \frac{2\beta_{t+1}^2\delta^2}{b\tau^2\eta_t}\big) \nonumber \\
 & \leq \frac{1}{T\gamma\eta_T} \big( F(x_1) - F^*+  \frac{9L^2_f\gamma}{\lambda\tau}\Delta_1 +\frac{2\gamma\delta^2}{b\tau^2\eta_0} \big) + \frac{1}{T\eta_T}\sum_{t=1}^T\big( \frac{2\alpha_{t+1}^2\delta^2}{b\tau^2\eta_t}  + \frac{2\beta_{t+1}^2\delta^2}{b\tau^2\eta_t}\big)  \nonumber \\
 & = \frac{F(x_1) - F^*}{T\gamma\eta_T} +  \frac{9L^2_f}{T\eta_T\lambda\tau}\Delta_1 +\frac{2\delta^2}{Tb\tau^2\eta_T\eta_0}
 + \frac{2(c_1^2+c_2^2)\delta^2}{Tb\tau^2\eta_t}\sum_{t=1}^T\eta_t^3 \nonumber \\
 & \leq \frac{F(x_1) - F^*}{T\gamma\eta_T} +  \frac{9L^2_f}{T\eta_T\lambda\tau}\Delta_1 +\frac{2\delta^2}{Tb\tau^2\eta_T\eta_0} + \frac{2(c_1^2+c_2^2)\delta^2}{Tb\tau^2\eta_T}\int^T_1\frac{k^3}{m+t} dt\nonumber \\
 & \leq \frac{F(x_1) - F^*}{T\gamma\eta_T}+  \frac{9L^2_f}{T\eta_T\lambda\tau}\Delta_1 +\frac{2\delta^2}{Tb\tau^2\eta_T\eta_0}
 + \frac{2(c_1^2+c_2^2)\delta^2 k^3}{Tb\tau^2\eta_T}\ln(m+T) \nonumber \\
 & = \bigg(\frac{F(x_1) - F^*}{T\gamma k} + \frac{9L^2_f}{T k \lambda\tau}\Delta_1+\frac{2m^{1/3}\delta^2}{Tb\tau^2k^2}
  + \frac{2(c_1^2+c_2^2)\delta^2 k^2}{Tb\tau^2}\ln(m+T)\bigg)(m+T)^{1/3},
\end{align}
where the second inequality holds by the above inequality \eqref{eq:K05}.
Let $M'' =  \frac{F(x_1) - F^*}{\gamma k} + \frac{9L^2_f\Delta_1}{k \lambda\tau}+\frac{2m^{1/3}\delta^2}{b\tau^2k^2} + \frac{2(c_1^2+c_2^2)\delta^2 k^2}{b\tau^2}\ln(m+T)$,
we have
\begin{align}
 \frac{1}{T} \sum_{t=1}^T \mathbb{E}\|\nabla F(x_t)\|^2  \leq \frac{2M''}{T}(m+T)^{1/3}. \nonumber
\end{align}
According to Jensen's inequality, we have
\begin{align}
 \frac{1}{T} \sum_{t=1}^T \mathbb{E}\|\nabla F(x_t)\|
 & \leq \big( \frac{1}{T} \sum_{t=1}^T\mathbb{E}\|\nabla F(x_t)\|^2 \big)^{1/2} \nonumber \\
 & \leq \frac{\sqrt{2M''}}{T^{1/2}}(m+T)^{1/6} \leq \frac{\sqrt{2M''}m^{1/6}}{T^{1/2}} + \frac{\sqrt{2M''}}{T^{1/3}},
\end{align}
where the last inequality is due to $(a+b)^{1/6} \leq a^{1/6} + b^{1/6}$ for all $a,b>0$.

\end{proof}

\section{ Comparison of Assumptions Used in Zeroth-Order Methods }
\label{appendix:B}

We admit that our methods (Acc-ZOM, Acc-ZOMDA, Acc-MDA) and the existing \textbf{variance-reduced} zeroth-order and
first-order methods (e.g., ZO-SPIDER-Coord, SPIDER-SZO, ZO-SREDA-Boost, SREDA and SREDA-boost) rely on a relative strong assumption \\ (\textbf{component function smoothness}), i.e., $\|\nabla f(x_1;\xi) - \nabla f(x_2;\xi)\| \leq L\|x_1-x_2\|$ for mini-optimization and $\|\nabla f(x_1,y_1;\xi) - \nabla f(x_2,y_2;\xi)\| \leq L(\|x_1-x_2\|+\|y_1-y_2\|)$ for minimax-optimization.

At the same time, we also argue that the comparison \textbf{non-variance-reduced} methods (such as ZO-AdaMM and ZO-Min-Max) in Table \ref{tab:1}
require stronger assumptions than the \textbf{component function smoothness} assumption. For example, ZO-AdaMM \citep{chen2019zo} method requires the following
two assumptions (Please see the page 4 of paper ``ZO-AdaMM: Zeroth-Order Adaptive Momentum
Method for Black-Box Optimization" https://arxiv.org/pdf/1910.06513.pdf):
\begin{itemize}
 \item[A1)] $f_t(\cdot)=f(\cdot,\xi_t)$ has $L_g$-Lipschitz continuous gradient, where $L_g>0$.
 \item[A2)] $f_t$ has $\eta$-bounded stochastic gradient $\|\nabla f_t(x)\|_{\infty} \leq \eta$.
\end{itemize}
In fact, the above assumption A1 is a component function smoothness assumption.
Clearly, the above assumptions A1 and A2 required in ZO-AdaMM method is more stronger than the component function smoothness assumption required in
our methods.

Meanwhile, ZO-Min-Max \citep{liu2019min} method requires a stronger \textbf{bounded
gradient} Assumption ( Please see Assumption A1 at the page 4 of paper
'Min-Max Optimization without Gradients: Convergence and Applications to
Black-Box Evasion and Poisoning Attacks' https://arxiv.org/pdf/1909.13806.pdf):
\begin{itemize}
 \item[A1)] $f(x,y)=\mathbb{E}_{\xi\sim p}[f(x,y;\xi)]$ has bounded
gradients $\|\nabla_x f(x,y;\xi)\| \leq \eta^2$ and $\|\nabla_y f(x,y;\xi)\|  \\ \leq \eta^2$ for stochastic optimization with $\xi\sim p$.
\end{itemize}
Clearly, this Assumption required in ZO-Min-Max method is stronger than the component function smoothness assumption required in
our methods.

\section{Query Complexity of ZO-Min-Max Method in \citep{liu2019min}}
\label{appendix:C}

\cite{liu2019min} do not provide the explicit query complexity of ZO-Min-Max method. 
However, the query complexity $O((d_1+d_2)\epsilon^{-6}))$ of ZO-Min-Max method given in \citep{wang2020zeroth} is incorrect (See Table 1 at page 10 of https://arxiv.org/pdf/2001.07819.pdf). Here, we give a correct complexity $O((d_1+d_2)\kappa^6_y\epsilon^{-6}))$ of ZO-Min-Max method based on the results in the original paper \citep{liu2019min}. The detailed proof is given as follows:

From Theorems 1-2 and Remarks 1-2 in \citep{liu2019min} ( Please see the pages 5-6 of paper:
``Min-Max Optimization without Gradients: Convergence and Applications to
Black-Box Evasion and Poisoning Attacks" https://arxiv.org/pdf/1909.13806.pdf ), we have  $\beta=\frac{\gamma}{8L^2_y}$, $\alpha = 1/(L_x + \frac{4L^2_x}{\gamma^2\beta} + \beta L^2_x)$, $\zeta=\min\big(\frac{2L^2_y}{\gamma}, \frac{2L^2_x}{\gamma}+\frac{L_x}{2}\big)$
and $c=\max\big(L_x+3/\alpha, 3/\beta\big)$, where $L_x$ and $L_y$ are the smooth parameters, $\gamma$ is the parameter about strongly concave $f(x,y)$ w.r.t. $y$.

For notational simplicity, let $L=L_x=L_y$ as in \citep{luo2020stochastic,xu2020enhanced} and $\kappa_y = L/\gamma$. It is easy verified that $\beta^{-1}=O(\kappa_y)$ and $\alpha^{-1}=O(\kappa_y^3)$, $c=O(\kappa_y^3)$ and $\zeta=O(\kappa_y)$. Thus, we have $\frac{c}{\zeta} = O(\kappa_y^2)$ in Theorem 1. Since Theorem 2 is similar to Theorem 1 in \citep{liu2019min}, we also have $\frac{c}{\zeta'} = O(\kappa_y^2)$. Then based on the remarks about Theorems 1-2 in \citep{liu2019min},
we have $\mathbb{E}\|\mathcal{G}(x^r,y^r)\|^2=O(\frac{\kappa_y^2}{T}+\frac{\kappa_y^2}{b}+\frac{\kappa_y^2\tilde{d}}{q})$, where $(x^r, y^r)$ randomly picked
from $\{(x^t,y^t)\}_{t=1}^T$, and  $\tilde{d}=d_1+d_2$,
$b$ is mini-batch size, and $q$ is the number of
random direction vectors for estimating zeroth-order gradient.

Considering $\mathbb{E}\|\mathcal{G}(x^r,y^r)\|=O(\frac{\kappa_y}{\sqrt{T}}+\frac{\kappa_y}{\sqrt{b}}+\frac{\kappa_y\sqrt{\tilde{d}}}{\sqrt{q}}) \leq \epsilon$, let $T=b=q/\tilde{d}$, then we have $T=b=q/\tilde{d} = O(\kappa^2_y\epsilon^{-2})$.
Since the ZO-Min-Max algorithm requires query $4bq$ function values to estimate zeroth-order gradients  $\hat{\nabla}_xf(x,y)$ and
$\hat{\nabla}_yf(x,y)$ at each iteration, and need $T$ iterations, it requires a query complexity of
$4bqT=O\big(\tilde{d}\kappa_y^6\epsilon^{-6}\big) = O\big((d_1+d_2)\kappa_y^6\epsilon^{-6}\big)$
for finding an $\epsilon$-stationary point (i.e., $\mathbb{E}\|\mathcal{G}(x^r,y^r)\| \leq \epsilon$).
At the same time, the  mini-batch size is  $\max(b,q) = O((d_1+d_2)\kappa_y^2\epsilon^{-2})$ in ZO-Min-Max method.


\vskip 0.2in
\bibliography{AccZOM}

\end{document}